\documentclass[reqno]{amsart}
\usepackage{amssymb,amsmath,amsfonts,latexsym}
\usepackage{pb-diagram}
\numberwithin{equation}{section}
\newcommand{\A}{{\mathcal{A}}}

\newcommand{\tr}{{\operatorname{tr}}}
\newcommand{\orbit}{{\mathcal O}}

\newcommand{\Hom}{{\operatorname{Hom}}}
\newcommand{\tot}{{\operatorname{Tot}}}

\newcommand{\id}{{\operatorname{id}}}
\newcommand{\topotimes}{\overline{\otimes}\,}

\newcommand{\complex}{{\mathbb C}}

\newcommand{\C}{{\mathbb{C}}}
\newcommand{\R}{{\mathbb{R}}}
\newcommand{\Z}{{\mathbb{Z}}}

\newcommand{\N}{{\mathbb{N}}}
\newcommand{\h}{{\mathbb{H}}}

\newcommand{\calc}{{\mathcal C}}

\newcommand{\fraka}{\mathfrak{a}}
\newcommand{\frakb}{\mathfrak{b}}

\theoremstyle{plain}
        \newtheorem{theorem}{Theorem}[section]
        \newtheorem{lemma}[theorem]{Lemma}
        \newtheorem{proposition}[theorem]{Proposition}
        \newtheorem{corollary}[theorem]{Corollary}

\theoremstyle{definition}
        \newtheorem{definition}[theorem]{Definition}
        \newtheorem{remark}[theorem]{Remark}
        \newtheorem{example}[theorem]{Example}

\newenvironment{sketchproof}{{\noindent \it Sketch of Proof. }}{
        \mbox{ }\hfill$\Box$\vspace{1.5ex}\par}
\newenvironment{lemproof}{{\noindent \it Proof of the Lemma. }}{
        \mbox{ }\hfill$\Box$\vspace{1.5ex}\par}

\date{\today}
\title{Homology of formal deformations of proper \'etale Lie groupoids}
\author{N.~Neumaier, M.J.~Pflaum, H.B.~Posthuma \textrm{and} X.~Tang}
\begin{document}
\mbox{ } \vspace{-15mm}\\
\begin{abstract}
  In this article, the cyclic homology theory of formal deformation 
  quantizations of the convolution algebra associated to a proper \'etale 
  Lie groupoid is studied.
  We compute the Hochschild cohomology of the convolution algebra and express
  it in terms of alternating multi-vector fields on the associated inertia 
  groupoid.
  We introduce a noncommutative Poisson homology whose computation enables 
  us to determine the Hochschild homology of formal deformations of the 
  convolution algebra.
  Then it is shown that the cyclic (co)homology of  such formal deformations
  can be described by an appropriate sheaf cohomology theory. This enables
  us to determine the corresponding cyclic homology groups in terms of
  orbifold cohomology of the underlying orbifold.  
  Using the thus obtained description of cyclic cohomology of the deformed 
  convolution algebra, we give a complete classification of all traces on 
  this formal deformation, and provide an explicit construction.
\end{abstract}
\address{\newline
   Nikolai Neumaier, {\tt nikolai.neumaier@physik.uni-freiburg.de }\newline
   \indent {\rm Physikalisches Institut, Universit\"at Freiburg, 
           Germany } \newline
   Markus J. Pflaum, {\tt pflaum@math.uni-frankfurt.de}\newline
   \indent {\rm Fachbereich Mathematik, Goethe-Universit\"at Frankfurt/Main, 
           Germany } \newline
   Hessel Posthuma, {\tt posthuma@math.uni-frankfurt.de}\newline
   \indent {\rm Fachbereich Mathematik, Goethe-Universit\"at Frankfurt/Main, 
           Germany } \newline
   Xiang Tang, {\tt xtang@math.udavis.edu}   \newline   
   \indent {\rm  Department of Mathematics, University of California, Davis, 
           USA } 
}
\maketitle
\tableofcontents
\newpage
\section{Introduction}
In symplectic geometry and mathematical physics one often encounters, 
for example by reduction, Poisson spaces which are singular. One of the 
easiest examples of such singular spaces is given by a symplectic orbifold; 
this an orbifold which admits a covering by orbifold charts equipped with an 
invariant symplectic structure. Therefore, the study of deformation 
quantization and index theory of such spaces appears naturally.

To address such questions, one first has to decide what ``algebra 
of smooth functions'' on an orbifold one wants to consider. Any orbifold 
has a natural sheaf of functions which locally can be lifted to smooth 
invariant functions on any orbifold chart. For this algebra on a symplectic
orbifold, a deformation quantization was constructed in \cite{P2}, 
generalizing Fedosov's method \cite{fe:book} on smooth manifolds. 
This case was further studied by Fedosov--Schulze--Tarkhanov in 
\cite{fst:conjecture}, making several interesting conjectures on the 
related index problem.

However from the point of view of noncommutative geometry \cite{connes}, 
an orbifold presents one of the prime examples of a 
``noncommutative manifold'': its ``algebra of smooth functions'' is given 
by the noncommutative convolution algebra on a proper \'{e}tale Lie
groupoid, whose quotient space identifies with the underlying orbifold
(cf.~\cite{moerdijk}). This construction generalizes the crossed product 
of the algebra of smooth functions on a manifold by a finite group.
As shown in \cite{ta:lemma}, a Poisson structure on the orbifold induces 
a natural noncommutative Poisson structure on this algebra in the sense 
of Block--Getzler \cite{bl-ge:quantization} and Xu \cite{xu}, which admits 
a deformation quantization. It is the properties of this deformed algebra 
that we study in this paper. Notice that it contains the deformed algebra 
of \cite{pflaum} as the subalgebra of invariants under the groupoid. 

The first step in understanding the deformation quantization of a groupoid
algebra is to count how many noncommutative Poisson structures it has. 
We partially answer this question by calculating the Hochschild cohomology 
of the groupoid algebra of a proper \'etale groupoid. When the groupoid 
is a manifold, this is given by the Hochschild--Kostant--Rosenberg Theorem.
In the case of an orbifold, partial results have already appeared in the
literature, e.g.~\cite{giaquinto} for the case of a global quotient of an 
algebraic variety by a finite group.
In Section 3, we present a calculation in the case of a proper \'etale Lie
groupoid. By Teleman's localization technique, we relate this cohomology 
to the sheaf cohomology of the multivector fields on the corresponding 
inertia groupoid. This is the first step to classify the Poisson structures
on a groupoid algebra. We leave the study of the Gerstenhaber bracket and 
possible extensions to a ``noncommutative formality theorem'' for future 
research. 

The second step in understanding the quantization of an algebra 
is its semiclassical geometry, or in other words its noncommutative 
Poisson geometry. In this paper, we introduce a noncommutative Poisson 
homology generalizing Brylinski's definition on a Poisson manifold 
\cite{byi:poisson}. This noncommutative Poisson homology appears as the 
$E^1$-term of the spectral sequence associated to the $\hbar$-adic 
filtration in the Hochschild homology of the deformed algebra. We calculate 
this Poisson homology in case of the noncommutative Poisson structure on 
the convolution algebra of a groupoid associated to an orbifold with a Poisson
structure. Our calculation uses the methods developed by Connes, Burghelea, 
Brylinski, Nistor, and Crainic in calculating the cyclic homology of an 
\'etale groupoid. We track the change of the Poisson differential in the 
various steps of the calculation of the Hochschild homology of the 
groupoid algebra, and relate it to the sheaf homology of Brylinski's 
complex on the corresponding inertia groupoid. 

Next, we compute the Hochschild and cyclic (co)homology of the deformation
quantization of the convolution algebra of a symplectic orbifold. Our 
computation draws upon two ideas: one is the spectral sequence introduced 
in \cite{brg} and \cite{nt} associated to the $\hbar$-adic filtration of 
the Hochschild complex and its relation to the noncommutative Poisson 
homology of the previous section. Second is the localization to the 
inertia groupoid, as used in the calculation of the Hochschild and 
cyclic homology of the ``classical'' groupoid algebra of 
Brylinski--Nistor \cite{bn} and Crainic \cite{crainic}. The main difficulty 
here is that, due to the noncommutative nature of the sheaf of 
``quantized functions'', one has to add ``quantum corrections'' to this 
localization map. Interestingly, all our results are given in terms of the 
orbifold cohomology of Chen--Ruan \cite{cr}.

The calculation of the cyclic cohomology gives in particular a complete 
classification of the traces of the deformed groupoid algebra. In the 
last section we give an explicit construction of all traces, building on 
earlier work by Fedosov \cite{fe:g-index}. Finally, our constructions also 
clarify a conjecture in \cite{fst:conjecture} on a certain 
``Picard group'' acting on the space of traces.

Since our computations use quite some machinery of groupoids and cyclic 
homology developed by several people over the years, we have included, 
for the convenience of the reader, a rather detailed section devoted to 
these subjects. Its main purpose is to give an introduction into the 
ideas of Connes, Burghelea, Brylinski, Nistor, and Crainic in calculating 
the cyclic homology of \'etale groupoids, and to set up the notation.

Our paper is related to the recent paper 
\cite{dolget} by Dolgushev--Etingof, where for the quotient of a
smooth complex affine symplectic variety $X$ by a finite group 
the Hochschild cohomology of the convolution algebra and of the algebra of
invariant regular functions is computed. 
Note however, that due to the algebraic nature of their setup the methods 
used there are quite different to the ones used here.

Let us also mention at this point, that some of the results presented here
have been obtained by the fourth author in his PhD-thesis and, 
independently, by the collaboration of the remaining authors. 
During a conference in Luminy, where three of us,
namely M.J.P., H.P.~and X.T.~met, we then decided to continue our work 
together and write a joint paper about our results.

\vspace{3mm}
\paragraph{\bf Acknowledgements}
We would like to thank K.~Behrend and P.~Xu for organizing 
the conference on ``Groupoids and Stacks'' in Luminy in June 2004,
where our collaboration started. X.T.~would like to thank the 
University of Frankfurt for hospitality during his visit there.
He would also like to thank Alan Weinstein, his Ph.D.~advisor, for many 
stimulating discussions, M.~Karoubi and R.~Nest for hosting his visit of 
IHP in summer 2004, and A.~Gorokhovsky, R.~Nest, B.~Tsgyan, and V.~Dolgushev 
for helpful discussions. M.J.P.~and H.P.~gratefully acknowledge
financial support by the Deutsche Forschungsgemeinschaft. 
M.J.P.~would like to thank B.~Fedosov for helpful discussions about
traces on star product algebras over orbifolds.
\section{Preliminaries}
\subsection{Notation}
For clarity, we collect here some of the notation used throughout the paper. 
\begin{itemize}
\item $X$ denotes an orbifold,
\item $\tilde{X}$ its inertia orbifold as constructed in \cite{cr}.
\item By $G$ we denote a smooth \'{e}tale groupoid modeling $X$, 
\item $\A$ is the $G$-sheaf of smooth functions on $G_0$, if not stated
      otherwise,
\item and $\A^\hbar$ is the $G$-sheaf of a formal deformation of $\A$. 
\item Concerning products, $f \cdot g$ denotes the pointwise product 
      of two functions $f,g\in C^\infty(G_1)$,
\item whereas $f * g$ is the convolution product on the groupoid $G$,
\item and, finally, $f\star g$ is the star product associated to 
      a formal deformation.
\end{itemize}

 \subsection{Orbifolds}
 The notion of an orbifold was introduced by Satake \cite{satake}. 
 Roughly  speaking, an orbifold is a
 second countable Hausdorff space $X$ which is locally modeled on the quotient 
 of open subsets of $\R^n$ by a finite group. In this article we will 
 use the language of groupoids to model orbifolds, following the approach by 
 Moerdijk \cite{moerdijk} (see also \cite[Chap.~5]{moerdijk-mrcun} for
 an introduction). To set up notation and for the convenience of the 
 reader let us recall some basic notions from the theory of groupoids.
 
 A groupoid $G$ is a small category in which every morphism is invertible. 
 More explicitly, denote by $G_0$ the space of objects and by $G_1$ the
 space of arrows in $G$. The groupoid structure is encoded by the following 
 five maps:
 \begin{displaymath}
   G_1 \times_{G_0} G_1  \overset{m}{\rightarrow} G_1 \overset{i}{\rightarrow}
   G_1 \overset{s}{\underset{t}{\rightrightarrows}}
   G_0 \overset{u}{\rightarrow} G_1
 \end{displaymath}
 Here, $s$ and $t$ are the source and target map, $m$ is the
 multiplication resp.~composition  
 $(g,h) \mapsto  m(g,h) := gh := g\circ h$, $i$ denotes
 the inverse which is given by $g \mapsto i(g):=g^{-1}$ and finally $u$ is the 
 inclusion of 
 objects by identity arrows, i.e., $u(x):= \id_x$ for all $x\in G_0$. 
 An arrow $g\in G_1$ with $s(g)=x$ and $t(g)=y$ will often be denoted by 
 $g:x\rightarrow y$. Moreover, if no confusion can arise, we write
 $G$ instead of $G_1$.

 A groupoid $G$ is a Lie groupoid (also called a smooth groupoid) if
 both $G_0$ and $G_1$ are smooth manifolds, all the structure maps are
 smooth and $s$ and $t$ submersions. It then follows that $u$ is an
 immersion and that $i$ is a diffeomorphism. A Lie groupoid is called
 \textit{proper} when the map $(s,t):G_1\rightarrow G_0\times G_0$ is
 proper. It is called a foliation groupoid, if every isotropy group $G_x$ is
 discrete. An \textit{\'{e}tale} groupoid is a special type of foliation
 groupoid for which $s$ and $t$ are local diffeomorphisms.
 \begin{definition} {\cite[Def.~3.1]{moerdijk}}
   An orbifold groupoid is a proper foliation groupoid.
 \end{definition}
 More precisely, one proves that the orbit space $X:=G_0/G_1$ has a canonical
 orbifold structure. This description of orbifolds by groupoids goes back to 
 \cite{mp}. The fundamental idea is that the orbifold structure on $X$ in fact 
 only depends on the \textit{Morita equivalence class} of the groupoid $G$. 
 For an introduction to the theory of Morita equivalence of groupoids we
 refer to \cite{moerdijk}. We merely remark that this allows us to choose a 
 proper \'{e}tale groupoid $G$ representing the orbifold $X$. This groupoid
 has the property that for each $x\in G_0$ there exists a neighbourhood 
 $U_x\subset G_0$ such that the restriction $G|_{U_x}$ is isomorphic 
 to a translation groupoid $\Gamma_x\ltimes U_x$, with $\Gamma_x$ a finite 
 group. This property gives the connection with the usual definition of 
 orbifolds in terms of local charts. In the following, we will denote by 
 $\pi:G\rightarrow X$ the projection onto the orbit space of $G$.
 \subsection{Sheaves on orbifolds}
\label{sheaves}
 The theory of sheaves on orbifolds is discussed in \cite{mp}. 
 One can use an orbifold groupoid to model the
 category of sheaves. In general a $G$-sheaf $\mathcal S$ on an 
 \'etale groupoid  $G$ is a sheaf on $G_0$ with a right action of $G$. 
 This means that any arrow $g:x\rightarrow y$ induces a morphism on 
 stalks $\mathcal S_y\rightarrow\mathcal S_x$ satisfying the
 obvious properties. A section $a$ of a sheaf is said to be invariant, if
 on the level of germs one has 
 $[a]_y g=[a]_x$ for every arrow $g:x\rightarrow y$. The abelian category of
 $G$-sheaves of abelian groups is denoted by $\mathsf{Sh}(G)$. 
 There is a left exact functor 
 $\Gamma_{\text{\tiny inv}}:\mathsf{Sh}(G)\rightarrow \mathsf{Ab}$, 
 where $\mathsf{Ab}$ is the category of abelian groups, given by 
 associating to a $G$-sheaf $\mathcal S$ its global invariant sections. 
 Its right derived functors define the groupoid cohomology groups 
 $H^\bullet(G,\mathcal S)$.

 Likewise, we have the compactly supported cohomology groups: consider the
 functor 
 $\Gamma_{\text{\tiny inv},\text{c}}:\mathsf{Sh}(G)\rightarrow \mathsf{Ab}$,
 defined by
 $$
  \Gamma_{\text{\tiny inv},\text{c}}(\mathcal S)=
  \{ a\in\Gamma_{\text{\tiny inv}}(\mathcal S)\mid
  \pi(\operatorname{supp}(a)) \text{ is compact in } X\}.
 $$ 
 Since $G$ is a proper groupoid, this functor is
 left exact and the compactly supported cohomology groups
 $H^\bullet_\text{c}(G,\mathcal S)$ are defined as the right derived functors
 of $\Gamma_{\text{\tiny inv},\text{c}}$. 
 This definition extends in the usual way to define the hypercohomology 
 groups $\h^\bullet(G,\mathcal S^\bullet)$ and 
 $\h_\text{c}^\bullet(G,\mathcal S^\bullet)$ of any cochain complex 
 $\mathcal S^\bullet$ of $G$-sheaves.

\begin{remark}
 The cohomology groups $H^\bullet_\text{c}(G,\mathcal S)$ are isomorphic 
 to the (re-indexed) homology groups of \cite{cm}, see Section 4.9 and 
 4.10 of that paper. 
 In view of the cohomological indexing, we use the invariant 
 instead of the coinvariant sections; the two are isomorphic as one can show 
 by an averaging argument
 (properness is essential for this, cf.~e.g.~\cite[p.~280]{loday}). 
 The homology theory of \cite{cm} is defined for any \'{e}tale groupoid,
 but not a derived functor in general.
\end{remark}

 More generally, one can associate to every morphism 
 $f:G\rightarrow H$ of \'{e}tale groupoids a functor 
 \begin{displaymath}
    f_!: \: \mathsf{Sh}(G)\rightarrow \mathsf{Sh} (H).
 \end{displaymath}
 For its construction, observe first that the functors 
 $f_*$ and $f^{-1}$ can be extended to the category of sheaves on 
 \'etale groupoids in the obvious manner. The sheaf 
 $f_!\mathcal S\in \mathsf{Sh}(H)$ then has stalk at $x\in G_0$ given by
 \begin{equation}
 \label{stalk} 
   (f_!\mathcal S)_x:=H_\text{c}^0(x/f,\pi_x^{-1}\mathcal S),
 \end{equation}
 where $x/f$ denotes the comma groupoid over $x$, that is the fiber of 
 $f$ over $x$:
 \begin{displaymath}
 \begin{diagram}
   \node{x/f}\arrow{e,t}{\pi_x}\arrow{s}\node{H}\arrow{s,r}{f}\\
   \node{1}\arrow{e,t}{x}\node{G.}
 \end{diagram}
 \end{displaymath}
 Analogously, the right derived functors $R^kf_!$ have a similar
 construction using the higher compactly supported cohomology
 groups of the comma groupoids. The Leray spectral sequence
 generalizes to the category of sheaves on \'etale groupoids. In
 particular, in the case of orbifolds this spectral sequence
 degenerates for the projection $\pi:G\rightarrow X$ and induces an
 isomorphism
\begin{equation}
\label{d}
  H^\bullet_\text{c}(G,\mathcal S)\cong H^\bullet_\text{c}(X,\pi_!\mathcal S).
\end{equation}

 Let us briefly discuss the Bar-complexes used in explicit computations 
 of these cohomology groups. Denote by $G^{(k)}$ the space of $k$-composable
 arrows: 
\begin{displaymath}
  G^{(k)}=\{(g_1,\ldots,g_k)\in G^k \mid s(g_i)=t(g_{i+1}),~i=1,\ldots,k-1\}.
\end{displaymath} 
 These spaces are part of a simplicial manifold with face maps 
 $d_i:G^{(k)}\rightarrow G^{({k-1)}}$, $i=0,\ldots,k$ defined as usual by
\begin{displaymath}
\begin{split}
   d_i(g_1,\ldots,g_k)=
   \begin{cases}
     (g_2,\ldots,g_k), & \text{ for $i=0$,} \\
     (g_1,\ldots,g_{i}\cdot g_{i+1},g_{i+2},\ldots,g_k), & 
     \text{ for $1\leq i\leq k-1$,} \\
     (g_1,\ldots,g_{k-1}), &  \text{ for $i=k$.} 
   \end{cases}
\end{split}
\end{displaymath}
Note that $d_0,d_1:G_1\rightarrow G_0$ are simply the source
and target map. The geometric realization of this simplicial space 
is a model for the classifying space $BG$ of the groupoid $G$.

We have two maps $\epsilon_k,\tau_k:G^{(k)}\rightarrow G_0$, which send 
a string
\begin{equation*}
 x_0\stackrel{g_1}{\longleftarrow}x_1\stackrel{g_2}{\longleftarrow}
 \ldots\stackrel{g_k}{\longleftarrow}x_k
\end{equation*}
to $x_k$ resp.~$x_0$. Let $\mathcal S$ be a $G$-sheaf. Define 
$\mathcal S^k:=\tau^{-1}_k\mathcal S$ and put 
$$ B_k(G,\mathcal S):=\Gamma_\text{\rm c}(G^{(k)},\mathcal S^k). $$ 
These vector spaces can be turned into a simplicial space by observing 
that there are isomorphisms $d_i^*\mathcal S^{k-1}\cong \mathcal S^k$ 
which act on the stalks as identity for $i\neq 0$, but by $g_1$, if $i=0$. 
Using this isomorphism, the simplicial maps $d_i$ induce 
differentials in the obvious way, and its associated homology groups 
compute the homology $H_\bullet(G,\mathcal S)$, in case 
$\mathcal S$ is c-soft. 

\begin{definition}
\label{hyperhom}
 (Cf.~\cite[Sec.~3]{cm} and \cite[Sec.~2.3]{crainic}.)
 Let $\mathcal S_\bullet$ be a (bounded below) chain complex of 
 c-soft $G$-sheaves. The {\it hyperhomology} of $\mathcal S_\bullet$ on $G$
 then is defined as the total homology of the double complex
 $B_\bullet (G,\mathcal S_\bullet)$, i.e.,
 $$
   \h_\bullet (G,\mathcal S_\bullet ) := 
   H_\bullet (\operatorname{Tot} (B_\bullet (G,\mathcal S_\bullet)).
 $$ 
\end{definition}

\subsection{Orbifold cohomology}
The notion of twisted sectors of an orbifold plays an important role 
in index theory \cite{kawasaki}. 
Loosely speaking, it is a geometric way of dealing with the
``stacky aspects'' of an orbifold, that is, the automorphisms of
points. As before, we let $X$ be an orbifold, represented by a groupoid
$G$. Denote by $B^{(0)}\subset G_1$ the ``space of loops'' 
\begin{displaymath}
  B^{(0)}=\{g\in G \mid s(g)=t(g)\}.
\end{displaymath}
 The groupoid $G$ acts on
$B^{(0)}$ by conjugation and one defines 
$$\Lambda G :=B^{(0)}\rtimes G.$$
This groupoid has ``loops'' in $G$ as objects, and the space of arrows
can be identified with
$$\Lambda G_1=\big\{(g_1,g_2)\in G^{(2)} \mid g_1 \in B^{(0)}\big\}.$$
One observes that $\Lambda G$ is again an orbifold groupoid which comes
equipped with a canonical morphism $\beta:\Lambda G\rightarrow G$. The
orbifold underlying $\Lambda G$ is denoted $\tilde{X}$, and is called 
the inertia orbifold \cite{cr}. There is a canonical open-closed embedding 
$G\hookrightarrow\Lambda G$ of groupoids. In fact, this embedding 
is induced by the partition of $B^{(0)}$ into so-called sectors of $G$:
\begin{equation}
\label{orbits}
  B^{(0)}=\coprod \mathcal O,
\end{equation}
where each $\mathcal O$ is a $G$-saturated open-closed subset of $B^{(0)}$
and minimal among such sets with respect to set-theoretic inclusion. 
We denote the set of sectors of $G$ by $\operatorname{Sec} (G)$.
Note that the above decomposition of $B^{(0)}$  
induces similar decompositions of the inertia groupoid $\Lambda G$
and the orbifold $\tilde{X}$, where each component of $\Lambda G$
can be identified as $\orbit\rtimes G$.
The components besides $G_0\subset \Lambda G_0$ are called 
``twisted sectors''. The twisted sectors play an important role in 
orbifold cohomology, which we will now define. In the following, let 
$\ell(\orbit)$ denote the codimension of a sector $\mathcal{O}$ inside $G$.
\begin{definition}
  Let $X$ be an orbifold represented by a groupoid $G$.
  The orbifold cohomology groups and the orbifold cohomology groups
  with compact support of $X$ are defined as
\begin{eqnarray*}
    H^\bullet_{\text{\rm\tiny orb}}(X,\C)\!\!&:=&\!\! H^\bullet(\Lambda G,\C)
    \: = \:
    \bigoplus_{\orbit}H^\bullet(\orbit\rtimes G,\C),\\
    H^\bullet_{\text{\rm\tiny orb,c}}(X,\C)\!\!&:=&\!\!
    \bigoplus_{\mathcal{O}}H^{\bullet-\ell(\orbit)}_\text{\tiny c}
    (\orbit\rtimes G,\C).
\end{eqnarray*}
\end{definition}
Notice that the right hand sides are equal to  
cohomology groups of the inertia orbifold $\tilde{X}$. 
Our shifting of degrees differs from 
\cite{cr}, but likewise we have a Poincar\'{e} duality 
$H^k_{\text{\rm\tiny orb}}(X,\C)\times 
 H^{\dim X-k}_{\text{\rm\tiny orb,c}}(X,\C)\rightarrow \C$.
In \cite{cr}, a remarkable cup product is defined on these 
orbifold cohomology groups, which however will not be used in this paper. 
\subsection{Cyclic homology}
\label{ch}
Here we give a quick review of the basic definitions of Hochschild and 
cyclic (co)homology. For more details, one should consult e.g.~\cite{loday}.

\textit{Cyclic objects.} Let $r\in \mathbb{N}^* \cup\{ \infty\}$. An 
$r$-cyclic object in a category then is a simplicial object 
$(X_\bullet,d,s)$ together with automorphisms 
(cyclic permutations) $t_k:X_k\rightarrow X_k$ satisfying the identities
\begin{eqnarray*}
  d_it_{k+1}&=&
  \begin{cases} 
    t_{k-1}d_{i-1}& \text{for $i\neq 0$,}\\
    d_k & \text{for $i=0$,}  
  \end{cases}
  \\
  s_it_k&=&
  \begin{cases} 
    t_{k+1}s_{i-1} & \text{for $i\neq 0$,} \\
    t^2_{k+1}s_k & \text{for $i=0$,}
  \end{cases}\\
  t_k^{r\,(k+1)}&=&1, \quad \text{if $r\neq \infty$}.
\end{eqnarray*}
Just like a simplicial object is a contravariant functor from the 
simplicial category, an $r$-cyclic object is nothing but a functor from 
the $r$-cyclic category (cf.~\cite{connescyclic}).

\textit{Mixed complexes.} A mixed complex $(X_\bullet,b,B)$ in an 
abelian category is a graded object $(X_k)_{k\in \N}$ equipped with maps 
$b:X_k\rightarrow X_{k-1}$ of degree $-1$ and $B:X_k\rightarrow X_{k+1}$ 
of degree $+1$ such that $b^2=B^2=bB+Bb=0$. 
A mixed complex gives rise to a first quadrant double complex
$B_{\bullet,\bullet} (X)$
\begin{displaymath}
\begin{diagram}
\node{}\arrow{s,r}{b}
\node{}\arrow{s,r}{b}
\node{}\arrow{s,r}{b}
\node{}\arrow{s,r}{b}
\\
\node{X_3}\arrow{s,r}{b}
\node{X_2}\arrow{s,r}{b}\arrow{w,t}{B}
\node{X_1}\arrow{s,r}{b}\arrow{w,t}{B}
\node{X_0}\arrow{w,t}{B}
\\
\node{X_2}\arrow{s,r}{b}
\node{X_1}\arrow{s,r}{b}\arrow{w,t}{B}
\node{X_0}\arrow{w,t}{B}
\\
\node{X_1}\arrow{s,r}{b}
\node{X_0}\arrow{w,t}{B}
\\
\node{X_0}
\end{diagram}
\end{displaymath}
\begin{definition}
\label{hch}
The Hochschild homology $HH_\bullet(X)$ of a mixed complex 
$X=(X_\bullet,b,B)$ is the homology
of the $(X_\bullet,b)$-complex, which sometimes will also be denoted by 
$(C_\bullet (X),b)$. The cyclic homology $HC_\bullet(X)$ is
defined as the homology of the total complex associated to the double 
complex $B_{\bullet,\bullet}(X)$.
\end{definition}
Looking at the double complex above, it is clear that there is a short 
exact sequence of complexes 
$$
  0\longrightarrow (X_\bullet,b)\overset{I}{\longrightarrow}
  (\tot_\bullet B_{\bullet,\bullet} (X),b+B)\overset{S}{\longrightarrow}
  (\tot_\bullet B_{\bullet,\bullet} (X)[-2],b+B)
  \longrightarrow 0.
$$ 
With the definitions above, the associated long exact sequence in homology 
reads 
$$
  \ldots\overset{B}{\longrightarrow} HH_k(X)\overset{I}{\longrightarrow}
  HC_k(X)\overset{S}{\longrightarrow}HC_{k-2}(X)\overset{B}{\longrightarrow} 
  HH_{k-1}(X)\overset{I}{\longrightarrow}\ldots \: .
$$ 
This sequence is called the SBI sequence and relates Hochschild and cyclic 
homology. Stabilizing with respect to the shift operator $S$, the homology 
of the inverse limit complex 
$$
  \lim_{\overset{\longleftarrow}{k}}\tot_\bullet B_{\bullet,\bullet} (X)[-2k],
$$ 
is called the periodic cyclic homology $HP_\bullet(X)$. 
Note that periodic cyclic homology is only $\Z_2$-graded. Alternatively, 
it is the homology of 
the super complex $\hat{X}=\prod_k X_k$ with differential $B-b$.

In case $r\neq \infty$, an $r$-cyclic object $(X_\bullet,d,s,t)$ gives 
rise to a mixed complex (cf.~\cite[A.3.2]{feitsy} and 
\cite[3.1.2]{crainic}). Define 
$$
  b':=\sum_{i=0}^{k-1}(-1)^id_i, \quad
  b:=\sum_{i=0}^k(-1)^id_i, \quad
  N=\sum_{i=0}^{(k+1)r-1} (-1)^{ik}t_k^i,
$$ 
and finally put $B=(1+(-1)^kt_k)sN$. One checks that the triple 
$(X_\bullet,b,B)$ is a mixed complex.
\begin{example}
\label{ex}
Let $A$ be a unital algebra equipped with an automorphism $\alpha$. Define 
the simplicial object $A^\natural_\alpha= (A^\natural_{\alpha,\bullet},d,s)$
by 
$ A^\natural_{\alpha,k} := A^{\otimes(k+1)}$ with face maps given by
$$
  d_i(a_0\otimes\ldots\otimes a_k)=
  \begin{cases}
    a_0\otimes\ldots\otimes a_ia_{i+1}\otimes\ldots\otimes a_k, & 
    \text{ if $0\leq i\leq k-1$,}\\ 
    \alpha(a_k)a_0\otimes\ldots\otimes a_{k-1},&\text{ if $i=k$,}
  \end{cases}
$$ 
and degeneracies given by
$$ 
  s_i(a_0\otimes\ldots \otimes a_{k+1})=
  a_0\otimes\ldots\otimes a_i\otimes 1\otimes a_{i+1}\otimes\cdots\otimes a_k.
$$ 
The differential on the associated Hochschild complex 
$\big( C_\bullet (A_\alpha) ,b_\alpha\big) :=
 \big( A^\natural_{\alpha,\bullet },b_\alpha \big)$ then reads as follows:
\begin{equation}
\begin{split}
\label{twhs}
  b_\alpha & \: (a_0\otimes\cdots\otimes a_k) = \\
  & = \sum_{i=0}^ka_0\otimes\cdots\otimes a_ia_{i+1}\otimes\cdots\otimes a_k
  +(-1)^k\alpha(a_k)a_0\otimes\cdots\otimes a_{k-1}.
\end{split}
\end{equation}
Its homology is denoted by  $HH_\bullet(A_\alpha)$. The map
$$
  t_k(a_0\otimes\cdots\otimes a_k)=
  \alpha(a_k)\otimes a_0\otimes\cdots\otimes a_{k-1}
$$ 
defines an $r$-cyclic structure, where $r$ denotes the order of $\alpha$. 
When $r$ is finite, one defines  Hochschild, cyclic and periodic cyclic 
homology as in Definition \ref{hch} using this cyclic object. When 
$\alpha=1$, the chain complex 
$\big( C_\bullet(A),b \big):=\big(A^\natural_\bullet,b\big)$ is
nothing but the usual Hochschild complex. The double complex 
$B_{\bullet,\bullet} (A)$ associated to the mixed complex 
$\big( A^\natural_\bullet,b,B\big)$ is called Connes' $(b,B)$-complex. 
In this case one denotes the homologies simply by 
$HH_\bullet(A)$, $HC_\bullet(A)$, $HP_\bullet(A)$.
\end{example}
\subsection{Cyclic homology of orbifold groupoids}
\label{chog}
Here we briefly review the computations of \cite{bn,crainic} of the
Hochschild and cyclic homology of \'etale groupoids. To be precise,
this is the homology of the convolution algebra of the groupoid that means 
of $\mathcal C_\text{c}^\infty(G)$ with the multiplication 
\begin{equation}
\label{convolution}
(a_1*a_2)(g)=\sum_{g_1g_2=g} a_1(g_1)a_2(g_2), \quad \text{where 
$a_1,a_2 \in \mathcal C_\text{c}^\infty(G)$, $g \in G$}.
\end{equation}

Let us first describe the general idea behind the computation, which we
will need in particular in Sec.~\ref{cyhom} when computing the Hochschild and
cyclic homology of deformations of the convolution algebra.  

The convolution algebra is a special case of the ``crossed product''
algebra $\A\rtimes G$ associated to any c-soft $G$-sheaf $\mathcal A$
of unital algebras. As a vector space, one has
$\mathcal A\rtimes G:=\Gamma_\text{c} (G_1,s^*\A)$, and the multiplication $*$
is determined by 
\begin{equation}
\label{crossedprd}
  [a_1 * a_2 ]_g= \sum_{g_1 \, g_2=g}\big( [a_1]_{g_1}g_2 \big) [a_2]_{g_2}, 
  \quad 
  \text{for $a_1,a_2 \in \Gamma_\text{c}(G_1,s^*\A)$, $g\in G$},
\end{equation}
where $[a]_g$ denotes the germ of $a$ at $g$.
Indeed, for the sheaf of smooth functions, one recovers the 
convolution algebra \eqref{convolution} in this way. 
The computations of the cyclic homology involve one piece of machinery,
introduced in full generality in \cite{crainic}, that we briefly introduce 
now.

\medskip

\textit{Cyclic groupoids}. A cyclic groupoid is an \'etale groupoid $G$
equipped with a continuous map 
$\theta:G_0\rightarrow G_1,~x\mapsto \theta_x$ such that 
$s(\theta_x)=t(\theta_x)=x$ and $g\theta_x=\theta_yg$ for all $g\in G_1$ 
with $s(g)=x$ and $t(g)=y$. If $\operatorname{ord} (\theta_x)<\infty$ for 
all $x\in G_0$, we say that $G$ resp.~$\theta$ is elliptic.
Of course, any \'etale groupoid has a cyclic 
structure with $\theta_x=u(x)$. The main example of a nontrivial cyclic 
groupoid is $\Lambda G$ with $\theta(g)=(g,g)$ for a loop $g\in B^{(0)}$. 
Notice that $\theta$ is elliptic in this case, since we 
deal only with proper \'etale groupoids, i.e., orbifolds.
Next, define a $\theta$-cyclic sheaf on a cyclic groupoid $G$ to be
an $r$-cyclic object $\mathcal X_\bullet$ in $\mathsf{Sh}(G)$ such that for 
all $x\in G_0$, the morphism 
$[t_k]_x:(\mathcal X_k)_x\rightarrow(\mathcal X_k)_x$
is given by the action of $\theta_x$. As explained in Sec.~\ref{ch}, 
this gives rise to a mixed complex $(\mathcal{X}_\bullet,b,B)$ of sheaves 
on $G$, if $\theta$ is elliptic. 
We then denote the associated Hochschild complex by 
$(\mathcal{C}_\bullet(\mathcal X),b)$ and the associated double complex by 
$\mathcal{B}_{\bullet,\bullet}(\mathcal X)$.
We now define, using the hyperhomology 
from Def.~\ref{hyperhom}.
\begin{definition}
\label{definition}
  The Hochschild and cyclic homology groups of a $\theta$-cyclic sheaf 
  $\mathcal X$ on an elliptic cyclic groupoid $G$ are defined by 
$$
    HH_\bullet(G,\theta;\mathcal X):=\mathbb{H}_\bullet
    (G,(\mathcal C_\bullet (\mathcal X),b)),
    \:\: HC_\bullet(G,\theta;\mathcal X):=\mathbb{H}_\bullet
    (G,(\operatorname{Tot}_\bullet \mathcal{B}_{\bullet,\bullet} 
    (\mathcal X),b+B)).
$$
\end{definition}
Notice that the boundary operators $b$ and $B$ involve the twisting by the 
cyclic structure $\theta$. Of course, for an \'etale groupoid with the 
trivial cyclic structure this twisting is trivial. As explained in 
Sec.~\ref{sheaves}, for a complex of $c$-soft sheaves, homology is 
computed from the Bar-complex of $G$. Therefore, we have the following
consequence of the Eilenberg--Zilber Theorem which will be used 
implicitly several times throughout this paper.
\begin{proposition}
\label{diagsimp}
{\rm (}Cf.~\cite[3.2.8]{crainic}.{\rm )}
  If $\mathcal X$ is a $\theta$-cyclic sheaf on an elliptic cyclic \'etale 
  Lie groupoid $G$ such that each $\mathcal X_k$
  is c-soft, then $HH_\bullet (G,\theta; \mathcal X)$, 
  $HC_\bullet (G,\theta; \mathcal X)$ (and 
  $HP_\bullet (G,\theta; \mathcal X)$) are computed by the diagonal of the 
  bisimplicial vector space $B_\bullet (G,\mathcal X_\bullet)$, 
  i.e.~by the cyclic vector space
  \begin{displaymath}
   \Gamma_\bullet (G, \theta;\mathcal X) : \:\: \cdots \:
   {\substack{\longrightarrow \\ \longrightarrow \\ \longrightarrow \\ 
   \longrightarrow}} \:
   \Gamma_\text{\rm c} (G^{(2)}, \mathcal X_2)  \:
   {\substack{\longrightarrow \\ \longrightarrow \\ \longrightarrow }} \:
   \Gamma_\text{\rm c} (G^{(1)}, \mathcal X_1) \:
   {\substack{\longrightarrow \\ \longrightarrow}}  \:
   \Gamma_\text{\rm c} (G^{(0)}, \mathcal X_0) . 
  \end{displaymath}
\end{proposition} 
Let us now recall the basic idea behind the computations of Hochschild 
and cyclic homology. We have groupoids and maps as follows:
\begin{equation}
\label{diagram}
\begin{diagram}
  \node{}\node{\Lambda G}\arrow{sw,t}{\alpha}\arrow{se,t}{\beta}\node{}\\
  \node{NG}\node{}\node{G}
\end{diagram}
\end{equation}
where $NG$ is the groupoid obtained from $\Lambda G$ by dividing out the 
action of the cyclic structure. It has the same space of objects 
$B^{(0)}$, but $NG_1=\Lambda G_1/\Z$, where 
$n\cdot g=\theta^n_{s(g)}g$, $n\in\Z$. Then, consider the functor
\begin{equation}
\label{transf}
  \tau :=\alpha_!\circ\beta^{-1}: \mathsf{Sh} (G)\rightarrow \mathsf{Sh}(NG).
\end{equation}
Any sheaf of algebras $\A$ gives rise to a cyclic sheaf 
$\A^\natural$ in $\mathsf{Sh}(G)$ by putting
$\A^\natural_k:=\Delta_{k+1}^{-1}\A^{\boxtimes (k+1)}$, where
$\Delta_k:G_0\rightarrow (G_0)^k$ is the diagonal embedding. 
\begin{remark}
Notice that
in the case of sheaves of locally convex topological algebras, 
e.g.~the sheaf of smooth functions, one has to make a choice 
with respect to the topology on the tensor product. 
In most cases, one chooses the projective tensor product topology
\cite[Chap.~1, \S 1, Def.~1]{grothendieck},
but in our setting the inductive tensor product topology 
\cite[Chap.~1, \S 3, Def.~3]{grothendieck} is more natural.
To circumvent such subtleties bornological instead of topological tensor 
products could be used. For a discussion on this point see for instance 
\cite{meyer}.
\end{remark}
The pull-back sheaf
$\A^\natural_{\text{\tiny \rm tw}}:=\beta^{-1}\A^\natural$ 
to the cyclic groupoid $\Lambda G$ carries a natural $\theta$-cyclic 
structure; the stalk of $\A^\natural_{\text{\tiny tw}}$ at 
$g\in B^{(0)}$ is the $r$-cyclic vector space $(\A_{s(g)})^\natural_{\theta_g}$ 
as in Example \ref{ex}, with the 
automorphism given by $\theta_g$. 
Applying the functor $\alpha_!$, one can kill 
the twisting and compute the twisted cyclic homology of 
$\A^\natural_{\text{\tiny tw}}$ as sheaf homology over $NG$. 
(This follows from \cite[Prop 3.3.12]{crainic} and the fact that 
$\Lambda G$ is elliptic for $G$ proper.)

The general idea in the computation of Hochschild and cyclic homology of 
the crossed product $\A\rtimes G$ is to relate it to sheaf homology of 
$\A^\natural_\text{\tiny tw}$ over $\Lambda G$. This proceeds in two steps:

{\bf Step I.} \textit{Reduction to loops}.
Consider Burghelea's space 
$$
  B^{(k)}:=\{(g_0,\cdots,g_k)\in G^{(k+1)},~t(g_0)=s(g_k)\}.
$$ 
There is a map $\sigma_k:B^{(k)}\rightarrow
(G_0)^{k+1},~\sigma_k(g_0,\cdots,g_k)=(s(g_0),\cdots,s(g_k))$. 
For any c-soft sheaf of unital algebras $\mathcal A$, 
we define the vector spaces 
\begin{equation}
\label{simpl}
 \Gamma_\text{\rm c}\Lambda^\natural_k \mathcal A := 
 \Gamma_c(B^{(k)},\sigma_k^{-1}\A^{\boxtimes(k+1)}).
\end{equation}
By construction, $\Gamma_\text{\rm c}\Lambda^\natural_k \mathcal A$ is the 
space of global sections (with compact support) of the sheaf
$\Lambda^\natural_k \mathcal A$ on $B^{(k)}$ which  for 
$(g_0,\cdots,g_k) \in B^{(k)}$ has stalks 
$\big( \Lambda^\natural_k \mathcal A  \big)_{(g_0,\cdots,g_k)}$
given by germs
\begin{equation}
  [a_0\otimes \cdots \otimes a_k]_{(g_0,\cdots,g_k)},
\end{equation}
where each $a_i$ is an element
of $\mathcal A (U_i)$ defined over an open neighborhood $U_i$ of $s(g_i)$.
As explained in \cite[Sec.~3.4]{crainic}, the vector spaces 
$\Gamma_\text{\rm c}\Lambda^\natural_k \mathcal A$ carry a 
canonical cyclic structure, 
combining the structure maps from the underlying cyclic manifold 
$B^{(\bullet)}$ with the structure maps from the cyclic sheaf 
$\A^\natural$. 
The associated Hochschild and cyclic homology is denoted by 
$HH_\bullet (\Gamma_\text{\rm c}\Lambda^\natural_\bullet \mathcal A)$ 
and $HC_\bullet (\Gamma_\text{\rm c}\Lambda^\natural_\bullet \mathcal A)$. 

In case of the convolution algebra, that means in case $\mathcal A$ is the 
sheaf of smooth functions with compact support, we use, as stated above, 
the completed inductive topological tensor product $\topotimes$  
in the definition of the cyclic vector space. The completed inductive tensor 
product has the crucial property that 
$$\mathcal C^\infty_c(M)\topotimes \mathcal C_c^\infty(N)\cong 
\mathcal C_c^\infty(M\times N)$$
for two smooth manifolds $M$ and $N$. We
therefore have topological linear isomorphisms
\begin{displaymath}
  (\A\rtimes G)^\natural_k= \Gamma_c(G,s^{-1}\A)\topotimes \cdots
  \topotimes \Gamma_c(G,s^{-1}\A)
  \overset{{ }_\sim\hspace{1em}}{\longrightarrow}
  \Gamma_c(G^{k+1},s^{-1}_{k+1}\A^{\boxtimes{(k+1)}}),
\end{displaymath}
where $s_k=s\times\ldots\times s.$ The ``reduction to loops'' now is the
natural projection
\begin{equation}
\label{reductiontoloops}
  p:\Gamma_c(G^{k+1},s_{k+1}^{-1}\A^{\boxtimes (k+1)})\rightarrow 
  \Gamma_c(B^{(k)},
  \sigma_k^{-1}\A^{\boxtimes (k+1)}),
\end{equation}
given by the restriction $B^{(k)}\subset G^{k+1}$. 
This defines a map of cyclic vector spaces
$(\A\rtimes G)^\natural_k\rightarrow \Gamma_\text{\rm c}
\Lambda^\natural_k\mathcal A$. 
The essential point proved in \cite[Prop.~3.2]{bn} and 
\cite[Prop.~4.1.1]{crainic}
is that for $\mathcal A$ the sheaf of smooth functions this map
induces isomorphisms in Hochschild and cyclic homology.

{\bf Step II.} \textit{Relation to sheaf cohomology}.
The Hochschild and cyclic homology of the cyclic vector space 
$\Gamma_c\Lambda^\natural_\bullet\A$, defined in \eqref{simpl}, turn out 
to be related to sheaf homology on the inertia groupoid $\Lambda G$:
\begin{proposition}
\label{propcyhomcycsheaf}
  {\rm (}Cf.~\cite{crainic}.{\rm )}
  For any c-soft sheaf of unital algebras $\A$ on a proper \'etale
  groupoid $G$, there are natural isomorphisms 
  \begin{equation}
  \begin{split}
    & HH_\bullet (\Gamma_\text{\rm c}\Lambda^\natural_\bullet \mathcal A)
    \cong HH_\bullet (\Lambda G,\theta;\A^\natural_{\text{\rm \tiny tw}})
    \cong HH_\bullet (NG,\tau (\A^\natural)),
    \\
    & HC_\bullet (\Gamma_\text{\rm c}\Lambda^\natural_\bullet \mathcal A)
    \cong HC_\bullet (\Lambda G,\theta;\A^\natural_{\text{\rm \tiny tw}})
    \cong HC_\bullet (NG,\tau (\A^\natural)).
  \end{split}
  \end{equation}
\end{proposition}
\begin{sketchproof} 
  There is an isomorphism $B^{(k)}\cong (\Lambda G)_k$ of cyclic manifolds, 
  which induces an isomorphism of $\infty$-cyclic vector spaces
  $$
    \Gamma_\text{\rm c} \Lambda^\natural_k \mathcal A =
    \Gamma_c(B^{(k)},\sigma_k^{-1}\A^{\boxtimes(k+1)})\rightarrow
    \Gamma_c(\Lambda G_k,
    \beta^{-1}\A^\natural_k).
  $$ 
Over the stalk at $(g_0,\cdots ,g_k) \in B^{(k)}$, 
the isomorphism is given as follows: 
\begin{equation}
\label{isocyclicvs}
\begin{split}
  \big( \Lambda^\natural_k\mathcal A\big)_{(g_0,\cdots ,g_k)}
  & \rightarrow 
  \big( \mathcal A^\natural_{\text{tw},k} 
  \big)_{(g_1\cdots g_k g_0, g_1,\cdots, g_k)} 
  = \big( \beta^{-1} \mathcal A^\natural_k 
    \big)_{(g_1\cdots g_k g_0, g_1,\cdots, g_k)}  \\
   [a_0 \otimes \cdots  \otimes a_k]_{(g_0,\cdots, g_k)}
  & \mapsto
  [(a_0 g_1\cdots g_k g_0)\otimes \cdots \otimes
  (a_k g_0)]_{(g_1\cdots g_k g_0, g_1,\cdots, g_k)}.
\end{split}
\end{equation}
By Proposition \ref{diagsimp}, the cyclic vector spaces on the right hand 
side computes the $\theta$-twisted Hochschild and cyclic homology of 
$\A^\natural_{\text{\tiny tw}}$ on $\Lambda G$.
By applying $\alpha_!$ one obtains the second isomorphism
in Hochschild and cyclic homology.
\end{sketchproof}

Applied to the sheaf of smooth functions, the twisted 
Hochschild--Kostant--Rosenberg Theorem \cite[Lem.~3.1.5]{crainic}
shows that $\mathcal A^\natural_\text{\tiny tw}$ is quasi-isomorphic
$(\mathcal C_{\Lambda G}^\infty)^\natural$, the usual
cyclic sheaf associated to the commutative sheaf of smooth functions on
$\Lambda G$. 
Applying the Hochschild--Kostant--Rosenberg--Connes quasi-isomorphism 
$(\mathcal C_\bullet (\mathcal C_{\Lambda G}^\infty),b)\cong 
 (\Omega^\bullet_{\Lambda G},0)$, 
which turns the $B$-operator into the de Rham differential, one finds
the additive isomorphism
\begin{displaymath}
  HC_\bullet (\mathcal A \rtimes G)
  \cong H^\bullet_\text{\tiny c}(\Lambda G ,\C) =
  H^\bullet_\text{\tiny c}(\tilde X,\C) \: .
\end{displaymath}

\textit{Localization}. 
As explained in \cite[Prop.~3.3]{bn}, the partition
\eqref{orbits} of $B^{(0)}$ into sectors induces decompositions 
\begin{displaymath}
\begin{split}
  & HH_\bullet(\mathcal A \rtimes G)\cong
  \bigoplus_\orbit HH_\bullet(\mathcal A \rtimes G)_\orbit 
  := \bigoplus_\orbit HH_\bullet ( \Lambda G_{|\mathcal O}, \theta ; 
  \mathcal A^\natural_\text{\tiny tw}),
  \\
  & HC_\bullet(\mathcal A \rtimes G)\cong
  \bigoplus_\orbit HC_\bullet(\mathcal A \rtimes G)_\orbit 
  := \bigoplus_\orbit HC_\bullet ( \Lambda G_{|\mathcal O}, \theta ; 
  \mathcal A^\natural_\text{\tiny tw}),
\end{split}
\end{displaymath} 
and similar for periodic cyclic homology. Of course, this can be read 
off from the final result of the computation in terms of orbifold cohomology,
but it also follows by acting with idempotents 
$e_\orbit \in (\beta^{-1} \mathcal A) (\Lambda G)$ having support 
$\mathcal O$ on the complexes 
$\Gamma_\bullet (\Lambda G,\theta;\A^\natural_{\rm \tiny tw})$. In fact, this 
works for any fine sheaf of unital algebras.
\subsection{Quantization of proper \'etale groupoids}
\label{sec:intr-ncq} 
In \cite{ta:lemma}, the last author 
considered deformation quantization of a pseudo \'etale
groupoid and proved that one can construct star products for such
groupoids. As a special case one obtains that every proper \'etale 
Lie groupoid with an invariant Poisson structure has a formal deformation 
quantization.
In this section, we will recall the basic concepts and constructions from
\cite{ta:lemma}.

\begin{definition}
\label{dfn:poisson} (Cf.~\cite{bl-ge:quantization} and \cite{xu}.)
  A Poisson structure on an associative complete locally convex
  topological algebra $A$ is an element $[\Pi]$
  of the (continuous) Hochschild cohomology group $H^2(A, A)$ such 
  that the cohomology class of the Gerstenhaber bracket
  $[\Pi, \Pi]$ vanishes.
\end{definition}

\begin{remark}
  If $\Pi \in Z^2 (A,A)$ is a Hochschild cocycle representing a 
  Poisson structure on $A$, one has $[\Pi, \Pi] = \delta (\Theta)$ for some 
  Hochschild cochain $\Theta$ on $A$; in the following we will 
  occasionally make use of this fact. By slight abuse of language, 
  we sometimes also call $\Pi$ a Poisson structure on $A$.
\end{remark}

\begin{definition}
\label{dfn:quant-noncommutative-poisson}
  (Cf.~\cite{bl-ge:quantization,ta:lemma}.)
  Let $(A, [\Pi])$ be a noncommutative Poisson algebra,  and $A[[\hbar]]$ the
  space of formal power series with coefficients in $A$. A
  {\it formal deformation quantization} of $(A, [\Pi])$  
  (or in other words {\it star product})
  then is an associative product 
\begin{displaymath}
  \star : A[[\hbar ]] \times A[[\hbar ]] \rightarrow A[[\hbar ]],
  \quad (a_1,a_2) \mapsto a_1 \star a_2 = \sum_{k=0}^{\infty}
  \hbar^k c_k (a_1,a_2)
\end{displaymath}
satisfying the following properties:
\begin{enumerate}
\item 
  Each one of the maps $c_k: A[[\hbar]]\otimes A[[\hbar]]\to A[[\hbar]]$ is 
  $\complex[[\hbar]]$-bilinear.
\item 
  One has $c_0(a_1, a_2)=a_1 \cdot a_2$ for all $a_1, a_2\in A$.
\item
  The relation
  $$
    a_1\star a_2 -c_0 (a_1,a_2) - \frac i2 \hbar  \Pi (a_1, a_2) 
    \in \hbar^2 A[[\hbar]]
  $$
  holds true for some representative $\Pi\in Z^2(A, A)$ of the Poisson 
  structure and  all $a_1,a_2\in A$.
\end{enumerate}
\end{definition}
From now on we consider a proper \'etale Lie groupoid $G$ and let 
$\mathcal A$ denote the sheaf of smooth functions on $G_0$.
\begin{definition}
\label{def:poisson-gpd} A Poisson (resp.~symplectic) structure on 
$G$ is a Poisson (resp. symplectic) structure $\Pi$ on the
unit space $G_0$ which is invariant under the local diffeomorphisms
induced by the source and target maps.
\end{definition}
One easily checks that in the symplectic case, this notion is equivalent 
to the definition of a symplectic orbifold. Note that an invariant Poisson 
bivector on $G_0$ has a canonical lift to a Poisson bivector on $G_1$. 
Having this in mind define a Hochschild 2-cochain on $\mathcal A \rtimes G$
by 
\begin{equation}
\label{def:pi} 
  \widetilde\Pi ( a_1 , a_2)(g)= \sum_{g_1 \, g_2 =g} \, 
  \Pi(g )\big( [da_1]_{g_1} \otimes [da_2]_{g_2} \big), \quad g\in G_1, \: \:
  a_1,a_2 \in \mathcal A \rtimes G ,
\end{equation}
where $[da_1]_{g_1}$ and $[da_2]_{g_2}$ have been pulled back to $g$ along 
the maps $t$ and $s$.
In \cite{ta:lemma}, it was proved that this Hochschild 2-cochain
gives rise to a Poisson structure on the convolution algebra 
indeed. For convenience, we will simply denote the Poisson structure
$\widetilde\Pi$ on $\mathcal A \rtimes G$  by $\Pi$ as well; 
this will not cause any confusion. As proved in \cite{bl-ge:quantization}, 
the center of a noncommutative Poisson algebra carries a natural Poisson 
structure in the commutative sense. For proper \'etale Lie groupoids,
the center equals $\mathcal C^\infty (X)$, the smooth functions on the 
orbifold with the Poisson structure considered in \cite{P2}.

In \cite{ta:lemma} it has been shown that the above Poisson structure 
on the groupoid algebra of a proper \'etale groupoid admits a formal
deformation quantization. Such a deformation can be constructed as follows: 
first construct a deformation quantization of the Poisson manifold $G_0$, 
invariant under the action of the groupoid. In the symplectic case this 
can be done by Fedosov's construction \cite{fe:const} associated to an 
invariant symplectic connection. This defines a fine, so in particular
c-soft, sheaf of noncommutative algebras $\A^\hbar\in\mathsf{Sh}(G)$. 
The associated 
crossed product algebra $\A^\hbar\rtimes G$, as in \eqref{crossedprd}, 
quantizes the convolution algebra $\mathcal A \rtimes G$ with the Poisson 
structure \eqref{def:pi}. We denote the multiplication on 
$\A^\hbar\rtimes G$ which combines the star product $\star$ on 
$\mathcal A^\hbar$ with convolution on 
$G$ by $\star_c$. Notice that $\Gamma_{{\rm\tiny inv ,c}}(\A^\hbar)$ is the 
deformation quantization of $\mathcal C_\text{c}^\infty(X)$ with 
the induced Poisson structure as the center of 
$\mathcal C_\text{c}^\infty(G)$; the thus obtained star product algebra
coincides with the formal deformation quantization studied in \cite{P2}.
\section{Hochschild cohomology of \'etale groupoids}
\label{hceg}
Given a proper \'etale Lie groupoid $G$, we will determine in this section 
the Hochschild cohomology 
$H^\bullet (\A \rtimes G, \A \rtimes G)$, where $\A$ is the $G$-sheaf of 
smooth functions on $G_0$. Recall that the (continuous) Hochschild cohomology
of $\A \rtimes G$ with values in a locally convex topological 
$(\A \rtimes G)$-bimodule $\mathcal M$
is defined as the cohomology of the cochain complex 
$\big( C^\bullet (\A \rtimes G , \mathcal M ), \beta \big)$,
where 
\begin{displaymath}
  C^k (\A \rtimes G , \mathcal M ) = \Hom_{(\A \rtimes G) - (\A \rtimes G)} 
  \big( (\A \rtimes G)^{\topotimes (k+2)} , \mathcal M \big),
\end{displaymath}
and $\beta$ is the standard Hochschild coboundary map 
(see {\cite[Sec.~1.5]{loday}).
Hereby, we have denoted by 
$\Hom_{(\A \rtimes G) - (\A \rtimes G)} (\mathcal N, \mathcal M)$
the vector space of continuous $(\A \rtimes G)$-bimodule maps between two 
locally convex topological $(\A \rtimes G)$-bimodules $\mathcal N$ and 
$\mathcal M$. 
\begin{remark}
  Even though $\A \rtimes G$ is usually nonunital,
  the standard complexes derived from the Bar resolution can be used
  to compute Hochschild (co)homology, since $\A \rtimes G$
  has local units, hence is H-unital (cf.~\cite[Prop.~2]{cm2001}). 
  The same holds for the deformed convolution algebra 
  $\A^{\hbar} \rtimes G$. Later, we will tacitly make use of this fact
  when we determine the Hochschild homology of $\A^{\hbar} \rtimes G$.
\end{remark}
  
\noindent
For the computation of $H^\bullet (\A \rtimes G, \A \rtimes G)$
we proceed in several steps.

{\bf Step 1.} In the following we provide a more convenient description 
of the cochain complex $C^\bullet := C^\bullet (\A \rtimes G ,\A\rtimes G)$
and identify it with the complex of global section spaces of some 
sheaf complex $\mathcal K^\bullet$ on the orbit space $X := G_0/G_1$.
To this end first note 
that the Fr\'echet space $\mathcal C^\infty (G_1)$ inherits from the 
convolution algebra
$\A \rtimes G$ the structure of a locally convex $(\A\rtimes G)$-bimodule. 
More generally, observe that for every open subset
$U \subset  X$ the Fr\'echet space $\mathcal C^\infty (U_1)$ of smooth 
functions on the preimage $U_1 := (\pi \circ s)^{-1} (U) \subset G_1$
becomes a topological 
$(\A\rtimes G)$-bimodule by the formula
\begin{equation}
\label{bimodstr}
  (a_1 * a * a_2 )\, (g) =
  \sum_{h_1 \, h h_2 = g}\, a_1 (h_1)\, a(h)\, a_2 (h_2),
\end{equation} 
where $a_1, a_2 \in \A\rtimes G$, $a \in \mathcal C^\infty (U_1)$ and 
$g\in G_1$.
Now, choose for every compact $K\subset G_1$ a smooth function 
$\varphi_K : G_0 \rightarrow [0,1]$ with compact support such that
$\varphi_K (x) = 1$ for all $x \in s(K) \cup t(K)$.
Let $\varphi_K \delta_u: G_1 \rightarrow [0,1]$ be the function which coincides
with $\varphi_K$ on $u(G_0)$  and vanishes elsewhere. Then $\varphi_K \delta_u$
is an element of the convolution algebra, hence for every cochain
$\Phi \in C^k$ one obtains a continuous linear map
$\check\Phi : (\A \rtimes G)^{\topotimes k} 
\rightarrow \mathcal C^\infty (G_1)$
by putting 
\begin{displaymath}
 \check\Phi (a_1 \otimes \cdots \otimes a_k)_{| K} = 
 \Phi (\varphi_K \delta_u \otimes
 a_1 \otimes \cdots \otimes a_k \otimes \varphi_K \delta_u)_{| K},
\end{displaymath}
where $K$ runs through the compact subsets of $G_1$. One checks immediately,
that $\check\Phi$ is well-defined and continuous indeed. Moreover, it is 
easy to prove that the map $\check{\hspace{1em}}$ identifies $C^k $ with 
$\Hom ( (\A \rtimes G)^{\topotimes k}, \mathcal C^\infty (G_1))$. 
Having this identification in mind we  now put for every open $U\subset X$:
\begin{equation}
  \mathcal K^k (U) := \Hom (( \mathcal A \rtimes G)^{\topotimes k}, 
  \mathcal C^\infty (U_1)), \quad U_1:= (\pi \circ s)^{-1} (U),
\end{equation} 
where $\mathcal C^\infty (U_1)$ carries the $\mathcal A \rtimes G$-bimodule
structure given by Eq.~(\ref{bimodstr}). Since the Hochschild coboundary 
is functorial with respect to restriction maps, and since the smooth
functions on $G_1$ form a sheaf, $\mathcal K^\bullet$ is a complex of sheaves
on $X$. By the above identification it is clear that $\mathcal K^\bullet (X)$
can be naturally identified with $C^\bullet$.

{\bf Step 2.}
In this part we prove a localization result for 
Hochschild (co)homology of the convolution algebra.
To this end we need some more notation. First let us fix a smooth function 
$\varrho : \R \rightarrow [0,1]$ which has support
in $(-\infty , \frac 34]$ and which satisfies $\varrho (r) =1$ for 
$r\leq \frac 12$. For $\varepsilon >0$ we denote by $\varrho_\varepsilon$ 
the rescaled function $r \mapsto \varrho (\frac s \varepsilon )$. 
Next choose a $G$-invariant metric 
$d$ on $G_0 $ such that $d^2$ is smooth, and set for every 
$k\in \N \cup \{ -1 \}$, $i=1,\cdots, k1$ and $\varepsilon >0$:
\begin{displaymath}
  \Psi_{k,i,\varepsilon} (g_0,g_1,\cdots, g_k) =  \prod_{j=0}^{i-1} 
  \varrho_\varepsilon \big( d^2( s( g_j ) , t (g_{j+1})) \big),
  \quad \text{where $g_{k+1}:= g_0$}. 
\end{displaymath} 
Moreover, put $\Psi_{k,\varepsilon} := \Psi_{k,k+1,\varepsilon}$.

Given a Hochschild chain $c$ 
resp.~a Hochschild cochain  $F$ (of degree $k$) we now define
$\Psi_{k,\varepsilon} c \in  C_k := C_k (\A \rtimes G , \A \rtimes G)$ 
resp.~$\Psi^{k,\varepsilon} F\in  C^k := C^k (\A \rtimes G , \A \rtimes G)$ 
as follows:
\begin{displaymath}
\begin{split}
 &( \Psi_{k,\varepsilon} c) (g_0,g_1, \cdots, g_k) :=
  \Psi_{k,\varepsilon} (g_0,g_1, \cdots, g_k) \cdot
  c (g_0,g_1, \cdots, g_k), \\
 &( \Psi^{k,\varepsilon} F (a_1 \otimes \cdots \otimes a_k ) \, (g_0) :=
  F \big( \Psi_{k,\varepsilon} (g_0^{-1}, -, \cdots,- ) \cdot
 (a_1 \otimes \cdots \otimes a_k) \big) \, (g_0),\\
 & \hspace{12em} (g_0,g_1, \cdots, g_k) \in G^{k+1}, \quad
   a_1, \cdots, a_k \in \mathcal C^\infty_{\text{c}} (G).
\end{split}
\end{displaymath}
One immediately checks then that the operations
$\Psi_{\bullet,\varepsilon}$  and $\Psi^{\bullet,\varepsilon}$
are both chain maps on the Hochschild chain resp.~cochain complex.

Let us now construct a homotopy between the identity operator and
$\Psi_{\bullet,\varepsilon}$ resp.~$\Psi^{\bullet,\varepsilon}$.
To this end define maps 
$\eta_{k,i,\varepsilon}: C_k \rightarrow C_{k+1} $ for $1\leq i \leq k+1$ 
and maps 
$\eta^{k,i,\varepsilon}: C^k \rightarrow C^{k-1} $ for $1\leq i \leq k$ 
as follows:
\begin{displaymath}
\begin{split}
  &\eta_{k,i,\varepsilon} (c)  \, (g_0,g_1 , \cdots , g_{k+1}) =  \\
  & \: = 
  \begin{cases}
     \Psi_{k+1,i,\varepsilon} (g_0,g_1,\cdots, g_{k+1}) \cdot c(g_0,\cdots,
     g_{i-1},g_{i+1}, \cdots , g_{k+1} ) \cdot \delta_u (g_i), & 
     \text{$i< k+1$},\\
     \Psi_{k+1,k+1,\varepsilon} (g_0,g_1,\cdots, g_{k+1}) \cdot 
     c(g_0,\cdots, g_{k} )
     \cdot \delta_u (g_{k+1}), & \text{$i= k+1$},\\
  \end{cases}
\end{split}
\end{displaymath}
and 
\begin{displaymath}
\begin{split}
  &\eta^{k,i,\varepsilon} (F)  (a_1 \otimes \cdots \otimes a_{k-1}) \, 
 (g_0) =  \\  & \: = 
  \begin{cases}
     F \big( \Psi_{k,i,\varepsilon} (g_0^{-1},-,\cdots, -) \cdot 
     (a_1 \otimes \cdots \otimes a_{i-1} \otimes \delta_u \otimes a_{i}
     \otimes \cdots \otimes a_{k-1})\big) (g_0) , & \text{$i<k$},\\
     F \big( \Psi_{k,k,\varepsilon} (g_0^{-1} , - , \cdots, -) \cdot 
     (a_1 \otimes \cdots \otimes a_{k-1} \otimes \delta_u )\big) (g_0), 
     & \text{$i=k$}.\\
  \end{cases}
\end{split}
\end{displaymath}
Hereby, $\delta_u \in \mathcal C^\infty (G)$ denotes the function
\begin{displaymath}
 g \mapsto 
\begin{cases}
   1, & \text{if $g = u(x)$ for some $x\in G_0$}, \\
   0, & \text{else}.
\end{cases}
\end{displaymath} 
By a somewhat lengthy, but straightforward computation one then proves the 
following result.
\begin{proposition}
\label{prophomtploc}
  The maps 
  \begin{displaymath}
  \begin{split}
    H_{k,\varepsilon} := \sum_{i=1}^{k+1} \, (-1)^{i+1} \,
    \eta_{k,i,\varepsilon} : \: &
    C_k \rightarrow C_{k+1} \quad \text{and}\\
    H^{k,\varepsilon} := \sum_{i=1}^{k} \, (-1)^{i+1} \, 
    \eta^{k,i,\varepsilon} : \: &
    C^{k} \rightarrow C^{k-1}
  \end{split}
  \end{displaymath}
form a homotopy between the identity and the localization morphism 
$\Psi_{\bullet,\varepsilon}$ resp.~$\Psi^{\bullet,\varepsilon}$. 
More precisely, 
\begin{align}
    \label{homotopy1}
    (b_{k+1} H_{k,\varepsilon} + H_{k-1,\varepsilon} b_k) c & = 
    c - \Psi_{\bullet,\varepsilon} c \quad 
    \text{for all $c\in C_{k}$},\\
    \label{homotopy2}
    (\beta^{k-1} H^{k,\varepsilon}_\gamma + H^{k+1,\varepsilon} \beta^k) F & = 
    F - \Psi^{\bullet,\varepsilon} F \quad 
    \text{for all $F\in C^{k}$}.
\end{align}
\end{proposition}
\begin{sketchproof}
  Let $d_{k,j} :C_k \rightarrow C_{k-1}$ be the face maps of Example \ref{ex}.
  Then one easily checks the following commutation relations 
  for $i=1$ 
  \begin{displaymath}
    ( d_{k+1,j} \, \eta_{k,1,\varepsilon} \, c) \, (g_0,\cdots , g_k) =
    \begin{cases}
       c (g_0,\cdots , g_k) , & \text{ if $j=0$},\\
       ( \Psi_{k,1,\varepsilon}\, c)\, (g_0,\cdots, g_k),& 
       \text{ if $j=1$}, \\
       ( \eta_{k -1,1,\varepsilon} \, d_{k,j-1} \, c )\, (g_0,\cdots, g_k), 
       & \text{ if $1 < j \leq k+1$},
    \end{cases}
  \end{displaymath}
 for $i=2,\cdots, k$
\begin{displaymath}
    ( d_{k+1,j} \, \eta_{k,i,\varepsilon} \, c) \, (g_0,\cdots , g_k) =
    \begin{cases}
       (\eta_{k-1,i-1,\varepsilon} \, d_{k,j} \, c) \, (g_0,\cdots, g_k),& 
       \text{ if $0 \leq j < i-1$}, \\
       ( \Psi_{k,i- 1,\varepsilon} \, c)\, (g_0,\cdots, g_k) ,  & 
       \text{ if $j=i-1$}, \\
       ( \Psi_{k,i,\varepsilon} \, c)\, (g_0,\cdots, g_k) ,  & 
       \text{ if $j=i$}, \\
       (\eta_{k -1, i ,\varepsilon} \, d_{k,j-1} \, c )\, (g_0,\cdots, g_k),& 
       \text{ if $i < j \leq k+1$},
    \end{cases}
  \end{displaymath}
  and  for $i=k+1$
\begin{displaymath}
    ( d_{k+1,j} \, \eta_{k,k+1,\varepsilon} \, c) \, (g_0,\cdots , g_k) =
    \begin{cases} 
       (\eta_{k -1,k,\varepsilon} \, d_{k,j} \, c)\, (g_0,\cdots, g_k) ,  & 
       \text{ if $0 \leq j <k$}, \\
       (\Psi_{k,k,\varepsilon} \, c ) \, (g_0,\cdots , g_k) ,& 
       \text{ if $j=k$},\\
       ( \Psi_{k,k+1,\varepsilon} \, c ) \, (g_0,\cdots, g_k) ,& 
       \text{ if $j=k+1$}.
    \end{cases}
\end{displaymath}
From these commutation relations one immediately derives
Eq.~(\ref{homotopy1}).

Now let us consider the dual case. Let $\sigma^{k,j}:C^k \rightarrow C^{k+1}$, 
$j=0,\cdots, k+1$ be the face 
maps of the cosimplicial vector space $C^\bullet$, i.e., let
\begin{displaymath}
  \sigma^{k,j} F (a_1\otimes \cdots \otimes a_{k+1}) = 
  \begin{cases}
    a_1 * F  (a_2\otimes \cdots \otimes a_{k+1}), & \text{ if $j=0$},\\
    F(d_{k+1,j} (a_1\otimes \cdots \otimes a_{k+1})),  & 
    \text{ if $1 < j < k+1$},\\
    F(a_1\otimes \cdots \otimes a_{k}) * a_{k+1}, & \text{ if $j=k+1$}.
  \end{cases}
\end{displaymath} 
For $i=2,\cdots,k$ and $j=0$ one then computes
\begin{displaymath}
\begin{split}
  \big( &\eta^{k+1,i,\varepsilon}  \, \sigma^{k,j} \, F) 
  (a_1 \otimes \cdots \otimes a_k)
  \, (g_0) = \\
  & = \big( \sigma^{k,0}\, F \big) (\Psi_{k+1,i,\varepsilon} 
  (g_0^{-1},-\cdots,-) \cdots 
  (a_1 \otimes \cdots \otimes a_{i-1} \otimes 
  \delta_u \otimes a_i \otimes \cdots 
  \otimes a_k ) \, (g_0 ) \\
  & = \sum_{h\,h' = g_0} \, a_1 (h) \cdot F 
  (\Psi_{k+1,i,\varepsilon} (g_0^{-1},h , - ,\cdots,-) \cdot
  (a_2 \otimes \cdots \otimes a_{i-1} \otimes \delta_u \otimes a_i 
  \otimes \cdots 
  \otimes a_k ) \, (h' ) \\
  & = \sum_{h\,h' = g_0} \, a_1 (h) \cdot F 
  (\Psi_{k,i-1,\varepsilon} ((h')^{-1} , - ,\cdots,-) \cdot
  (a_2 \otimes \cdots \otimes a_{i-1} \otimes \delta_u \otimes a_i 
  \otimes \cdots 
  \otimes a_k ) \, (h' ) \\
  & = \big( \sigma^{k-1,j} \, \eta^{k,i-1,\varepsilon} \, F \big) 
  (a_1 \otimes \cdots \otimes a_k)  \, (g_0) .
\end{split}
\end{displaymath}
By computations of this type and the corresponding relations in the 
homology case one obtains for $i=1$ 
  \begin{displaymath}
  \begin{split}
    ( \eta^{k+1,1,\varepsilon} & \, \sigma^{k,j} \, F) 
    (a_1\otimes \cdots \otimes a_k) \,(g_0) =
    \\  = &
    \begin{cases}
       F (a_1\otimes \cdots \otimes a_k) \,(g_0), & \text{ if $j=0$},\\
       F (\Psi_{k,1,\varepsilon} (g_0^{-1},-,\cdots,-) \cdot 
       (a_1\otimes \cdots \otimes a_k) ) (g_0) ,
       & \text{ if $j=1$}, \\
       \big( \sigma^{k-1,j-1} \, \eta_{k,1,\varepsilon} \, F \big) 
       (a_1\otimes \cdots \otimes a_k) 
       \, (g_0) ,  & \text{ if $1 < j \leq k+1$},
    \end{cases}
  \end{split}
  \end{displaymath}
for $i=2,\cdots, k$
\begin{displaymath}
  \begin{split}
    ( \eta^{k+1,1,\varepsilon} & \, \sigma^{k,j} \, F) 
    (a_1\otimes \cdots \otimes a_k) \,(g_0) =
    \\  = &
    \begin{cases}
    \big( \sigma^{k-1,j} \, \eta^{k,i-1,\varepsilon} \, F\big) 
    (a_1 \otimes \cdots \otimes a_k)  \, (g_0), & 
    \text{ if $0 \leq j < i-1$}, \\
     F(\Psi_{k,i-1,\varepsilon} (g_0^{-1},-,\cdots,-) \cdot 
       (a_1\otimes \cdots \otimes a_k) ) (g_0) ,
       & \text{ if $j=i-1$}, \\
     F(\Psi_{k,i,\varepsilon} (g_0^{-1},-,\cdots,-) \cdot 
       (a_1\otimes \cdots \otimes a_k) ) (g_0) ,
       & \text{ if $j=i$}, \\
    \big( \sigma^{k-1,j-1} \, \eta^{k,i,\varepsilon} \, F\big) 
    (a_1 \otimes \cdots \otimes a_k)  \, (g_0), & \text{ if $i < j \leq  k+1$},
    \end{cases}
  \end{split}
  \end{displaymath}
and for $i=k+1$
\begin{displaymath}
  \begin{split}
  ( \eta^{k+1,1,\varepsilon} & \, \sigma^{k,j} \, F) 
  (a_1\otimes \cdots \otimes a_k) \,(g_0) =
    \\  = &
    \begin{cases}
      \big( \sigma^{k-1,j} \, \eta^{k,k,\varepsilon} \, F\big) 
      (a_1 \otimes \cdots \otimes a_k)  \, (g_0), & 
      \text{ if $0 \leq j < k$}, \\
      F(\Psi_{k,k,\varepsilon} (g_0^{-1},-,\cdots,-) \cdot 
      (a_1\otimes \cdots \otimes a_k) ) (g_0) ,
      & \text{ if $j=k$}, \\
      F(\Psi_{k,k+1,\varepsilon} (g_0^{-1},-,\cdots,-) \cdot 
      (a_1\otimes \cdots \otimes a_k) ) (g_0) ,
      & \text{ if $j=k+1$}.
    \end{cases}
 \end{split}
\end{displaymath}
 Using $\beta^k = \sum (-1)^j \sigma^{k,j}$, these commutation relations 
 immediately entail Eq.~(\ref{homotopy2}).
\end{sketchproof}  

Denote by $C_k^\varepsilon$ the subspace of all Hochschild chains with 
support in the complement of
$$
 U_{k+1,\varepsilon} :=\{ (g_0,\cdots,g_k) \in G^{k+1} \mid 
 d^2 (s(g_0),t(g_1)) + \cdots + d^2 (s (g_k) , t(g_0)))<\varepsilon\}
$$
and by $C^k_\varepsilon$ the space of all Hochschild cochains having 
support in in the complement of
$$
 \widetilde{U}_{k+1,\varepsilon} :=\{ (g_0,\cdots,g_k) \in G^{k+1} \mid 
 d^2 (s(g_0^{-1}),t(g_1)) + \cdots + d^2 (s (g_k) , t(g_0^{-1})))<
 \varepsilon\}.
$$
Moreover, let $C_k^0$ resp.~$C^k_0$ be the union of all
$C_k^\varepsilon$ resp.~$C^k_\varepsilon$, where $\varepsilon$ runs
through all positive real numbers. Then the proposition entails
\begin{corollary}
  The subcomplexes $C_\bullet^0$ and $C^\bullet_0$ are acyclic.
  In particular, the quotient maps
  \begin{displaymath}
    C_\bullet \rightarrow C_\bullet / C_\bullet^0 \quad \text{ and } \quad
    C^\bullet \rightarrow C^\bullet / C^\bullet_0
  \end{displaymath}   
  are quasi-isomorphisms.
\end{corollary}
\begin{remark}
  Originally, Brylinski--Nistor have shown in \cite[Prop.~3.2]{bn} that  
  $C_\bullet \rightarrow C_\bullet / C_\bullet^0$ is a quasi-isomorphism 
  and used this result to compute the Hochschild homology 
  $HH_\bullet (\A\rtimes G)$.
\end{remark}
\begin{remark}
\label{remcommloc}
  In the case, where $G$ is the Lie groupoid whose objects and arrows 
  are given by the points of a smooth manifold $M$, one recovers the 
  well-known localization scheme for Hochschild homology 
  \`a la Teleman \cite{teleman}. In the following, we will 
  freely make use of this fact.
\end{remark}

{\bf Step 3.} In the third step we restrict our considerations to the case, 
where $G$ is a transformation groupoid $\Gamma \ltimes M$ of a 
finite group $\Gamma$ acting on a smooth manifold $M$.
Recall that then $G_1 = \Gamma \times M$, $G_0=M$ and that every 
$\gamma \in \Gamma$ acts on $A := \mathcal C^\infty (M)$ by 
\begin{displaymath} 
  \gamma a (p) = a(\gamma^{-1} p), \quad 
  \text{where $a \in A$, $p\in M$}.
\end{displaymath}
Moreover, every element $a$ of the convolution algebra 
$A \rtimes \Gamma $ has a unique representation of the form
\begin{equation}
\label{repelconvalg}
  a = \sum_{\gamma \in \Gamma} f_\gamma \, \delta_\gamma ,
\end{equation}
where $f_\gamma \in A $ and  where  $f_\gamma \, \delta_\gamma$
is the function which satisfies 
$f_\gamma \, \delta_\gamma (\gamma,p) = f_\gamma (\gamma p)$
and vanishes elsewhere. One easily computes that then
\begin{equation}
\label{deltaprod}
  f_1 \delta_{\gamma_1} * f_2 \delta_{\gamma_2} = f_1 (\gamma_1 f_2) \,
  \delta_{\gamma_1 \gamma_2} \quad \text{for all $f_1,f_2\in A $ and 
  $\gamma_1,\gamma_2\in \Gamma$}.
\end{equation}
Concerning the topological tensor product considered, one should observe
in the following that the completed inductive tensor product
$A \topotimes A$ and the completed projective tensor product
$A \hat\otimes A$ coincide, since $A$ is a (nuclear) Fr\'echet space.

\begin{lemma}
\label{gammaactionhom}
  Let $\Gamma$  act on the space 
  $\Hom (A^{\hat\otimes k}, A \rtimes \Gamma)$ 
  (of continuous linear maps) as follows:
  \begin{displaymath}
    (\gamma \phi) ( f_1 \otimes \cdots \otimes f_k) 
    = \delta_\gamma *
    \phi ( \gamma^{-1} f_1 \otimes \cdots \otimes \gamma^{-1} f_k) * 
    \delta_{\gamma^{-1}}.
  \end{displaymath}
  Then the relation 
  \begin{equation}
     f \delta_e * \gamma \phi * f' \delta_e = 
     \gamma \big( (\gamma^{-1} f) \delta_e *\phi  * 
     (\gamma^{-1} f') \delta_e \big)      
  \end{equation}
  holds true for all 
  $\phi \in \Hom (A^{\hat\otimes k}, A \rtimes \Gamma )$ 
  and $f,f' \in A$.
\end{lemma}
\begin{proof}
  The claim is an immediate consequence of Eq.~(\ref{deltaprod}).
\end{proof}
  Consider now the vector spaces
  $C^{n,m}_\Gamma = \Hom (\C \Gamma^m , \Hom (A^{\hat\otimes n}, 
   A \rtimes \Gamma ))$, where $n,m\in \N$. Using the $\Gamma$-action on 
  $\Hom (A^{\hat\otimes k}, A \rtimes \Gamma )$ from above,
  the simplicial structures coming from group cohomology and 
  Hochschild cohomology then induce on $C_\Gamma^{\bullet,\bullet}$
  the structure of a bicosimplicial vector space as follows 
  (where $\Psi \in  C^{m,n}_\Gamma$):
\begin{align}
  \nonumber
  d_\text{v}^i: \: & C^{m,n}_\Gamma \rightarrow C^{m+1,n}_\Gamma , \quad
  d_\text{v}^i \Psi (\gamma_1, \ldots , \gamma_{m+1}) = \\
  & =
  \begin{cases}
    \gamma_1 ( \Psi (\gamma_2, \ldots , \gamma_{m+1})), & \text{if $i=0$},\\
    \Psi (\gamma_1, \ldots , \gamma_i \gamma_{i+1}, \ldots , \gamma_{m+1}), & 
    \text{if $1\leq i \leq m$},\\
    \Psi (\gamma_1, \ldots , \gamma_{m} ), & \text{if $i=m+1$},
  \end{cases}\\[2mm]
  \nonumber
  d_\text{h}^j: \: & C^{m,n}_\Gamma \rightarrow C^{m,n+1}_\Gamma , \quad
  d_\text{h}^j \Psi (\gamma_1, \ldots , \gamma_m) \, 
  (f_1 \otimes \cdots \otimes f_{n+1} ) = \\
  & =
  \begin{cases}
    f_1 \delta_e * \Psi (\gamma_1, \ldots , \gamma_m) 
    (f_2 \otimes \cdots \otimes f_{n+1} ), & \text{if $j=0$},\\
    \Psi (\gamma_1, \ldots , \gamma_m) 
    (f_1 \otimes \cdots \otimes f_j \, f_{j+1} 
    \otimes \cdots \otimes f_{n+1} ),
    & \text{if $1\leq j \leq n$},\\
    \Psi (\gamma_1, \ldots , \gamma_m ) (f_1 \otimes \cdots \otimes f_n) *
    f_{n+1} \delta_e, & \text{if $j=n+1$},
  \end{cases} \\[2mm]
\nonumber
  s_\text{v}^i: \: & C^{m,n}_\Gamma \rightarrow C^{m-1,n}_\Gamma , \quad
  s_\text{v}^i \Psi (\gamma_1, \ldots , \gamma_{m-1}) = \\
  & = 
 \begin{cases} 
   \Psi (e, \gamma_1, \ldots , \gamma_{m-1}),  & \text{if $i=0$}, \\
   \Psi (\gamma_1, \ldots , \gamma_i , e ,\gamma_{i+1}, \ldots ,
   \gamma_{m-1}),  & \text{if $1\leq i \leq m-1$},
 \end{cases}  \\[2mm]
\nonumber
  s_\text{h}^j: \: & C^{m,n}_\Gamma \rightarrow C^{m,n-1}_\Gamma , \quad
  s_\text{h}^j \Psi (\gamma_1, \ldots , \gamma_{m}) 
 (f_1 \otimes \cdots \otimes f_{n-1} )
  = \\
  & = 
 \begin{cases} 
   \Psi (\gamma_1, \ldots , \gamma_{m})
   (1 \otimes f_1 \otimes \cdots \otimes f_{n-1} ), & \text{if $j=0$}, \\
   \Psi (\gamma_1, \ldots , \gamma_m)
   (f_1 \otimes \cdots \otimes f_j \otimes 1 \otimes f_{j+1} 
   \otimes \cdots \otimes f_{n-1} ) , & \text{if $1\leq j \leq n-1$}.
 \end{cases}  
\end{align}
The $d_\text{v}^i$ (resp.~$d_\text{h}^j$) form the vertical (resp.~horizontal)
face maps of the bicosimplicial vector space, the $s_\text{v}^i$ 
(resp.~$s_\text{h}^j$) the vertical (resp.~horizontal) degeneracies.
Using  Lemma \ref{gammaactionhom} it is easy to show that every vertical
structure map commutes with every horizontal structure map, hence
$C^{\bullet,\bullet}_\Gamma$ is a bicosimplicial vector space indeed.
For example, let us show that $d_\text{v}^0 $ and $d_\text{h}^0$ commute:
\begin{align}
  \nonumber
  \big( d_\text{v}^0 d_\text{h}^0 \Psi \big) & \,
  (\gamma_1, \ldots , \gamma_{m+1})
  (f_1 \otimes \cdots \otimes f_{n+1} ) = \\ \nonumber
  = & \, \delta_{\gamma_1} * (\gamma_1^{-1} f_1 \delta_e) *
  \Psi (\gamma_2, \ldots , \gamma_{m+1} ) (\gamma_1^{-1} f_2 \otimes 
  \cdots \otimes \gamma_1^{-1} f_{n+1}) * \delta_{\gamma_1^{-1}} \\ \nonumber
  = &  \, f_1 \delta_e * \delta_{\gamma_1} *
  \Psi (\gamma_2, \ldots , \gamma_{m+1} ) (\gamma_1^{-1} f_2 \otimes 
  \cdots \otimes \gamma_1^{-1} f_{n+1}) * \delta_{\gamma_1^{-1}}  \\
  = &  \, \big( d_\text{h}^0 d_\text{v}^0 \Psi \big) \,
  (\gamma_1, \ldots , \gamma_{m+1})
  (f_1 \otimes \cdots \otimes f_{n+1} )
\end{align}
At this point recall that the bicosimplicial space
$C_\Gamma^{\bullet,\bullet }$ induces  the structure of a 
cosimplicial space on the diagonal 
$C^\bullet_\Gamma := 
\operatorname{diag} (C_\Gamma^{\bullet,\bullet })$ 
(see \cite[Sec.~8.5]{weibel}). Its structure maps are given 
by $ d^i = d^i_{\text{h}} \, d^i_{\text{v}}$ and
$ s^i = s^i_{\text{h}} \, s^i_{\text{v}}$. 
\begin{proposition}
\label{propsimpmap}
  Define for every $\Phi \in C^k (A\rtimes \Gamma, A \rtimes \Gamma)$ 
  an element $\hat{\Phi} \in C^k_\Gamma = 
  \Hom \big(\C \Gamma^k ,\Hom (A^k,A\rtimes \Gamma) \big) $ as follows:
  \begin{align}
    \nonumber \hat{\Phi} & \,
    (\gamma_1,\ldots,\gamma_k)(f_1 \otimes \cdots \otimes f_k)= \\  & = 
    \Phi \big( (\gamma_1 \, \gamma_2 \ldots \gamma_k f_1) \delta_{\gamma_1}
    \otimes (\gamma_2 \, \ldots \gamma_k f_2) \delta_{\gamma_2} \otimes
    \cdots \otimes (\gamma_k f_k) \delta_{\gamma_k} \big).
  \end{align}   
  Then $\hat{\hspace{1em}}: C^\bullet (A\rtimes \Gamma, A \rtimes \Gamma)
  \rightarrow C^\bullet_\Gamma$ is a cosimplicial map.
\end{proposition} 
\begin{proof}
  Denote by $b^i$ the face maps of Hochschild cohomology, which act on
  a cochain $\Phi \in C^k (A \rtimes \Gamma , A \rtimes \Gamma) $  as follows:
  \begin{align}
    \nonumber
     b^i \Phi & \, (f_1 \delta_{\gamma_1} \otimes \cdots \otimes f_{k+1} 
     \delta_{\gamma_{k+1}} ) = \\ & = 
     \begin{cases}
     f_1 \delta_{\gamma_1} * \Phi(f_2 \delta_{\gamma_2} \otimes \cdots 
     \otimes f_{k+1} \delta_{\gamma_{k+1}} )
     & \text{ if $ i = 0$} \\
     \Phi(f_1 \delta_{\gamma_1} \otimes \cdots \otimes f_i (\gamma_i f_{i+1} ) 
     \delta_{\gamma_i\gamma_{i+1}}
     \otimes \cdots \otimes f_{k+1} \delta_{\gamma_{k+1}} ) & 
     \text{ if $ 1\leq i < k$}, \\
     \Phi(f_1\delta_{\gamma_1}\otimes \cdots \otimes f_k \delta_{\gamma_k} ) 
     * f_{k+1} \delta_{\gamma_{k+1}} & \text{ if $ i = n$} 
   \end{cases}
  \end{align}
  Then compute for $1\leq i <k$:
  \begin{align}
    \nonumber
    (b^i \Phi)\hat{\hspace{0.3em}} & 
    \, (\gamma_1,\ldots, \gamma_{k+1} ) \, (f_1 \otimes \cdots 
    \otimes f_{k+1}) = \\ \nonumber
    & = b^i \Phi \big( (\gamma_1 \ldots \gamma_{k+1} f_1) \delta_{\gamma_1} 
    \otimes \cdots \otimes
    (\gamma_{k+1} f_{k+1}) \delta_{\gamma_{k+1}} \big)  \\ \nonumber
    & = \Phi \big( (\gamma_1 \ldots \gamma_{k+1} f_1)\delta_{\gamma_1} \otimes 
    \cdots \otimes (\gamma_i \, \gamma_{i+1}  \ldots \gamma_{k+1} 
    (a_i a_{i+1})) \delta_{\gamma_i \gamma_{i+1}} \otimes 
    \ldots \otimes \gamma_k f_1)   \big)   \\ 
    & = d^i_{\text{v}} d^i_{\text{h}} \, \hat{\Phi} 
    (\gamma_1, \ldots ,\gamma_{k+1}) 
    (a_1 \otimes \cdots \otimes a_{k+1}).
  \end{align}
  By a similar computation one shows that $\hat{\hspace{1em}}$ preserves 
  all the other face and the degeneracy maps. This proves the claim.
\end{proof}

We now have the tools to show the following result.
\begin{proposition} 
\label{propisohcoh}
{\rm (}Cf.~\cite[Prop.~4.1]{giaquinto}{\rm )}
Let $\Gamma$ be finite group acting by diffeomorphisms on the manifold $M$,
and $A$ the Fr\'echet algebra of smooth functions on $M$.
Then the Hochschild cohomology $H^\bullet (A,A\rtimes \Gamma)$ carries a 
natural $\Gamma$-action such that 
\begin{equation}
\label{eqhhiso}
  H^\bullet (A\rtimes \Gamma, A\rtimes \Gamma) \cong 
  H^\bullet (A, A \rtimes \Gamma)^\Gamma.
\end{equation}
On the level of cochains, this isomorphism is induced by the following 
chain map:
\begin{equation}
  C^k (A \rtimes \Gamma , A \rtimes \Gamma) \rightarrow
  C^k (A , A \rtimes \Gamma), \hspace{0.5em} F \mapsto 
  \big( a_1\otimes \ldots \otimes a_k \mapsto 
  F(a_1 \delta_e \otimes  \ldots \otimes a_k \delta_e) \big),
\end{equation}
 where $a_i \delta_e$ denotes the smooth function on $\Gamma \times M$
 which coincides with $a_i$ on the unit space $u (M) $ and vanishes elsewhere.
\end{proposition}
\begin{proof}
  By the Eilenberg-Zilber Theorem one has 
  \begin{displaymath}
   H^\bullet  \operatorname{diag} (C^{\bullet,\bullet}_\Gamma) = H^\bullet 
   (\operatorname{Tot} C^{\bullet,\bullet}_\Gamma ).
  \end{displaymath}
  Moreover, there is a spectral sequence
  \begin{displaymath}
    E_1^{m,n} = H^n_{\text{h}}  ( C^{\bullet,m}_\Gamma ) , \quad
    E_2^{m,n} = H^m_{\text{v}} H^n_{\text{h}} ( C^{\bullet,\bullet}_\Gamma ) 
   \quad \Rightarrow \quad  H^{m+n} 
    \operatorname{diag} (C^{\bullet,\bullet}_\Gamma) .
  \end{displaymath}
  Now recall that the group cohomology of a finite group vanishes in degrees 
  $\geq 1$.
  Using the $\Gamma$-action from Lemma \ref{gammaactionhom} we thus obtain the 
  following chain of natural isomorphisms:
  \begin{equation}
    H^n (A\rtimes \Gamma, A \rtimes \Gamma) \cong H^n 
    \operatorname{diag} (C^{\bullet,\bullet}_\Gamma)
    \cong H^0_{\text{v}} H^n_{\text{h}} ( C^{\bullet,\bullet}_\Gamma ) 
    = \big( H^\bullet (A , A \rtimes \Gamma ) \big)^\Gamma  .
  \end{equation}  
  This proves the first claim; the second is a direct consequence of 
  Prop.~\ref{propsimpmap} and the spectral sequence argument leading 
  to Eq.~(\ref{eqhhiso}).
\end{proof}

\begin{corollary}
  Let $\Gamma, M, A$ as above and 
  $A_\text{c}:= \mathcal C^\infty_\text{c} (M)$ 
  the algebra of smooth functions with compact support on $M$. 
  Then there exists a commutative diagram of canonical isomorphisms:
  \begin{displaymath}
\begin{diagram}
\node{HH^k(A\rtimes \Gamma, A \rtimes \Gamma)}\arrow{s,r,T}{}\arrow{e,t,T}{}
\node{HH^k(A_\text{c}\rtimes \Gamma, A_\text{c} \rtimes \Gamma)}
\arrow{s,l,T}{}\\
\node{HH^k (A, A\rtimes \Gamma)^\Gamma}\arrow{e,t,T}{}
\node{HH^k (A_\text{c}, A_\text{c}\rtimes \Gamma)^\Gamma.}
\end{diagram} 
\end{displaymath}
\end{corollary}
\begin{proof}
  Using the identification $\check{\hspace{1em}}$ of the first step one checks 
  that the following morphisms are chain maps:
  \begin{displaymath}
  \begin{split}
    C^k (A \rtimes \Gamma, A \rtimes \Gamma) & \rightarrow 
    \Hom \big( (A_\text{c} \rtimes \Gamma)^{\topotimes k} , 
    A \rtimes \Gamma \big)  \cong 
    C^k (A_\text{c} \rtimes \Gamma, A_\text{c} \rtimes \Gamma) \\
    F & \mapsto \big( (A_\text{c} \rtimes \Gamma)^{\topotimes k} 
    \ni a_1 \otimes \ldots \otimes a_k \mapsto
    F ( a_1 \otimes \ldots \otimes a_k ) \in A \rtimes \Gamma \big) ,
  \end{split}
  \end{displaymath}
  and 
  \begin{displaymath}
  \begin{split}
    C^k (A , A \rtimes \Gamma) & \rightarrow 
    \Hom \big( A_\text{c}^{\topotimes k} , A \rtimes \Gamma \big)  \cong 
    C^k (A_\text{c} , A_\text{c} \rtimes \Gamma) \\
    F & \mapsto \big( A_\text{c}^{\topotimes k} 
    \ni f_1 \otimes \ldots \otimes f_k \mapsto
    F ( f_1 \otimes \ldots \otimes f_k ) \in A \rtimes \Gamma \big) .
  \end{split}
  \end{displaymath}
  With the help of the localization maps $\Psi^{\bullet,\varepsilon}$
  (associated to a complete $\Gamma$-invariant metric $d$ 
  on $M$) and an appropriate invariant smooth partition of unity on $M$
  one can construct quasi-inverses to the chain maps
  $C^\bullet (A \rtimes \Gamma, A \rtimes \Gamma)\rightarrow 
  C^\bullet (A_\text{c} \rtimes \Gamma, A_\text{c} \rtimes \Gamma)$ and
  $C^\bullet (A , A \rtimes \Gamma)  \rightarrow 
  C^\bullet (A_\text{c} , A_\text{c} \rtimes \Gamma) $.
  Thus, the two horizontal arrows in the above diagram are 
  isomorphisms. The left vertical arrow is an isomorphism by the preceding
  proposition, hence the induced right vertical arrow has to be an isomorphism
  as well.
\end{proof}

{\bf Step 4.}
According to Prop.~\ref{propisohcoh}, it suffices to compute  the 
(invariant part of the)
cohomology of the cochain complex  $C^\bullet (A, A\rtimes \Gamma)$, if
$G$ is a translation groupoid $\Gamma \ltimes M$.
To this end we specialize the situation further and assume that $M$ is an
open $\Gamma$-invariant neighborhood of the origin of some finite dimensional
linear $\Gamma$-representation space $V$. 
We choose a $\Gamma$-invariant scalar product on $V$ and orthonormal linear 
coordinates $x_1, \ldots,x_n$ of $V$ such that $x_1, \ldots , x_{l_\gamma}$ 
span the fixed point space $V^\gamma$ and $x_{l_\gamma +1}, \ldots , x_n$ 
span $W$, the subspace orthogonal to $V^\gamma$.
We assume further that $M$ has the form $M^\gamma \times N$ with 
$M^\gamma$ an open ball in $V^\gamma$ and $N$ an open ball in $W$.

For the computation of $H^\bullet (A, A\rtimes \Gamma)$ we will use the 
(topologically projective) resolution  of $A$ given by the Koszul 
complex $(K_\bullet, \partial)$ associated to the regular sequence
$(x_1 \otimes \id - \id \otimes x_1, 
\ldots, x_n \otimes \id - \id \otimes x_n)$ in $A\hat\otimes A$. 
More precisely, the resolution of $A$ by $K_\bullet$ has the form
\begin{displaymath}
   0 \longrightarrow A\hat\otimes A \otimes \Lambda ^n V^*
  \overset{\partial}{\longrightarrow} \cdots
  \overset{\partial}{\longrightarrow}
  A\hat\otimes A \otimes \Lambda ^k V^*
  \overset{\partial}{\longrightarrow} \cdots
  \overset{\partial}{\longrightarrow} 
  A \hat\otimes A \overset{m}{\longrightarrow A} \longrightarrow 0
\end{displaymath}
with differential 
$\partial : A\hat\otimes A \otimes \Lambda^k V^* \rightarrow
 A\hat\otimes A \otimes \Lambda^{k-1} V^*$ given by 
\begin{displaymath}
\begin{split}
  f_1\otimes f_2 \otimes & \, d x_{i_1} \wedge \ldots \wedge d x_{i_k} 
  \mapsto \\
  & \sum_{j=1}^k (-1)^j 
  ( x_{i_j} f_1\otimes f_2 - f_1 \otimes x_{i_j} f_2 ) \otimes 
  d x_{i_1} \wedge \ldots \wedge \widehat{d x_{i_j}} 
  \wedge \ldots \wedge d x_{i_k}. 
\end{split}
\end{displaymath} 
Let us provide another description of the Koszul complex
$(K_\bullet, \partial)$.
Denote by $E_k$ the pull-back bundle 
$\operatorname{pr}_2^* (\Lambda^k T^* M )$, where $\Lambda^k T^* M$
is the exterior product of the cotangent bundle of $M$, and
$\operatorname{pr}_2 : M \times M \rightarrow M$ is the projection on 
the second coordinate. Then the vector field 
\begin{equation}
  \xi : M \times M \rightarrow V , \quad (p,q) \mapsto \xi (p,q) = 
  \sum_{i=1}^n 
  \big( x_i (p)- x_i (q)\big) \, \frac{\partial}{\partial x_i},
\end{equation} 
comprises a section of $E^*_1$ which does not vanish outside the diagonal.
Moreover, $K_k$ can be naturally identified  with the sectional space
$\Gamma^\infty  (E_k)$, and $\partial$ is the insertion of the vector field
$\xi$. 

The cohomology $H^\bullet (A,A\rtimes \Gamma)$ now is given  as the direct sum
over the elements $\gamma \in \Gamma$ of the cohomologies of the cochain 
complexes $\big( \Hom (K_\bullet , A_\gamma) ,\partial^* \big)$, where 
$A_\gamma$ coincides with $A$ as a Fr\'echet space and carries the 
following $A$-bimodule structure: 
\begin{equation}
\label{Abimodulestr}
  (f_1 * a * f_2) (p) = f_1 (\gamma p) \, a(p) \, f_2 (p)
  \quad \text{for all $p\in M$, $a\in A_\gamma$, $f_1,f_2 \in A$}.
\end{equation}
This entails immediately that for every natural $k$ there 
is a canonical isomorphism
 \begin{displaymath} 
   \eta_k : \Gamma^\infty (\Lambda^k T M) \rightarrow 
   \Hom_{A - A} (K_k , A_\gamma), 
   \quad \tau \mapsto \eta (\tau ), 
 \end{displaymath} 
which is uniquely determined by the relation 
 \begin{displaymath}  
   \eta (\tau) (\omega) = \langle \Delta_\gamma^* \omega , \tau \rangle
   \quad \text{for all $\omega \in \Gamma^\infty (E_k)$}. 
 \end{displaymath} 
Hereby, $\langle -, -\rangle : \Omega^k (M) \times \Lambda^k T M 
\rightarrow \mathcal C^\infty (M)$ denotes the canonical  
fiberwise pairing, 
$\Delta_\gamma : M \rightarrow M \times M$ is the embedding  
$p \mapsto (\gamma p, p)$, and 
$\Delta_\gamma^* \omega$ is defined by
$\langle \Delta_\gamma^* \omega  (p),v\rangle = \langle \omega 
(\Delta_\gamma (p)),v \rangle$ for every 
$p\in M$ and $v \in T_pM \cong V$. 
Clearly, $\eta_k$ is injective. Let us show that it is surjective as well.
Let $F$ be a continuous $A$-bimodule map from $\Gamma^\infty (E_k)$ to 
$A_\gamma$ 
and define for all multiindices $1\leq i_1< \ldots < i_k \leq n$ coefficients
$\tau_{i_1,\ldots,i_k}$ by $\tau_{i_1, \ldots , i_k} := 
  F \big(\operatorname{pr}_2^* (dx_{i_1} \wedge \ldots \wedge dx_{i_1})\big)$.
Then $\eta$ 
maps the multivectorfield
\begin{displaymath}
  \tau := \sum_{i_1< \ldots < i_k} \, \tau_{i_1, \ldots , i_k} 
  \frac{\partial}{\partial x_{i_1}} \wedge \ldots \wedge 
  \frac{\partial}{\partial x_{i_n}}
\end{displaymath}
to $F$, hence $\eta$ is surjective.

Now let $\kappa$ be the vector field on $M$ defined by 
\begin{equation}
\label{defkappa}
  \kappa (p) =  \xi (\gamma p ,p ) = \sum_{i = l_\gamma +1}^n 
  \big( x_i (\gamma p)- x_i (p)\big) \, \frac{\partial}{\partial x_i} .
\end{equation}
Under the isomorphism $\eta$, the cohomological differential  
$\partial^*$ corresponds to the operation $\kappa \wedge -$. To check this, let
$\omega \in \Gamma^\infty (E_{k+1})$ and compute:
\begin{displaymath}
  \big( \partial^* \eta (\tau) \big) \, (\omega) = \eta (\tau) \, 
  (i_\xi \omega) =
  \langle \Delta^*_\gamma i_\xi \omega , \tau \rangle = 
  \langle i_\kappa^* \Delta^*_\gamma \omega , \tau \rangle = 
  \langle \Delta^*_\gamma \omega , \kappa \wedge \tau \rangle
  = \eta (\kappa \wedge \tau) (\omega) ,
\end{displaymath}   
which proves the claim. Hence it remains to determine the cohomology of 
the cochain complex
\begin{equation}
\label{dualkoszul}
  ( \Gamma^\infty (\Lambda^\bullet T M), \kappa \wedge - ).
\end{equation}
But this complex is a dual Koszul complex. To compute its cohomology observe
first that  the decomposition $V = V^\gamma \oplus W$
induces a decomposition of the alternating multivector fields on $M$ as 
follows:
\begin{equation}
\label{decaltvector}
  \Gamma^\infty \Lambda^k ( TM  ) \cong
  \bigoplus_{p=0}^k \Lambda^p ( V^\gamma ) \otimes \Gamma^\infty 
  ( \Lambda^{k-p} \operatorname{pr}_{N}^* TN )), 
\end{equation}
where $\operatorname{pr}_{N}:M \rightarrow N$ is the projection onto $N$
along $M^\gamma$.
Under this isomorphism, the differential $\kappa \wedge -$ acts only on 
the second components.
Hence one can interpret the cohomology  of (\ref{dualkoszul}) as the total 
cohomology of the
double complex $D^{p,q} = \Lambda^p ( V^\gamma ) \otimes \Gamma^\infty 
 ( \Lambda^q \operatorname{pr}_{N}^* TN ))$, which has $0$-differential
in the $p$-direction and differential $\kappa \wedge -$ in $q$-direction.
Since $\operatorname{pr}_{N}^* TN$ is a trivial vector bundle with fiber 
dimension $n-l_\gamma$, the sectional space 
$\Gamma^\infty (\Lambda^q \operatorname{pr}_{N}^* TN )$ is isomorphic to 
$\mathcal C^\infty (M) \otimes \Lambda^q \R^{n-l_\gamma}$. Together with
Eq.~(\ref{defkappa}) this implies that 
$\big( \Gamma^\infty (\Lambda^q \operatorname{pr}_{N}^* TN ),
\kappa\wedge -\big)$
is the dual Koszul complex of the algebra $\mathcal C^\infty (M)$ associated
to the regular sequence 
\begin{equation}
\label {regseq}
  (\gamma^{-1} x_{l_\gamma + 1} - x_{l_\gamma + 1} ,\ldots, 
  \gamma^{-1} x_n - x_n).
\end{equation}
The cohomology of this dual Koszul complex is well-known 
(cf.~\cite[Sec.~17.2]{eisenbud}).
It does not vanish only for $q= n-l_\gamma$, 
where it is given as the quotient of $\mathcal C^\infty (M)$ by the (closed)
ideal generated by the regular sequence (\ref{regseq}), i.e., by the algebra
$\mathcal C^\infty (M^\gamma)$. 
Using the spectral sequence of the double complex $D^{p,q}$
one then concludes that    
\begin{equation}
\label{cohomdualkoszul}
  H^k ( A, A_\gamma) \cong 
  H^k ( \Gamma^\infty (\Lambda^\bullet T M), \kappa \wedge - )
  \cong  \Gamma^\infty (M^\gamma , \Lambda^{k-n +l_\gamma} TM^\gamma ).
\end{equation}
Using a standard localization argument for Hochschild cohomology 
(see Remark \ref{remcommloc}) one now infers from this equation and 
Prop.~\ref{propisohcoh} the following result.
\begin{proposition}
  Let $\Gamma$ be a finite group acting on a smooth manifold $M$. Then 
  the Hochschild cohomology $H^k (A \rtimes \Gamma ,A \rtimes \Gamma)$
  can be naturally identified as follows with spaces of invariant 
  multivectorfields:
  \begin{equation}
  \begin{split}
     H^k & \, (A \rtimes \Gamma ,A \rtimes \Gamma ) \cong 
     H^k (A ,A \rtimes \Gamma )^\Gamma
     = \\
     & = \bigoplus_{<\gamma >  \in \operatorname{Conj} (\Gamma)} \: 
     \bigoplus_{M^\gamma_\alpha  \in \operatorname{Comp} (M^\gamma)}
     \Gamma^\infty \big( M^\gamma_\alpha , 
     \Lambda^{k-\dim M + \dim M^\gamma_\alpha} 
     TM^\gamma_\alpha \big)^{Z(\gamma)},
  \end{split}
  \end{equation}
  where $\operatorname{Conj} (\Gamma)$ is the set of conjugacy classes of 
  $\Gamma$, $\operatorname{Comp} (M^\gamma)$ the set of connected components
  of $M^\gamma$, and $Z(\gamma) \subset \Gamma$ the centralizer of 
  $\gamma$ in $\Gamma$. 
\end{proposition}

{\bf Step 5.}
From now on we consider again the general case of a proper \'etale Lie 
groupoid and use all the previous results to prove the following main theorem.
\begin{theorem}
  Let $G$ be a proper \'etale Lie groupoid. 
  Then the Hochschild cohomology of the convolution algebra 
  $\mathcal A \rtimes G$ with
  values in $\mathcal A \rtimes G$ is naturally given as follows:
  \begin{equation}
     H^k (\mathcal A \rtimes G ,\mathcal A \rtimes G) \cong 
     \bigoplus_{\mathcal O \in \operatorname{Sec} (G)} \,
     \Gamma^\infty_\text{\rm \tiny inv} \big(  
     \Lambda^{k  -\ell (\mathcal O)} T\mathcal O 
     \big), 
  \end{equation}
  where the sum is taken over the sectors of $G$.
\end{theorem}  
\begin{proof}
 Consider the complex $\mathcal K^\bullet$ of sheaves on the orbit space $X$
 constructed in Step 1., and define a second sheaf complex 
 $\mathcal H^\bullet$ on $X$ (with differential the zero map) as follows:
 \begin{displaymath}
   \mathcal H^k (U) := 
   \bigoplus_{\mathcal O \in \operatorname{Sec} (G)} \,
   \Gamma^\infty_\text{\rm \tiny inv} \big( \mathcal O \cap \Lambda U_0 , 
   \Lambda^{k  -\ell (\mathcal O)} T\mathcal O
   \big),
 \end{displaymath}
 where $U$ runs through the open subsets of $X$ and 
 $\Lambda U_0 := (\pi \circ \beta \circ s)^{-1} (U)$. 
 Observe now that both sheaf complexes $\mathcal K^\bullet$ and 
 $\mathcal H^\bullet$ are fine, since the sheaf of smooth functions 
 on the orbifold $X$ is fine. Moreover, note that the global section space
 of the cohomology sheaf of 
 $\mathcal K_\bullet$ is the cohomology we want to compute and that
 $\mathcal H^\bullet (X)$ is the graded vector space we claim the
 cohomology to coincide with.
 Hence, if one can construct 
 a morphism of sheaf complexes 
 $\Xi^\bullet : \mathcal K^\bullet \rightarrow \mathcal H^\bullet$ 
 which locally is a quasi-isomorphism, the claim is proved by 
 \cite[Chap.~6, Sec.~8, Thm.~9]{spanier}. 
 Thus it remains to construct $\Xi$ and prove that, locally, $\Xi$ is a 
 quasi-isomorphism. Before we come to the details of the construction we need 
 two lemmas.
\begin{lemma}
  Assume $U\subset X$ to be open, let $U_0:=  \pi^{-1} (U)$ and
  $U_1:= ( \pi \circ s)^{-1} (U)$. Denote by 
  $G_{|U_1}$ the restriction of the groupoid $G$ to $U_1$
  and let $\mathcal C_{U_0}^\infty$ be the $G_{|U_1}$-sheaf of smooth
  functions on $U_0$.
  Then the embedding 
  $\mathcal C_{U_0}^\infty  \rtimes G_{|U_1} \hookrightarrow 
  \mathcal A \rtimes G$ 
  induces a quasi-isomorphism
  \begin{equation}
  \label{eqshmor}
     \mathcal K^\bullet (U) \rightarrow C^\bullet 
     \big(\mathcal C_{U_0}^\infty \rtimes G_{|U_1}, 
     \mathcal C^\infty (U_1) \big).
  \end{equation}
\end{lemma}
\begin{lemproof}
   Note first that $C^\bullet 
   \big(\mathcal C_\text{c}^\infty (U_0) \rtimes G_{|U_1},
   \mathcal C^\infty (U_1) \big)$
   is the global section space of the sheaf complex $\mathcal K_U^\bullet$,
   which for $V\subset U$ open and $k\in \N$  has section space 
   \begin{displaymath}
      \mathcal K_U^k (V) := \Hom \big( 
      (\mathcal C_{U_0}^\infty  \rtimes G_{|U_1})^{\hat\otimes k}, 
      \mathcal C^\infty (V_1) \big) , \quad V_1 := ( \pi \circ s)^{-1} (V),
   \end{displaymath}
   and which has the Hochschild coboundary as its differential. Note also that
   $\mathcal K_U^\bullet$ and the restriction $\mathcal K_{|U}^\bullet$ of 
   $\mathcal K^\bullet$ to $U$ are both complexes of fine sheaves, since 
   the sheaf of smooth functions on $X$ is fine.
   If we can now show that the natural morphism of sheaf complexes 
   $\mathcal K_{|U}^\bullet \rightarrow \mathcal K_U^\bullet$  is locally a 
   quasi-isomorphism the claim is proved by \cite[Chap.~6, Sec.~8, Thm.~9]{spanier}.
    
   To verify this it suffices to check that
   $\mathcal K_{|U}^\bullet (V) \rightarrow  \mathcal K_U^\bullet (V)$ is 
   a quasi-isomorphism 
   for every relatively compact connected open subset $V \subset U$. 
   To this end choose a complete metric
   $d_X$ on $X$ and a complete metric $d$ on $G_0$ (for which $d^2$ is smooth)
   such that $d(x,y) \geq d_X(\pi (x),\pi (y))$ for all $x,y \in G_0$.
   By the assumptions on $V$, there exists an $\varepsilon >0$ such that 
   \begin{equation}
   \label{ineq}
    d(t(g),t(h))\geq d_X(\pi (t(g) , \pi (t(h)) > \varepsilon \quad 
    \text{for all $g\in V_1$ and $h \in  G_1 \setminus U_1$}.
   \end{equation}
   Moreover, since the preconditions of Step.~2 are satisfied,  
   we have the localization functions $\Psi_{k,\varepsilon}$ at our disposal.
   With their help define now for every $N\in  \N^*$ chain maps 
   $\Theta_{V}^{\leq N}:\mathcal K^{\leq N}_U (V)\rightarrow\mathcal 
   K^{\leq N} (V)$ between the cut-off chain complexes as follows:
   \begin{displaymath}
     \Theta_{V}^k (F) (a_1\otimes \ldots \otimes a_k) \, (g_0) =
     F \big( \Psi_{k,\varepsilon /N}
     (g_0^{-1}, -,\ldots,-) \cdot (a_1 \otimes \ldots \otimes a_k) \big) (g_0),
   \end{displaymath}
   where $k\leq N$, $F \in \mathcal K_U^k (V)$, 
   $a_1,\ldots,a_k\in \mathcal A \rtimes G$ and  $g_0\in V_1$.
   By (\ref{ineq}) one concludes that 
   \begin{displaymath}
     \Psi_{k,\varepsilon/N} 
     (g_0^{-1},g_1,\ldots ,g_k) \cdot a_1(g_1)\cdot\ldots\cdot a_k (g_k) =0,
   \end{displaymath}
   if  $g_0\in V_1$ and $g_1,\ldots ,g_k\in G_1$ with some 
   $g_i \in G\setminus U_1$, hence $\Theta_{V}^k (F)$ is well-defined 
   indeed for $k\leq N$.
   Prop.~\ref{prophomtploc} now entails that 
   $\mathcal K_{|U}^\bullet (V) \rightarrow  \mathcal K_U^\bullet (V)$ is a 
   quasi-isomorphism in degrees $k<N$. Since $N$ was arbitrary, it is a 
   quasi-isomorphism in all degrees, and the claim follows. 
\end{lemproof}
  Since $G$ is a proper \'etale Lie  groupoid, there exists for every point 
  $\tilde x\in G_0$ an open contractible neighborhood 
  $M_{\tilde x}\subset G_0$,
  a smooth action of the isotropy $G_{\tilde x}$ on  
  $M_{\tilde x}\subset G_0$ and a monomorphism of groupoids 
  \begin{displaymath}
    \iota_{\tilde x} : G_{\tilde x} \ltimes M_{\tilde x} \hookrightarrow 
    G
  \end{displaymath}  
  which induces a Morita equivalence of groupoids from 
  $G_{\tilde x} \ltimes M_{\tilde x}$ to the restricted groupoid 
  $G_{|U_{\tilde x,1}}$, where $U_{\tilde x} := \pi (M_{\tilde x})$
  and $U_{\tilde x,1} = (\pi\circ s)^{-1} (U_{\tilde x})$. 
\begin{lemma}
  The Morita equivalence $\iota_{\tilde x}$ gives rise to a 
  quasi-isomorphism
  \begin{equation}
    \iota^* : \: \mathcal K^\bullet_{U_{\tilde x}} ({U_{\tilde x}})
    \rightarrow C^\bullet \big(\mathcal C^\infty_\text{c} ( M_{\tilde x}) 
    \rtimes G_{\tilde x} , 
    \mathcal C^\infty (M_{\tilde x})\rtimes G_{\tilde x}\big)
  \end{equation} 
  which associates to every $F \in \mathcal K^k_{U_{\tilde x}} ( U_{\tilde x})$
  the cochain
  \begin{displaymath}
  \begin{split}
    \iota^* (F) : \: & \big(\mathcal C^\infty_\text{c} ( M_{\tilde x})
    \rtimes G_{\tilde x}\big)^{\hat \otimes k}
    \rightarrow \mathcal C^\infty (M_{\tilde x})\rtimes G_{\tilde x},
    \\  &  a_1 \otimes \ldots \otimes a_k \mapsto
    F (\iota_* a_1 \otimes \ldots \otimes \iota_* a_k)_{|M_{\tilde x}},
  \end{split}
  \end{displaymath}
  where we have put for $a \in \mathcal C^\infty_\text{c} ( M_{\tilde x})
    \rtimes G_{\tilde x}$ and $g\in U_{\tilde x,1}$
  \begin{displaymath}
    \iota_* (a) (g) = 
    \begin{cases}
       a (\iota_{\tilde x,1}^{-1} (g)), & \text{if $g\in \operatorname{im}
       \iota_{\tilde x}$,}\\
       0, & \text{else}.
    \end{cases}
  \end{displaymath}
\end{lemma} 
\begin{lemproof}
  For the proof of the claim we use the language of Hilsum-Skandalis maps
  (i.e~biprincipal bundles) and their associated Morita bimodules as
  explained in \cite{mrcun}. As shown in \cite[Sec.~1]{mrcun}, the Morita
  equivalence 
  $\iota: G_{\tilde x} \ltimes M_{\tilde x} \hookrightarrow G_{U_{\tilde x}} $
  induces a principal $G_{\tilde x}$-$G_{|U_{\tilde x,1}}$-bibundle 
  $\langle \iota \rangle$ as follows:
  \begin{displaymath}
  \begin{split}
    \langle \iota \rangle & =\{ (g,p) \in  G_1\times M_{\tilde x} \mid
    s(g) =p \} \cong \{ g \in G_1 \mid s(g) \in M_{\tilde x} \},\\
    & G_1 \times_{G_0} \langle \iota \rangle \rightarrow \langle \iota 
    \rangle , \quad (g' ,g) \mapsto g' g , \\
    & \langle \iota \rangle \times_{M_{\tilde x}} 
    (G_{\tilde x} \ltimes M_{\tilde x} ) \rightarrow \langle \iota \rangle , 
    \quad (g,\gamma,p) \mapsto g \cdot \iota (\gamma , p) .
  \end{split}
  \end{displaymath}  
  Then, by \cite[Sec.~2]{mrcun}, the locally convex topological vector space
  $\mathcal C_\text{c}^\infty (\langle \iota \rangle)$ carries the structure
  of a
  $\big( \mathcal C_\text{c}^\infty (U_{x,0}) \rtimes G_{|U_{x,1}}\big)$-$\big(
  \mathcal C_\text{c}^\infty (M_{\tilde x}) \rtimes G_{\tilde x}
  \big)$-bimodule 
  and forms a Morita equivalence between these two algebras.
  Now, the chain map $\iota^*$ is induced by this Morita equivalence,
  as one checks by an immediate but somewhat tedious computation
  (see \cite[Sec.~1.2]{loday}). Hence, $\iota^*$ is a quasi-isomorphism.
\end{lemproof} 

 Now we come back to the construction of the morphism of sheaf complexes 
 $\Xi$. Let $U\subset X$ be open, and choose for every $x\in U$ an open 
 neighborhood $U_x \subset U$, 
 a point $\tilde x \in \pi^{-1} (x) $ together with an open neighborhood 
 $M_{\tilde x} \subset G_0$, an action of the isotropy group $G_{\tilde x}$ 
 on $M_{\tilde x}$ and, finally, a Morita equivalence of Lie groupoids 
 \begin{displaymath}
  \iota_{\tilde x} : G_{\tilde x} \ltimes M_{\tilde x} \rightarrow 
  G_{|U_{x,1}}. 
 \end{displaymath}
 Then one has for every one of the $U_x$ a sequence of natural chain maps:
 \begin{equation}
 \begin{split}
   \mathcal K^\bullet & (U)  \longrightarrow \mathcal K^\bullet (U_x) 
   \overset{\text{\tiny q.i.}}{\longrightarrow} 
   \mathcal K^\bullet_{U_x} (U_x) \overset{\iota_{\tilde x}^*}\longrightarrow 
   C^\bullet \big( \mathcal C^\infty_{\text{c}} (M_{\tilde x}
   \rtimes G_{\tilde x} ),
   C^\infty (M_{\tilde x} ) \rtimes G_{\tilde x}  \big) \longrightarrow \\ 
   & \longrightarrow
   C^\bullet \big( \mathcal C^\infty_{\text{c}} (M_{\tilde x} ),
   C^\infty (M_{\tilde x} ) \rtimes G_{\tilde x}  \big) \longrightarrow 
   \Gamma^\infty (\Lambda^\bullet TM_{\tilde x} ) \longrightarrow \\
   & \longrightarrow \bigoplus_{\gamma \in G_{\tilde x}}  \:
   \bigoplus_{M_{\tilde x,\alpha}^\gamma \in \operatorname{Comp} 
   (M_{\tilde x}^\gamma)}
   \Gamma^\infty (\Lambda^{\bullet-\dim G + \dim M_{\tilde x,\alpha}^\gamma} 
   TM_{\tilde x,\alpha}^\gamma )
   \longrightarrow \\
   & \longrightarrow \bigoplus_{<\gamma> \in \operatorname{Conj} G_{\tilde x}}  \:
   \bigoplus_{M_{\tilde x,\alpha}^\gamma \in \operatorname{Comp} 
   (M_{\tilde x}^\gamma)}
   \Gamma^\infty (\Lambda^{\bullet-\dim G + \dim M_{\tilde x,\alpha}^\gamma} 
   TM_{\tilde x,\alpha}^\gamma )^{Z (\gamma)}
   \cong \mathcal H^\bullet (U_x), 
 \end{split}
 \end{equation} 
 where the arrow in the last line is the projection onto the invariant part 
 obtained by averaging over $G_{\tilde x}$.
 By naturality of all constructions involved one checks that for all $x,y\in U$ 
 the following diagram commutes:
 \begin{displaymath}
 \begin{diagram}
  \node{\mathcal K^\bullet (U) }\arrow{s,r,T}{}\arrow{e,t,T}{ } 
  \node{\mathcal H^\bullet (U_y) }
  \arrow{s,l,T}{}\\
  \node{\mathcal H^\bullet (U_x) }\arrow{e,t,T}{}\node{\mathcal H^\bullet 
  (U_x \cap U_y).}
 \end{diagram} 
 \end{displaymath}
 Hence, by the sheaf property of $\mathcal H^\bullet$ one can glue together 
 these maps to a 
 chain map $\Xi (U): \mathcal K^\bullet (U) \rightarrow \mathcal H^\bullet (U)$.
 By construction, the $\Xi (U)$ commute with restriction maps, hence one obtains a
 morphism of sheaf complexes 
 $\Xi: \mathcal K^\bullet \rightarrow \mathcal H^\bullet$.
 By the above lemmas and Steps 3.~and 4.~it is clear that for every $x\in U$ 
 the chain map 
 $\Xi (U_x):\mathcal K^\bullet (U_x)\rightarrow\mathcal H^\bullet (U_x)$ has 
 to be a quasi-isomorphism.  This finishes the last part of the proof and 
 thus entails the claim.
\end{proof}
\section{Noncommutative Poisson homology}
\label{ncph}
This section is divided into two parts. In the first part, we
introduce a Poisson homology for a noncommutative Poisson algebra.
In the second part, we calculate this homology for the Poisson
algebra constructed from a Poisson structure on a proper \'etale Lie
groupoid (recall Sec.~\ref{sec:intr-ncq} for definitions).
\subsection{Poisson homology} 
In \cite{byi:poisson}, Brylinski defined Poisson
homology on a Poisson manifold $(M,\Pi)$ as the homology of the complex 
$b_\Pi :\Omega^\bullet(P)\to \Omega ^{\bullet -1}(P)$, where 
\begin{equation}
\label{poissondiff}
\begin{split}
  b_\Pi (f_0 \, df_1 & \wedge df_2\wedge \cdots \wedge df_k)= 
  \sum _{j=1}^k (-1)^{j-1}d\{ f_0, f_j\} \, df_1\wedge\cdots\wedge
  \widehat {df}_j
  \wedge \cdots \wedge df_k + \\
  & + \, \sum _{i<j} (-1)^{i+j-1}f_0 \, d\{f_i, f_j\}\wedge df_1\cdots \wedge
 \widehat {df}_i \wedge \cdots \wedge\widehat {df}_j \cdots \wedge df_k.
\end{split}
\end{equation}
In the following, we define a noncommutative analog of this
Poisson homology.
Like in the manifold case, we start from a Poisson structure
$[\Pi] \in HH^2 (A,A)$ of the algebra $A$ and let it act on 
$HH_\bullet (A)$, the noncommutative analog of differential forms. 
Before we introduce the precise definition, we first recall a
well-known action of the Hochschild cochains on  Hochschild chains.

\begin{definition} 
\label{dfn:L-infty-rep} 
  For an element $\varphi\in C^k(A, A)$, define 
  $d_{\varphi}: C_{n}(A)\to C_{n-k+1}(A)$ by
  $d_{\varphi}  (a_0\otimes a_1\otimes \cdots \otimes a_n):= $
  \begin{displaymath}
  \begin{split}
    & \sum _{i=0}
    ^{n-k+1}(-1)^{i(k-1)} a_0 \otimes \cdots\otimes \varphi
   (a_i \otimes \cdots \otimes a_{i+k-1})
    \otimes\cdots\otimes a_{n} \: +\\
    & +\sum_{i=2} ^{k}(-1)^{(n-k+i)(k-i+1)}\varphi(a_{n-k+i} \otimes 
    \cdots \otimes a_{n} \otimes  a_0 \otimes \cdots \otimes a_{i-2})
    \otimes \cdots \otimes a_{n-k+i-1}.
  \end{split}
\end{displaymath}
\end{definition}

\begin{theorem}
\label{prop:action} 
For $\phi \in C^k(A, A),\ \psi \in C^l(A, A)$, we have
\begin{displaymath}
  d_{\phi}\circ d_{\psi}-(-1)^{(k-1)(l-1)}d_{\psi}\circ
  d_{\phi}=d_{[\phi, \psi]}.
\end{displaymath}
\end{theorem}

\begin{proof}
We prove this property on $a_0\otimes a_1\otimes \cdots \otimes a_n$. \\[2mm]
(1) $ d_{\phi}\circ d_{\psi}(a_0\otimes a_1\otimes \cdots \otimes a_n) =$
{\small
\begin{displaymath}
\begin{array}{cl}
 = & \sum _{i=0}^{n-l+1}(-1)^{i(l-1)}d_{\phi}
   \big( a_0 \otimes \cdots\otimes \psi(a_i\otimes a_{i+1}\otimes \cdots 
   \otimes a_{i+l-1})\otimes\cdots \big) \: + \\[1mm]
   & + \: \sum_{i=2}^l (-1)^{(n-l+i)(l-i+1)} \\
   & \hspace{1.5em} d_{\phi}\big( \psi(a_{n-l+i}\otimes \cdots \otimes a_{n}\otimes 
   a_0\otimes \cdots \otimes a_{i-2})\otimes a_{i-1}\otimes \cdots \otimes a_{n-l+i-1}\big)\\[1mm]
 = & \sum _{i=0} ^{n-l+1}(-1)^{i(l-1)}\Big( \sum_{j=0}^{i-k}(-1)^{j(k-1)} \\
   & \hspace{1.5em} \big( a_0\otimes \cdots \otimes \phi (a_j\otimes \cdots \otimes a_{j+k-1}) \otimes 
     \cdots \otimes \psi(a_i\otimes  \cdots\otimes a_{i+l-1})\otimes\cdots \big) \: +\\[1mm]
   & + \: \sum_{j=i-k+1} ^{i} (-1)^{j(k-1)} \\
   & \hspace{1.5em} \big( a_0 \otimes\cdots \otimes a_{j-1}\otimes 
   \phi(a_j\otimes \cdots \otimes  \psi(a_i\otimes \cdots \otimes a_{i+l-1}) \otimes \cdots \otimes a_{i+k+l-2})
   \otimes\cdots \big) \: +\\[1mm]
   & + \: \sum _{i=j+k+l-1} ^{n} (-1)^{(j-k+1)(l-1)} \\
   & \hspace{1.5em} \big( a_0\otimes \cdots \otimes \psi(a_j\otimes \cdots \otimes a_{j+l-1} )\otimes \cdots
   \otimes \phi(a_i\otimes \cdots \otimes a_{i+k-1})\otimes \cdots \big) \Big) \:  +\\[1mm]
   & + \: \sum_{i=2}^l (-1)^{(n-l+i)(l-i+1)} \\
   & \hspace{1.5em} \Big( \big( \phi(\psi(a_{n-l+i}\otimes \cdots \otimes a_{n} \otimes a_0 \otimes \cdots
   \otimes a_{i-2})\otimes a_{i-1}\otimes\cdots\otimes a_{i+k-2})\otimes a_{i+k-1} \otimes \cdots \big) \: + \\[1mm]
   & + \: \sum _{j=i-1}^{n-l+i-k}(-1)^{(i-j)(k-1)} \\
   & \hspace{1.5em} \big( \psi (a_{n-l+i}\otimes \cdots \otimes a_{n} \otimes a_0 \otimes \cdots \otimes a_{i-2})\otimes
   a_{i-1} \otimes \cdots \otimes \phi(a_j\otimes
   \cdots \otimes a_{j+k-1})\otimes \cdots \big) \: +\\[1mm]
   & + \: \sum_{j=n-l+i-k+1}^{n-l+i-1}(-1)^{(n-l+i-j)(j-i+1)} \\
   & \hspace{1.5em} \big( \phi (a_j\otimes \cdots \otimes \psi (a_{n-l+i}\otimes \cdots \otimes a_{i-2})\otimes 
   \cdots  \otimes a_{j+k+l-n-3})\otimes\cdots \big) \Big) .
\end{array}
\end{displaymath}}
(2) $ d_{\psi}\circ d_{\phi}(a_0\otimes a_1\otimes \cdots \otimes a_n) =$ \\
{\small
\begin{displaymath}
  \begin{array}{cl}
 = & \sum _{j=0} ^{n-k+1}(-1)^{j(k-1)}\Big( \sum_{i=0}^{j-l}(-1)^{i(l-1)} \\
   & \hspace{1.5em} \big( a_0\otimes \cdots \otimes \psi (a_i\otimes \cdots \otimes a_{i+l-1}) \otimes 
     \cdots \otimes \phi(a_j\otimes  \cdots\otimes a_{j+k-1})\otimes\cdots \big) \: +\\[1mm]
   & + \: \sum_{i=j-l+1} ^{j} (-1)^{i(l-1)} \\
   & \hspace{1.5em} \big( a_0 \otimes\cdots \otimes a_{i-1}\otimes 
   \psi(a_i\otimes \cdots \otimes  \phi(a_j\otimes \cdots \otimes a_{j+k-1}) \otimes \cdots \otimes a_{j+l+k-2})
   \otimes\cdots \big) \: +\\[1mm]
   & + \: \sum _{j=i+l+k-1} ^{n} (-1)^{(i-l+1)(k-1)} \\
   & \hspace{1.5em} \big( a_0\otimes \cdots \otimes \phi(a_i\otimes \cdots \otimes a_{i+k-1} )\otimes \cdots
   \otimes \psi(a_j\otimes \cdots \otimes a_{j+l-1})\otimes \cdots \big) \Big) \:  +\\[1mm]
   & + \: \sum_{j=2}^k (-1)^{(n-k+j)(k-j+1)} \\
   & \hspace{1.5em} \Big( \big( \psi(\phi(a_{n-k+j}\otimes \cdots \otimes a_{n} \otimes a_0 \otimes \cdots
   \otimes a_{j-2})\otimes a_{j-1}\otimes\cdots\otimes a_{j+l-2})\otimes a_{j+l-1} \otimes \cdots \big) \: + \\[1mm]
   & + \: \sum _{i=j-1}^{n-k+j-l}(-1)^{(j-i)(l-1)} \\
   & \hspace{1.5em} \big( \phi (a_{n-k+j}\otimes \cdots \otimes a_{n} \otimes a_0 \otimes \cdots \otimes a_{j-2})\otimes
   a_{j-1} \otimes \cdots \otimes \psi(a_i\otimes
   \cdots \otimes a_{i+l-1})\otimes \cdots \big) \: +\\[1mm]
   & + \: \sum_{i=n-k+j-l+1}^{n-k+j-1}(-1)^{(n-k+j-i)(i-j+1)} \\
   & \hspace{1.5em} \big( \psi (a_i\otimes \cdots \otimes \phi (a_{n-k+j}\otimes \cdots \otimes a_{j-2})\otimes 
   \cdots  \otimes a_{i+l+k-n-3})\otimes\cdots \big) \Big) .
\end{array}
\end{displaymath}}
(3) $d_{\phi}\circ d_{\psi}-(-1)^{(k-1)(l-1)}d_{\psi}\circ d_{\phi}$. \\
Note that there are two types of terms:
  \begin{enumerate}
  \item[(a)] \label{type1}
  $\big( \cdots \otimes \phi (\cdots )\otimes \cdots \otimes \psi(\cdots )
  \otimes\cdots \big) \: $ 
  or 
  $\: \big( \cdots \otimes \psi(\cdots)\otimes \cdots \otimes \phi( \cdots)
  \otimes\cdots \big)$.
  \item[(b)] \label{type2}
  $\big( \cdots \otimes \phi(\cdots \psi(\cdots))\otimes \cdots \big)\: $ or 
  $\: \big( \cdots \otimes \psi(\cdots \phi(\cdots))\otimes \cdots \big)$.
  \end{enumerate}
  For type (a), they appear in both $d_{\psi}\circ d_{\phi}$ and
  $d_{\phi}\circ d_{\psi}$, and differ by the sign
  $(-1)^{(k-1)(l-1)}$. Therefore, they get cancelled in the sum.
  For type (b), terms like 
  $$\big( \cdots \otimes \phi(\cdots \psi(\cdots)\cdots)\otimes \cdots \big)$$
  appear in $d_{\phi}\circ  d_{\psi}$, while terms like 
  $$\big( \cdots \otimes \psi(\cdots \phi(\cdots)\cdots)\otimes \cdots \big)$$ 
  appear in $d_{\psi}\circ d_{\phi}$. 
  They come out in pairs. It is straightforward to check that their signs 
  match with those of $d_{[\phi, \psi]}$.
  Therefore, we conclude that 
  $$d_{\phi}\circ d_{\psi}-(-1)^{(k-1)(l-1)}d_{\psi}\circ d_{\phi}=
  d_{[\phi,\psi]}.$$
\end{proof}

\begin{remark}
  Definition \ref{dfn:L-infty-rep} is a special case of a general theory 
  of Nest and Tsygan \cite{nt99} of operations on Hochschild and 
  cyclic complexes.
  Moreover, Theorem \ref{prop:action}  shows that by Definition
  \ref{dfn:L-infty-rep} one obtains an $L_{\infty}$-module structure for
  the differential graded Lie algebra $(C^\bullet(A, A), \beta, [\ ,\ ])$
  and also a Lie module structure  on $HH_\bullet(A)$
  for the super Lie algebra $(H^\bullet(A, A), [\ , \ ])$.
\end{remark}

On an associative algebra $A$, there is a natural $2$-cocycle $m$
associated to the multiplication defined by 
$m(a_1\otimes a_2)=a_1a_2$. It is easy to check that $d_m$ is the Hochschild
differential on $C_\bullet(A)$, so we will simply write $b$ instead of
$d_m$.

By taking $\phi=m$, $\psi=\Pi$, where $\Pi \in Z^2(A,A)$ is a representative 
of the Poisson structure $\pi$, Theorem \ref{prop:action} now gives
\begin{displaymath}
  b\circ d_{\Pi}+d_{\Pi}\circ b=d_{[m, \Pi]}=d_{\beta (\Pi)}=0.
\end{displaymath}
Therefore, $d_{\Pi}$ descends to homology and gives rise to a map
$HH_\bullet(A)\to HH_{\bullet-1}(A)$.
By taking $\phi=\psi=\Pi$ and using that $[\Pi,\Pi] = \beta (\Theta)$ for
some Hochschild 1-cochain $\Theta$, we get
\begin{displaymath}
  2d_{\Pi}\circ d_{\Pi}=d_{[\Pi, \Pi]}=d_{\beta (\Theta)}=b\circ
  d_{\Theta}+d_{\Theta}\circ b.
\end{displaymath}
This proves that $d_{\Pi} ^2=0$ in $ HH_\bullet(A)$.
Finally, we obtain for $\phi=m$, $\psi=\beta(\eta)$
\begin{displaymath}
  b\circ d_{\eta}-d_{\eta}\circ b=d_{[m,\eta]}  =d_{\beta (\eta)}=d_{\psi},
\end{displaymath}
which shows that two representatives of the Poisson structure
define the same homology in $ HH_\bullet(A)$. In other words this means
that noncommutative Poisson homology defined below depends only on the 
Poisson structure $[\Pi] \in H^2(A,A)$ and not on the particular choice of a 
representative $\Pi \in Z^2 (A,A)$.

\begin{definition}
\label{def:poisson} For a noncommutative Poisson structure $[\Pi]$
on an associative algebra $A$, its Poisson homology is defined as
the homology of the differential complex 
$d_{\Pi}:HH_\bullet(A) \to HH_{\bullet-1}(A)$, where $\Pi$ is a 
cocycle representing the Poisson structure.
\end{definition}

In the case of a Poisson manifold $M$, our definition of the
Poisson homology on the Hochschild homology of
$\mathcal C_\text{c}^{\infty}(M)$ is compatible with the one defined on
differential forms by Brylinski.

\begin{proposition}
\label{commutative} 
Let $M$ be a Poisson manifold with Poisson bivector $\Pi$,
and $A_\text{c}=\mathcal C_\text{c}^{\infty}(M)$ the algebra of 
compactly supported
smooth functions together with the noncommutative Poisson structure 
induced by $\Pi$. Then the following diagram commutes: 
\begin{displaymath}
\begin{diagram}
\node{HH_k(A_\text{c})}\arrow{s,l,T}{d_{\Pi}}\node{\Omega_\text{c}^k(M)}
\arrow{w,t,T}{\epsilon _k }\arrow{s,l,T}{2(k-1)!b_\Pi}\\
\node{HH_{k-1}(A_\text{c})}\arrow{e,t,T}{\pi _{k-1}}
\node{\Omega_\text{c}^{k-1}(M),}
\end{diagram}
\end{displaymath}
where $b_\Pi$ is Brylinski's Poisson differential,
$\epsilon_k $ is the antisymmetrization map defined as follows 
\begin{displaymath} 
  \epsilon_k ( f_0 df_1\wedge\cdots \wedge df_k ):=
  \sum_{\sigma \in S_k} \operatorname{sgn} (\sigma) \, f_0\otimes 
  f_{\sigma^{-1}(1)}\otimes 
  f_{\sigma ^{-1}(2)}\otimes \cdots \otimes f_{\sigma ^{-1}(k)} ,
\end{displaymath}
and $\pi_{k-1}$ is the projection defined by
\begin{displaymath}
  \pi_{k-1}(f_0\otimes f_1\otimes \cdots \otimes
  f_{k-1}):= f_0df_1\wedge\cdots \wedge df_{k-1}.
\end{displaymath}
Since $\epsilon_k$ resp.~$\pi_k$ gives rise to an isomorphism 
between $HH_k(A_\text{c}) $ and
$\Omega_\text{c}^k(M)$, these maps also induce an isomorphism between the 
Poisson homologies by the above diagram.
\end{proposition}
\begin{proof}
  See Theorem 3.1.1 in \cite{byi:poisson} and also \cite{brg}.
\end{proof}

Morita equivalence is an important notion in the
study of algebras by methods of noncommutative geometry. In the rest of 
this section, we will briefly look at  Morita invariance of 
Poisson homology. With respect to Poisson geometry, there exists
quite some work on the invariance of Poisson (co)homology under 
(weak) Morita equivalences between Poisson manifolds 
(cf.~\cite{xu}). 
In this paper we will now consider algebraic versions of Morita
invariance within noncommutative Poisson homology.

It is well-known that Morita equivalent algebras have isomorphic
Hochschild (co)homologies (see \cite{loday} 1.2.4, 1.2.7 and 1.5.6).

\begin{proposition}
\label{prop:morita-poisson} A Morita equivalence bimodule between
algebras with local units $A$ and $B$ defines an isomorphism between  
the sets of Poisson structures.
\end{proposition}

\begin{proof}
The Hochschild cohomology $H^k(A,A)$ is isomorphic to
$Ext_{A^e}^k (A, A)$. In \cite{keller:derived}, Keller shows that
for derived equivalent algebras (which is more general than Morita
equivalence), the canonical isomorphism defined by tensoring an
extension by the equivalent bimodule preserves the corresponding
$G$-brackets. This result implies that a Morita equivalence
bimodule between $A$ and $B$ defines an isomorphism between
$H^\bullet(A,A)$ and $H^\bullet(B,B)$ as Lie algebras, which induces an
isomorphism on the sets of corresponding Maurer-Cartan elements.
We know that in $H^\bullet(A,A)$ and $H^\bullet(B,B)$ Maurer-Cartan elements
are Poisson structures. Therefore, we have isomorphic sets of
Poisson structures.
\end{proof}

\begin{proposition}
\label{prop:morita} Under the isomorphism between the Hochschild
cohomologies of Prop.~\ref{prop:morita-poisson}, the
corresponding noncommutative Poisson structures have isomorphic
Poisson homologies.
\end{proposition}
\begin{proof}
Consider the following diagram
\begin{displaymath}
  \begin{diagram}\label{diag:morita}
  \node{HH_\bullet(A) \index{$HH_\bullet(A)$}}
  \arrow{e,t,T}{d^A_{\Pi}}\arrow{s,l,<>}{\sigma_\bullet}
  \node{HH_{\bullet -1}(A)}\arrow{s,l,<>}{\sigma_{\bullet-1}}\\
  \node{HH_\bullet(B)}\arrow{e,t,T}{d^B_{\Pi}}\node{HH_{\bullet-1}(B)}
\end{diagram}
\end{displaymath}
where $\sigma _\bullet$ is an isomorphism constructed by a Morita
equivalence. The claim of the proposition says that the above
diagram commutes. The proof of this goes along the same lines as
the proof of \cite{loday}(1.2.7). Instead of working out
the general case, we will only look at the special case where $B=M_{n}(A)$,
the $n\times n$ matrix algebra of $A$.

For $A$ and $M_n(A)$, following (1.2.4) of Loday \cite{loday}, we
define $\sigma_k: HH_k(A)\to HH_{k}(M_n(A))$
as follows (where $E_{11}^a$ denotes the matrix with $a$ at the $(1,1)$
position and $0$ elsewhere)
\begin{displaymath}
 \sigma_k(a_0\otimes a_1\otimes \cdots a_k):=
 E^{a_0}_{11}\otimes E^{a_1}_{11} \otimes \cdots \otimes E^{a_k}_{11},
\end{displaymath}
and a generalized trace map 
$\operatorname{tr}: HH_k(M_n(A))\to HH_k(A)$ by
\begin{displaymath}
 \operatorname{tr}(\fraka  _0\otimes \fraka _1 \cdots\otimes \fraka_k)
 :=(\fraka _0)_{i_0i_1}\otimes(\fraka_1)_{i_1i_2} \otimes 
 \cdots \otimes (\fraka_k)_{i_ki_0}.
\end{displaymath}

It is shown in (1.2.4) of \cite{loday} that both $\operatorname{tr}$ and
$\sigma_\bullet$ induce isomorphisms in Hochschild homology, and that
their dual versions give rise to isomorphisms in cohomology. Under the
corresponding isomorphisms in Hochschild cohomology, a Poisson
structure $\Pi$ on $A$ is transformed to $\tilde {\Pi}$, a Poisson
structure on $M_n(A)$, by the following formula
\begin{displaymath}
  \tilde {\Pi}\, (\fraka, \frakb )= 
  \Big(\sum_l \Pi \big( \fraka_{il}, \frakb_{lj} \big) \Big)_{ij}.
\end{displaymath}
One can easily show that $\tilde {\Pi}$ is a Poisson structure on $M_n(A)$
indeed.

Since $\sigma$ and $\operatorname{tr}$ are inverse to each other, the claim 
is proved, if one can show that 
$d_{\Pi} =\operatorname{tr}\circ d_{\tilde{\Pi}} \circ \sigma$. 
But this formula is obvious from the definition of $\tilde{\Pi}$.
\end{proof}

\begin{remark}
To prove Morita invariance for noncommutative Poisson structures
in general, one has to construct chain homotopies which entail
the above diagram to be commutative. The corresponding 
constructions are similar to those of (1.2.7) in \cite{loday}.
\end{remark}

\subsection{Poisson homology of the noncommutative Poisson algebra}
Assume to be given a proper \'etale Lie groupoid $G$ together with an invariant
Poisson bivector. For $\mathcal A$ the $G$-sheaf of smooth functions
we compute in this part the Poisson homology of the induced Poisson structure 
on the convolution algebra  $\mathcal A \rtimes G$ defined in the 
Preliminaries by Eq.~(\ref{def:pi}). 
We start by introducing the following notion.
\begin{definition}
\label{def:alm-bicom}
Let $(X_\bullet, b , d )$ be a triple consisting of 
a graded vector space $X_\bullet = \bigoplus_{k\in \N} X_k$ and two 
homogeneous maps 
$b : X_\bullet \rightarrow X_{\bullet \pm 1}$, 
$d: X_\bullet \rightarrow X_{\bullet \pm 1}$,
both either of degree $+1$ or $-1$. 
Then $(X_\bullet, b , d)$ is called an almost bicomplex, if 
the following relations hold true for some 
$h :X_\bullet\to X_{\bullet-1}$ (resp.~some
$h :X_\bullet\to X_{\bullet+1}$):
\begin{displaymath}
  b^2=0, \quad db+bd=0, \quad d^2=bh+hb.
\end{displaymath}
An almost bicomplex $(X_\bullet, b , d)$ gives rise to
two complexes $H_\bullet^b(X)$ and 
$H_\bullet^d(X)$, where $H_\bullet^b(X)$ is defined by $(X_\bullet, b)$,
while $H_\bullet^d(X)$ is defined by $(H_\bullet^b, d)$.
\end{definition}
\begin{example}
If $A$ is an algebra with a Poisson structure induced by a Hochschild cocycle
$\Pi$, and $d_{\Pi}$ is the differential of Poisson homology as defined above,
then the triple $(C_\bullet(A), b , d_{\Pi})$ is 
an almost bicomplex.
\end{example}
In the following computations, we will frequently use the next result.
\begin{lemma}
\label{lem:calculation} 
  Assume to be given two almost bicomplexes 
  \index{almost bicomplex} 
  $(X^i_\bullet, b^i , d^i )$, $i=1,2$ and a 
  quasiisomorphism $\Psi$ between $(X^1_\bullet, b^1 )$ and 
  $(X^2_\bullet,b^2)$.
  If $\Psi$ commutes with $d^i$ on the homologies 
  $H_\bullet(X^i_\bullet , b^i)$, i.e.~if 
  $\Psi d^2=d^1 \Psi$, then $ \Psi$ 
  induces an isomorphism between $(H^{b^1}_\bullet(X^1), d^1)$ and 
  $(H^{b^2}_\bullet(X^2), d^2)$.
\end{lemma}
\begin{proof}
  The proof of this lemma is obvious, since $\Psi$ is an isomorphism
  between the homologies and commutes with the differentials $d^i$.
\end{proof}
\begin{remark}
Of course we can allow for more general morphisms between bicomplexes to induce isomorphisms on Poisson homology. In particular, the quasi-isomorphism $\Psi$ in the lemma is allowed to ``commute up to homotopy'', i.e., $d^1\Psi=\Psi d^2+b^1H+Hb^2$, for some $H:X^1_\bullet\rightarrow X^2_\bullet$. The main difficulty of the computation below is to show that the ``reduction to loops'' morphism \eqref{reductiontoloops} which computes the Hochschild homology of the convolution algebra is a morphism of this kind: since the morphism does not preserve the Poisson differential, a homotopy as above is required. 
\end{remark}

With the above preparations, we are now ready to determine the Poisson
homology on the convolution algebra of the proper \'etale groupoid $G$. 
Our strategy is to track the change of the Poisson
differential in the various steps of the computation.

Recall that $HH_\bullet( \mathcal A \rtimes G)$ 
is calculated by the Bar complex 
$\big( (\mathcal A \rtimes G)^\natural_\bullet, b\big)$ the components of
which are isomorphic to the  vector spaces 
$\Gamma_c(G^{k+1}; s^*_{k+1}\mathcal A^{\boxtimes (k+1)})$. 
Under these isomorphisms, the Poisson differential on 
$(\mathcal A \rtimes G)^\natural_\bullet$ has the following form:
\begin{displaymath}
\begin{split}
  d_{\Pi} & \: := \sum_{i=0}^k(-1)^{i} {d^i_{\Pi}}, \quad \text{where} \\
  & d^i_{\Pi}(a_0, \cdots, a_k) := 
  \begin{cases}
     (a_0 \otimes \cdots \otimes \Pi(a_i, a_{i+1})\otimes \cdots 
     \otimes a_k),& 
     \text{if $0\leq i\leq k-1$},\\
     (\Pi(a_k, a_0)\otimes a_1 \otimes \cdots \otimes a_{k-1}), & 
     \text{if $i=k$}.
  \end{cases}
\end{split}
\end{displaymath}
We now proceed in three major steps.

{\bf Step I.} {\it Reduction to loops}.
Recall from the Preliminaries (Sec.~\ref{chog}, Step I.)
the method of reduction to loops for the computation of Hochschild and 
cyclic homology of an \'etale groupoid. This method shows that the 
cyclic vector space $(\mathcal A\rtimes G)^\natural$ 
is quasi-isomorphic to 
$\Gamma_\text{\rm c}\Lambda^\natural \mathcal A$ via the natural restriction
$p: (\mathcal A\rtimes G)^\natural_\bullet
 \to \Gamma_\text{\rm c}\Lambda^\natural_\bullet \mathcal A$. 
Now, one observes that $p$ naturally induces a Poisson differential 
$d_\Pi^\Lambda$ on 
$\Gamma_\text{\rm c}\Lambda^\natural_\bullet \mathcal A$ 
by putting
\begin{displaymath}
  d_\Pi^\Lambda(p(a)) := p (d_{\Pi}(a)) \quad
  \text{for all $a\in  (\mathcal A\rtimes G)^\natural_k$}.
\end{displaymath}
Note that $d_\Pi^\Lambda$ is well-defined indeed, since
$p(a)$ and $p(d_{\Pi}(a))$ are, respectively, the germs of $a$ and 
$d_{\Pi}(a)$ on $B^{(k)}$. 
Using Lemma \ref{lem:calculation}, we now conclude that the homology of 
$d_\Pi^\Lambda$ on 
$HH_\bullet(\Gamma_\text{\rm c}\Lambda^\natural_\bullet \mathcal A)$ 
is equal to the homology of $d_{\Pi}$ on $HH_\bullet(\mathcal A\rtimes G)$.

{\bf Step II.} {\it Homology of the cyclic groupoid}.
Recall from the proof of Prop.~\ref{propcyhomcycsheaf} 
in the Preliminaries that
the $\infty$-cyclic vector spaces 
$\Gamma_\text{\rm c}\Lambda^\natural_\bullet\mathcal A$ and
$\Gamma_\bullet (\Lambda G,\theta, \mathcal A^\natural_\text{tw})$ 
are isomorphic, with isomorphism over the stalk
at $(g_0,\cdots ,g_k) \in B^{(0)}$ given by Eq.~(\ref{isocyclicvs}).
By this isomorphism, the Poisson differential $d_\Pi^\Lambda$ gives rise
to a noncommutative Poisson differential $d_\Pi^\text{tw}$ on
$\big( \mathcal A^\natural_{\text{tw},\bullet}, b_\text{tw}\big)$, 
by which we can define a Poisson homology.  The explicit formulas are
(with $(g_0,g_1,\cdots,g_k) \in \Lambda G_k$)
\begin{displaymath}
\begin{split}
  (d_\Pi^\text{tw})_i &
  \big( [a_0\otimes \cdots \otimes a_k]_{(g_0,g_1, \cdots,g_k)}\big)
  =\\ & = 
  \begin{cases}
    [\{a_0 , a_1\}g_1 \otimes a_2g_1 \otimes \cdots \otimes
    a_k g_1]_{(g_1^{-1} g_0 g_1, g_2, \cdots, g_k)},&
    \text{ for $i=0$,}\\
    [a_0 \otimes \cdots \otimes \{a_i, a_{i+1}\} \otimes \cdots \otimes 
    a_k]_{(g_0, g_1,\cdots, g_ig_{i+1},
    \cdots, g_k )},& \text{ for $1\leq i\leq k-1$,}\\ 
    [\{ a_k , a_0\} \otimes a_1 \otimes \cdots \otimes 
    a_{k-1}]_{(g_0, g_1, \cdots, g_{k-1})},&
    \text{ for $i=k$,}
  \end{cases}
\end{split}
\end{displaymath}
where $\{\ ,\ \}$ is the Poisson bracket on $G_0$.

Since the above isomorphism is a local diffeomorphism which maps
Poisson structures naturally, we conclude by Lemma \ref{lem:calculation}
that 
$\big( \Gamma_\bullet(\Lambda G, \theta, \mathcal A^\natural_\text{tw}), 
  b_\text{tw}, d_\Pi^\text{tw} \big)$
calculates the Poisson homology of $\mathcal A \rtimes G$. 
We define the Poisson homology 
$H^{\Pi}_\bullet(\Lambda G)_\text{tw}$ as the homology of $d_\Pi^\text{tw}$
on $HH_\bullet(\Lambda G, \theta, \mathcal A^\natural_\text{tw} )= 
H_\bullet \big(\Gamma_\bullet (\Lambda G, \theta,\mathcal A^\natural_\text{tw}),
b_\text{tw}\big)$.

From the above considerations one can now immediately derive the 
following localization property similarly to the corresponding one 
for Hochschild homology (cf.~\cite{bn,crainic}):
\begin{theorem}
\label{thm:locality} 
  Let $B^{(0)}=\bigcup_{\mathcal O \in \operatorname{Sec} (G)} \mathcal O$ 
  be the decomposition of $B^{(0)}$ into sectors. Then
  \begin{displaymath}
     H^{\Pi}_\bullet(\mathcal A \rtimes G)
     =\bigoplus_{\mathcal O \in \operatorname{Sec} (G)}
     H^{\Pi}_\bullet(\mathcal A \rtimes G)_{\mathcal O}.
  \end{displaymath}
\end{theorem}

{\bf Step III.} {\it Inertia groupoid}.
Since the groupoid is proper, one knows that the 
Poisson structure on $G_0$ defines a natural invariant
Poisson structure on $B^{(0)}$, which gives rise 
to an invariant Poisson bivector $\Pi_N$ on $NG_0$
(see \cite{ta:thesis} for a detailed proof in the general case and 
 Lemma \ref{fp} in the following section for groupoids with a 
 symplectic structure).
Now, over an invariant open-closed subset $\mathcal O$ of $NG_0$, we consider 
$\calc^{\infty}_{\mathcal O}$, 
the sheaf of smooth functions, together with the Hochschild 
differential $b_{\mathcal O}$ and the Poisson differential 
$d_{\Pi_N}$ defined by the Poisson structure 
$\Pi_N$. It is straightforward to check that 
$\big( (\mathcal C^{\infty}_{\mathcal O})^{\natural}_\bullet, b_{\mathcal O}, 
d_{\Pi_N}\big)$ forms an almost bicomplex. We define the Poisson 
homology $H^{\Pi}_\bullet(NG_{|\mathcal O})$ to be the homology of 
$d_{\Pi_N}$ on 
$HH_\bullet( NG,(\mathcal C^{\infty}_{\mathcal O})^\natural)$. 
By Proposition \ref{commutative} and Lemma \ref{lem:calculation},  we 
conclude that 
$\big( (\mathcal C^\infty_{\mathcal O})^{\natural}_\bullet,
 b_{\mathcal O}, d_{\Pi_N} \big)$ 
is quasi-isomorphic (as an almost bicomplex) to the bicomplex 
$(\Omega^\bullet_{\mathcal O}, 0, b_{\Pi_{N}})$, hence 
$H^{\Pi}_\bullet(NG_{|\mathcal O})$ is given by the homology of 
$b_{\Pi_N}$ on 
$H_\bullet( NG,\Omega^\bullet _{\mathcal O})$, 
which we will denote by $H^{\Pi}_\bullet (NG)_{\mathcal O}$.

\begin{theorem}
\label{thm:poi-hom}
For an invariant open-closed subset $\mathcal O \subset B^{(0)}$, one has
\begin{displaymath}
  H^{\Pi}_\bullet(\mathcal A \rtimes G )_{\mathcal O} =
  H^{\Pi}_\bullet( NG_{|\mathcal O})= H^{\Pi}_\bullet(NG)_{\mathcal O}.
\end{displaymath}
\end{theorem}

\begin{proof}
The second equality in the claim has been shown above, so it remains
to prove the first one. 
To this end recall first the twisted Hochschild--Kostant--Rosenberg Theorem
\cite[Lem.~3.1.5]{crainic} which entails that the natural restriction of the 
germ of a smooth function to $\mathcal O$ induces 
a quasi-isomorphism $\rho : (\mathcal A_{\mathcal O}^{\natural},
b_\text{tw})\rightarrow ((\mathcal C^{\infty}_{\mathcal O})^{\natural},b)$.
Below, we will show that via $\rho$ one can pushforward $d_\Pi^\text{tw}$ 
to a Poisson differential $d'_{\Pi_N}$ on 
$HH_\bullet(NG,(\mathcal C^{\infty}_{\mathcal O})^{\natural})$, 
which then calculates the Poisson homology of the convolution algebra 
$\mathcal A \rtimes G$. Moreover, we will show that
$d'_{\Pi_N}$ is equal to $d_{\Pi_N}$ on 
$HH_\bullet(NG, (\mathcal C^{\infty}_{\mathcal O})^{\natural}))$. 
This will prove the claim.
Note that the problem one has to cover here is the fact that due to the existence 
of normal directions, $\rho$ does not induce a Poisson map
from $(\mathcal A (  \mathcal O),\Pi)$ to
$(\mathcal C^{\infty}_{\mathcal O} (\mathcal O),\Pi_N)$.

Let us now construct $d'_{\Pi_N}$ in detail.
For all $a\in \Gamma_\bullet (NG, 
(\mathcal C^\infty_{\mathcal O})^{\natural})$ with
$b(a)=0$ one can find an 
$x\in \Gamma_\bullet (\Lambda G, \mathcal A_{\mathcal O}^{\natural})$
with $\rho(x)=a$ and $b_\text{tw}(x)_{|\mathcal O}=0$. 
We define $d'_{\Pi_N}$ on $a$ by 
$$
  d'_{\Pi_N}(a):=\rho \big(d_\Pi^\text{tw} (x)\big).
$$
The following is a list of properties of $d'_{\Pi_N}$.

(i) 
Since $b_\text{tw}(x)_{|\mathcal O}=0$, we have that
$d_\Pi^\text{tw}(b_\text{tw}(x))$ vanishes on $\mathcal O$ as well, 
hence one has
$\rho (d_\Pi^\text{tw}(b_\text{tw}(x)))=0$. 
This implies the following equality:
\begin{displaymath}
  b( d'_{\Pi_N} (a))=b(\rho (d_\Pi^\text{tw} (x)))=\rho(b_\text{tw}(d_\Pi^\text{tw}(x)))
  =-\rho( d_\Pi^\text{tw} (b_\text{tw}(x)))=0,
\end{displaymath}
where in the second equality we have used that $\rho$ commutes with $b$,
and in the third one we have used that $d_\Pi^\text{tw}$ anti-commutes 
with $b_\text{\rm tw}$.

(ii) Since $\rho$ is an isomorphism on the Hochschild homology, we can choose
for any $[a]\in HH_\bullet(NG, (\calc^{\infty}_{\mathcal O})^{\natural})$
a representative 
$a\in \Gamma_\bullet ( NG, (\mathcal C^{\infty}_{\mathcal O})^{\natural})$
such that there is
$x\in \Gamma_\bullet (\Lambda G, \mathcal A_{\mathcal O}^{\natural})$ 
satisfying
$\rho(x)=a$ and $b_\text{\rm tw}(x)=0$. Thus, we obtain
$$
  b_\text{\rm tw} (d_\Pi^\text{tw}(x))=-d_\Pi^\text{tw}(b_\text{tw}(x))= 0.
$$ 
Hence, by construction,
$$
  d'_{\Pi_N} \circ d'_{\Pi_N} (a)= \rho (d_\Pi^\text{tw} (d_\Pi^\text{tw}(x)))=
  \rho \big( b_\text{tw} (h(x))+ h(b_\text{tw}(x)) \big)=b(\rho(h(x))),
$$
where $h$ is the homotopy associated to $(d^\text{tw}_{\Pi})^2$. This shows
that $d'_{\Pi_N}\circ d'_{\Pi_N}=0$ in 
$HH_\bullet(NG,(\mathcal C^{\infty}_{\mathcal O})^{\natural})$.

(iii) To prove that $d'_{\Pi_N}$ is well defined, we need to show that
our definition is independent of the choices of $a$ and $x$. We will
use the following lemmas:
\begin{lemma}
\label{lem:definition} 
  Let $y\in \Gamma_\bullet (\Lambda G ,\mathcal A_{\mathcal O}^{\natural}) $
  with $\rho(y)=0$ and $b_\text{\rm tw} (y)=0$. Then there exists 
  $z \in \Gamma_\bullet (\Lambda G ,\mathcal A_{\mathcal O}^{\natural}) $
  such that $y=b_\text{\rm tw}(z)$.
\end{lemma}
\begin{lemproof} 
  Since $\rho$ is a quasi-isomorphism and $\rho (y)=0$, $y$
  has to be a boundary in 
  $\Gamma_\bullet (\Lambda G ,\mathcal A_{\mathcal O}^{\natural})$. 
  Therefore, there is a 
  $z\in \Gamma_\bullet (\Lambda G ,\mathcal A_{\mathcal O}^{\natural})$
  such that $y=b_\text{\rm tw}(z)$.
\end{lemproof}

By the lemma one now concludes that for any 
$x, y\in \Gamma_\bullet (\Lambda G, \mathcal A_{\mathcal O}^{\natural})$ 
with $\rho(x)=\rho(y)=a$ and $b_\text{tw}(x)=b_\text{tw}(y)=0$
there exists $z$ in 
$\Gamma_\bullet (\Lambda G, \mathcal A_{\mathcal O}^{\natural})$ , 
such that $x-y= b_\text{\rm tw}(z) $. Therefore,
\begin{displaymath}
\begin{split}
\rho(d_\Pi^\text{tw} (x))&=\rho(d_\Pi^\text{tw} (y))+
\rho(d_\Pi^\text{tw} (x-y))\\
&=\rho(d_\Pi^\text{tw}(y))+\rho(d_\Pi^\text{tw} (b_\text{tw}(z)))\\
&=\rho(d_\Pi^\text{tw}(y))+\rho(d_\Pi^\text{tw} (b_\text{tw}(z)))\\
&=\rho(d_\Pi^\text{tw}(y))-\rho(b_\text{tw}(d_\Pi^\text{tw} (z)))\\
&=\rho(d_\Pi^\text{tw}(y))-b(\rho(d_{\Pi} (z))).
\end{split}
\end{displaymath}

Hence, the homology class of $d'_{\Pi_N} (a)$ is independent of
the lift of $a$. To show that it is also 
independent of the representative $a$ in the homology class $[a]$, 
we prove the following proposition.
\begin{lemma}
\label{lem:transition} 
 $d'_{\Pi_N}$ is equal to the Poisson differential 
 $d_{\Pi_N}$ on $HH_\bullet(NG, (\calc^{\infty}_{\mathcal O})^{\natural})$.
\end{lemma}

\begin{lemproof}
For any 
$a\in \Gamma_\bullet (NG, (\calc^{\infty}_{\mathcal O})^{\natural})$ 
with $b(a)=0$ we construct a particular lift $x$ in 
$\Gamma_\bullet (\Lambda G, (\mathcal A_{\mathcal O}^\natural))$, such
that $b_\text{tw}(x)_{|\mathcal O}=0$. To achieve this,
recall that $\mathcal O$ is embedded in $G ^{(k+1)}$ with normal
bundle being the set of all nontrivial representations of $G$
on $(TG^{(k+1)})_{| \mathcal O}$. By the tubular neighborhood theorem,
one can find a function $x \in \mathcal C_\text{c}^\infty(G^{(k+1)})$
which is equal to the pull back of $a$ in a tubular neighborhood of
$\mathcal O$. Hence, $b(x)=0$ in the tubular neighborhood, and 
$d_\Pi^\text{tw} (x)_{|\mathcal O}$ is equal to
$d_{\Pi_N}(a)$. Therefore, we have
$d'_{\Pi_N}(a)=d_\Pi^\text{tw}(x)_{| \mathcal O}=d_{\Pi_N}(a)$. 
\end{lemproof}

By Lemma \ref{lem:transition} and the fact that $d_{\Pi_N}$ is
well-defined on 
$HH_\bullet(NG, (\calc^{\infty}_{\mathcal O})^{\natural}) $, we 
obtain that $d'_{\Pi_N} $ acts on
$HH_\bullet(NG, (\calc^{\infty}_{\mathcal O})^{\natural})$ 
independent of the choice of representatives $a$.
In other words, $d'_{\Pi_N}$ is well-defined on 
$HH_\bullet(NG, (\calc^{\infty}_{\mathcal O})^{\natural})$ and equal to 
$d_{\Pi_N}$. Altogether, this finishes the proof.
\end{proof}

\begin{remark}
 Lemma \ref{lem:calculation} now entails that we can use 
 $\big(\Gamma_\bullet (NG , 
 (\calc^{\infty}_{\mathcal O})^{\natural}), b, d_{\Pi_N}\big)$ 
 to calculate the Poisson homology of the localized convolution algebra
 $\mathcal A \rtimes G_{\mathcal O}$.
\end{remark}
\section{Hochschild and cyclic homology of the quantized algebra}
\label{cyhom}
In this section we present the computation of Hochschild and cyclic homology 
of a formal deformation quantization of the convolution algebra 
on a proper \'etale Lie groupoid $G$ representing a symplectic orbifold $X$.
By $\omega$ we denote the symplectic form on $G_0$, and by $\mathcal A$,
as before, the $G$-sheaf of smooth functions on $G_0$. The deformation 
quantization is constructed as explained in Section \ref{sec:intr-ncq}.
We choose a $G$-invariant star product $\star$  on  the sheaf
$\mathcal A[[\hbar]]$ of formal power series. Notice that we can assume 
without loss of generality that this deformation is a Fedosov star product 
associated to an invariant symplectic connection. This gives rise to the 
$G$-sheaf $\mathcal A^\hbar=(\mathcal A[[\hbar]] ,\star)$ and to the 
deformed global crossed product algebra $\A^\hbar\rtimes G$ as in 
Eq.~\eqref{crossedprd}.
\subsection{Periodic cyclic homology}
The computation of the periodic cyclic homology groups of $\A^\hbar\rtimes G$ 
follows at once from the ``classical'' computations of the periodic cyclic 
homology of \'etale groupoids in \cite{bn,crainic} by the following 
rigidity property \cite{getzler}, \cite[Thm.~A2.2]{nt}: 
For any formal deformation quantization
$A^\hbar:=(A[[\hbar]],\star)$ of an algebra $A$, one has an isomorphism 
$$
  HP_\bullet(A^\hbar)\cong HP_\bullet(A)\otimes\C[[\hbar]].
$$ 
Therefore, one easily finds
\begin{proposition}
\label{pch}
  The periodic cyclic homology groups of 
  $\A^\hbar\rtimes G$ are given by:
  $$
    HP_\bullet(\A^\hbar\rtimes G)=
    \prod_k H^{2k+\bullet}_{\text{\tiny\rm orb,c}}(X,\C[[\hbar]]).
  $$
\end{proposition}
A similar rigidity property of (algebraic) $K$-theory was proved in
\cite{rosenberg}. Recall that the Chern--Connes character maps 
$K^{\mbox{\tiny alg}}_\bullet(\A\rtimes G)$ to $HP_\bullet(\A\rtimes G)$. 
The two rigidity isomorphisms are compatible with this character map.
\subsection{Computation of Hochschild homology}
The computation of Hochschild and cyclic homology is more involved. The 
main tools in the computation are: \\[1mm]
1) The ``quantum to classical'' spectral sequence induced by the $\hbar$-adic 
filtration introduced in \cite{brg}. \\
2) The ``classical'' computation of cyclic homology of \'etale groupoids of 
\cite{bn,crainic} using the language of sheaves. The computation exactly 
follows the steps of these computations: we first ``localize'' to a sheaf 
cohomology computation on the inertia groupoid, and then use the 
$\hbar$-filtration to reduce the outcome to orbifold cohomology.

In the following, we will occasionally work over the field 
$\C((\hbar))$, so let us put
$\A^{((\hbar))}:=\A^\hbar\otimes_\C\C[\hbar^{-1}]$. 
Then one has $\A^{((\hbar))} \in \mathsf{Sh} (G)$ as well. For the following,
notice that when $G$ has an invariant symplectic form $\omega$, the 
pull-back $\beta^*\omega$ defines an invariant symplectic structure on 
$\Lambda G$, which descends to $NG$, cf.~also Prop.~\ref{fp} below.
\begin{theorem}
\label{qhochschild}
Let $G$ be a proper \'etale Lie groupoid representing a symplectic
orbifold $X$ of dimension $2n$. Then the Hochschild homology of the 
deformed convolution algebra $\A^{((\hbar))}\rtimes G$ is given by
\begin{displaymath}
  HH_\bullet(\A^{((\hbar))}\rtimes G) \cong 
  H^{2n-\bullet}_{\text{\tiny\rm orb,c}}\left(X,\C((\hbar))\right)
\end{displaymath}
\end{theorem}
\begin{proof}
Consider the spectral sequence induced by the $\hbar$-adic filtration on 
the Hoch\-schild complex of a formal deformation $A^\hbar=(A[[\hbar]],\star)$
 of a noncommutative Poisson algebra $(A,[\Pi])$ in the sense of 
Definition \ref{dfn:quant-noncommutative-poisson}. Clearly, in degree zero one 
finds the classical Hochschild complex of $A$ and a straightforward 
computation shows that the differential 
$d^1:E^1_{p,q}\rightarrow E^1_{p-1,q}$ is given by 
the noncommutative Poisson differential 
$d_\Pi: HH_{p+q}(A)\rightarrow HH_{p+q-1}(A)$ of 
Definition \ref{def:poisson}. In our case, with 
$A^\hbar=\A^\hbar \rtimes G$, $A=\mathcal C^\infty_c(G)$, one finds
\begin{displaymath}
  E^1_{p,q} = HH_{p+q}(A) \cong \Omega^{p+q}_{\text{\tiny inv, c}}(NG). 
\end{displaymath}
The last isomorphism follows from the computations of \cite{crainic}. 
For general \'etale groupoids this also includes higher cohomology groups 
$H_\text{c}^k(NG,\Omega^l)$, but these all vanish here by the following 
argument.
Using the projection $\pi:G\rightarrow X$ onto the orbifold, one 
identifies $\pi_!\Omega^\bullet$ as the sheaf on $X$ of invariant forms, 
which is fine, see \cite{pflaum}, and vanishing of cohomology follows 
from \eqref{d}. 

As explained in Sec.~\ref{ncph}, the differential $d^1$ is nothing but 
Brylinski's Poisson differential \eqref{poissondiff} on the invariant 
differential forms on the groupoid $NG$, which is well-defined because the 
induced symplectic form is invariant.
Of course, this is the image of the sheaf version of Brylinski's complex 
$(\Omega^\bullet_{NG},d_\Pi)$ in $\mathsf{Sh}(NG)$, under the functor 
$\Gamma_{\text{\tiny inv,c}}$. Using the fact that $NG$ is symplectic, 
there is a quasi-isomorphism between this sheaf complex and the 
de Rham complex, as in \cite{byi:poisson}:
$$
  (\Omega^\bullet_{\orbit},d_\Pi)\cong
  (\Omega_\orbit^{\dim\orbit-\bullet},d_\text{\tiny dRh}).
$$ 
Here, we restrict to a connected component of the inertia groupoid because 
the components may have different dimensions, affecting the degree shift in 
the isomorphism. By the fact that $\Omega^\bullet$ is 
$\Gamma_{\text{\tiny inv,c}}$-acyclic, on therefore finds 
$$
   E^2_{p,q}=H^{2n-p-q}_{\text{\tiny orb,c}}(X,\C((\hbar))).
$$
As in the case of smooth manifolds, we claim that the spectral sequence 
degenerates at this stage. When $NG$ has a finite number of components, 
the argument is the same as that of \cite{nt}: If the spectral sequence 
does not degenerate at this stage, one has 
$$
   \dim_{\C((\hbar))}HH_i(\A^\hbar\rtimes G)<b^{2n-i}_{\text{\tiny orb,c}}(X),
$$ 
where $b^i_{\text{\tiny orb, c}}(X)=\dim H^i_{\text{\tiny orb,c}}(X)$. But 
this would imply that  
$$
  \dim_{\C((\hbar))}HP_i(\A^\hbar\rtimes G)<
  \sum_kb^{2k+i}_{\text{\tiny orb,c}}(X),
$$ 
contradicting Prop.~\ref{pch}. The last inequality hereby follows
from  the spectral 
sequence from Hochschild homology to cyclic homology obtained by 
filtering the $(b,B)$-complex by columns. When $\tilde{X}$ however does have 
infinite many connected components, one uses the decomposition of Hochschild
homology induced by the $G$-invariant decomposition \eqref{orbits} of $NG$: 
$$
  HH_\bullet(\A^\hbar\rtimes G)=
  \bigoplus_\orbit HH_\bullet(\A^\hbar\rtimes G)_\orbit.
$$ 
This decomposition can be deduced from the similar decomposition in 
Poisson homology, cf. Theorem \ref{thm:locality}, by the spectral sequence 
above. The star product is given in terms of local multidifferential 
operators on $G_0$, and therefore the decomposition is preserved by the 
higher differentials. By Proposition \ref{pch}, the periodic cyclic homology 
has a similar decomposition, and it is not too difficult to see that both are
compatible. Since $G$ is proper, each $\orbit$ has a finite 
number of connected components, and the above  dimension argument, 
applied to each component separately, proves that the 
spectral sequence degenerates at the second stage. The theorem follows.
\end{proof}
\begin{remark}
  We expect a similar result to hold for Poisson structures, provided one 
  can prove the corresponding analogue of Prop.~\ref{fp} below. In that case, 
  it follows by formality \cite{dg:formality-forms,sh:formality-forms} that 
  the Hochschild homology is given by the Poisson homology of the 
  groupoid $NG$.
\end{remark}
\subsection{Cyclic homology}
Next, we proceed to compute the cyclic homology, analogously to the 
computations of \cite{bn} and \cite{crainic}. We begin with a detailed 
analysis of some natural sheaves on the space of loops $B^{(0)}$. Consider 
the morphism $\beta^{-1}:\mathsf{Sh}(G)\rightarrow\mathsf{Sh}(\Lambda G)$. 
Applied to the quantum sheaf of algebras $\A^\hbar$ on $G$, we get a 
$\theta$-cyclic sheaf $\beta^{-1}\ \big( (\A^\hbar)^\natural \big)$ 
on $\Lambda G$.
Again, its stalk at $g\in B^{(0)}$ is given by 
$\big( \A^\hbar_g\big)^\natural_{\theta_g}$, cf.~\cite[3.4.1.]{crainic} and 
Sec.~\ref{chog}. 
\begin{proposition}
\label{fp}
  The inertia groupoid $\Lambda G$ carries a natural symplectic structure 
  $\omega_{(0)} := \beta^*\omega$, which descends to $NG$. 
  If $\mathcal A_{(0)}$ denotes 
  the sheaf of smooth functions on $B^{(0)}$, then there is a Poisson
  morphism of sheaves on $\Lambda G$
  \begin{equation}
   \phi_0 : \mathcal A_{(0)} \rightarrow \beta^{-1} \mathcal A ,
  \end{equation} 
  where $\mathcal A_{(0)}$ carries the Poisson structure induced by 
  $\omega_{(0)}$, and $\beta^{-1} \mathcal A $ inherits the Poisson 
  structure from  $\omega$.  Moreover, the sheaf 
  $\beta^{-1}(\A^\hbar)\in \mathsf{Sh}(\Lambda G)$ is isomorphic to a
  formal deformation quantization with coefficients of the symplectic 
  structure $\beta^*\omega$.
\end{proposition}
\begin{proof}
The first and last claim are essentially local statements and therefore 
we can restrict to the case of a 
translation groupoid $G=\Gamma\ltimes M$ by the action of a finite group. 
In this case one has a decomposition 
$$ 
  B^{(0)}=\coprod_{\gamma\in\Gamma}M^\gamma.
$$ 
As shown in \cite[Sec.~5]{fe:g-index}, the pull back of the symplectic form 
along the embedding $M^\gamma\hookrightarrow M$ for all $\gamma\in\Gamma$ 
gives $B^{(0)}$ a symplectic structure and $\beta^{-1}\A^\hbar$ is a formal 
deformation quantization of $(\Lambda G,\omega_{(0)} )$ with coefficients 
in the normal bundle with respect to this embedding. 

For the construction of $\phi_0$ choose a $G$-invariant Riemannian metric on 
$G_0$. Let $g\in B^{(0)}$ and consider a germ $[f]_g \in (\mathcal A_{(0)})_g$.
Let $U_g \subset B^{(0)}$ be an open neighborhood of $g$ on which $f$ is 
defined and put $M_g := s(U_g)$. Then $M_g$ is a submanifold of $G_0$ 
and there is a projection $\pi_g: T_g \rightarrow M_g$ from
a tubular neighborhood $T_g$ to $M_g$ along geodesics of the chosen
$G$-invariant  Riemannian metric. We now define
$\phi_0 ([f]_g)$ as the germ $[f \circ s^{-1} \circ \pi_g]_g$.
Since the normal bundle to $M_g$ (with respect to the above Riemannian metric)
is a symplectic bundle by the symplectic slice theorem, it is clear, 
that the thus defined sheaf morphism is Poisson.
\end{proof}
\begin{remark}
\label{right-inverse}
The morphism $\phi_0$ is a right inverse to the restriction map 
$f\mapsto f|_{B^{(0)}}$ used by Brylinski--Nistor, which is a 
quasi-isomorphism on Hochschild homology, cf.~\cite[Lemma 5.2.]{bn}, but 
is \textit{not} Poisson. As a right inverse, 
$\phi_0:(\A^\natural_k,b)\rightarrow(\beta^{-1}\A^\natural_k,
b_{\rm{\tiny tw}})$ 
is a quasi-isomorphism on the sheafified twisted Hochschild complexes 
as well. 
\end{remark}
Next, we want to lift this Poisson morphism to a ``quantum morphism'' 
$\phi:\A^\hbar_{(0)}\rightarrow\beta^{-1}\A^\hbar$, where $\A^\hbar_{(0)}$ 
is some suitable deformation quantization of $(B^{(0)},\omega_{(0)})$. 
Consider now the characteristic class 
  $[\star] \in \frac{[\omega]}{\hbar} + H^2(G_0,\mathbb{C} [[\hbar]])$
  of the star product on $\mathcal A^\hbar$, and let 
  $c \in \frac{[\omega_{(0)}]}{\hbar} + H^2(B^{(0)},\mathbb{C} [[\hbar]])$
  be the pullback class $\beta^* [\star]$. Since $\star$ is $G$-invariant,
  there exists a $G$-invariant Fedosov star product $\star_{(0)}$ on 
  $\mathcal A_{(0)}$ with characteristic class given by $[\star_{(0)}]= c$.
  Denote the
  resulting deformed sheaf of algebras by $\mathcal A_{(0)}^\hbar$.
  By construction of $\star_{(0)}$ and Thm.~5.4 
  in \cite{bordemann}, one now concludes that the Poisson morphism 
  $\phi_0:\mathcal A_{(0)} \rightarrow  \beta^{-1} \mathcal A$
  from Prop.~\ref{fp} indeed can be extended to a morphism of sheaves of 
  algebras 
  $$
   \phi = \sum_{k=0}^\infty \phi_k \hbar^k : 
   \mathcal A_{(0)}^\hbar \rightarrow  \beta^{-1} \mathcal A^\hbar,
  $$
  where each $\phi_k$ is a linear morphism of sheaves from 
  $\mathcal A_{(0)}$ to $\beta^{-1} \mathcal A$. We now have:
\begin{proposition}
\label{quantcycquasiiso}
The morphism $\phi:\A^\hbar_{(0)}\rightarrow\beta^{-1}\A^\hbar$ induces a 
quasi-isomorphism 
$$
  \phi_*:(\mathcal A_{(0)}^\hbar)^\natural \rightarrow 
  (\mathcal A^\hbar)^\natural_\text{\rm tw}
$$ of 
$\theta$-cyclic sheaves on $\Lambda G$, where the left hand side carries the
trivial cyclic structure.
\end{proposition}
\begin{proof}
Let us start with the Hochschild complex. Clearly, the map $\phi$ induces 
a map on the (twisted) Hochschild homology. Consider now the spectral
sequence induced by the $\hbar$-adic filtration. Since the morphism 
$\phi$ induces an isomorphism at level $E^1$, see Remark \ref{right-inverse},
it must induce an isomorphism on the level of Hochschild homology of the 
quantum sheaves by the spectral sequence comparison theorem, 
cf.~\cite[Thm. 5.2.12]{weibel}. As for cyclic homology, since the 
morphism of mixed complexes induced by $\phi$ is a quasi-isomorphism on 
Hochschild homology, it must induce an isomorphism on cyclic homology as 
well (cf.~Prop.~2.5.15 in \cite{loday}).
\end{proof}
Next, we consider the ``reduction to loops'' map 
\begin{equation}
\label{redloops}
 p: (\A^\hbar\rtimes G)^\natural_k\rightarrow 
 \Gamma_c(B^{(k)},\sigma_{k+1}^{-1}(\A^\hbar)^{\boxtimes (k+1)}),
\end{equation}
restricting sections over $G^{k+1}$ to $B^{(k)}$. As explained in 
Section \ref{chog}, the right hand side equals the Bar complex computing 
the $\theta$-twisted homology of the twisted cyclic sheaf 
$\beta^{-1}\A^\hbar$ on $\Lambda G$. Now we are in a position to prove:
\begin{proposition}
\label{qiso}
Reduction to loops induces a quasi-isomorphism on Hochschild and cyclic 
homology: 
\begin{eqnarray*}
 HH_\bullet\big(\A^{((\hbar))}\rtimes G\big)&\cong&
 HH_\bullet\big(\Lambda G,\theta ,
   \big(\A^{((\hbar))}\big)^\natural_\text{\rm tw}\big),\\
 HC_\bullet\big(\A^{((\hbar))}\rtimes G\big)&\cong&
 HC_\bullet\big(\Lambda G,\theta, 
  \big(\A^{((\hbar))}\big)^\natural_\text{\rm tw}\big).
\end{eqnarray*}
\end{proposition}
\begin{proof}
Again, we start with Hochschild homology. Using the quasi-isomorphism of 
Prop.~\ref{quantcycquasiiso}, the right hand side of the first equation 
above is isomorphic to the homology of sheafified Hochschild complex 
associated to a formal deformation quantization of the inertia groupoid
$(\Lambda G,\omega_{(0)})$. But this sheaf homology is readily computed, 
cf.~Prop.~\ref{quasi-iso} below, to give the orbifold cohomology as in 
Theorem \ref{qhochschild}. It is not difficult to show that in fact the 
isomorphism in Theorem \ref{qiso} is equal to the map induced by $p$ in 
\eqref{redloops}. Therefore we conclude that reduction to loops induces 
an isomorphism on Hochschild homology. However, given that it is a 
quasi-isomorphism on Hochschild homology, it induces an isomorphism on 
cyclic homology as well.
\end{proof}
Having reduced the computation cyclic homology to sheaf cohomology, we can 
do the computations locally for $\mathcal A_{(0)}$ and then take cohomology. 
The local results are given in the following proposition. For 
completeness, we also state the analogous results for Hochschild homology.
\begin{proposition}
\label{quasi-iso}
 On a proper \`etale Lie groupoid with symplectic structure of dimension $2n$
 there are quasi-isomorphisms of complexes of sheaves:
 \begin{displaymath}
 \begin{split}
   (\mathcal{C}_\bullet(\A^{((\hbar))}),b)&\cong \underline{\C}((\hbar))[2n],\\
   (\operatorname{Tot}_\bullet (\mathcal B_{\bullet,\bullet}
   (\A^{((\hbar))})),b+B)&\cong 
   \bigoplus_{k \in \mathbb{N}}\underline{\C}((\hbar))[2n+2k].
 \end{split}
 \end{displaymath}
\end{proposition}
\begin{proof}
  There are obvious inclusions of the right hand side into the sheaf 
  complexes on the left hand side. To prove that these are 
  quasi-isomorphisms, one needs to prove that they induce isomorphisms of 
  the stalk-wise cohomology, or, equivalently, of the homology sheaves. 
  However, locally the sheaf $\A^{((\hbar))}$ is isomorphic to the 
  Weyl algebra $\mathbb{W}^{((\hbar))}_n$ of formal Laurent series on 
  $\C^n$. Its Hochschild and cyclic homology are given as follows 
  (see \cite[Sec.~3.2]{nt}), where the second identity follows from 
  the first by the spectral sequence induced by filtering the 
  $(b,B)$-complex by columns.
\begin{displaymath}
\begin{split}
  HH_k\big( \mathbb{W}^{((\hbar))}_n \big) & =
  \begin{cases}
    \C((\hbar)), & \text{for $k=2n$}, \\
    0 , & \text{else},
  \end{cases} \\
  HC_k \big( \mathbb{W}^{((\hbar))}_n \big) & =
  \begin{cases}
    \C((\hbar)), & \text{for $k=2n+2l$},\\
    0, & \text{else}.
  \end{cases}
\end{split}
\end{displaymath}
This proves that the canonical inclusion is a quasi-isomorphism.
\end{proof}
\begin{theorem}
  Let $G$ be a proper \`etale Lie groupoid and $\mathcal A^\hbar$ like above.
  Then the cyclic homology of $\A^{((\hbar))}\rtimes G$ is given by 
  $$
    HC_\bullet(\A^{((\hbar))}\rtimes G)=
    \bigoplus_{k\geq 0}H_{\text{\rm \tiny orb,c}}^{2n+2k-\bullet}(X,\C((\hbar))).
  $$
\end{theorem}
\begin{proof}
 By Propositions \ref{fp}, \ref{qiso} and \ref{quantcycquasiiso},
 the cyclic homology of $\A^{((\hbar))}\rtimes G$ equals the hyperhomology 
 of the total complex of the sheafified Connes' $(b,B)$-complex of a 
 formal deformation of the sheaf of smooth functions on $\Lambda G_0$. 
 The preceding Proposition then entails the claim.
\end{proof}
\subsection{Cohomology}
Having computed the Hochschild and cyclic homology, the dual results may be 
computed analogously, as in \cite{crainic}. Recall that the Hochschild 
cohomology of an algebra is computed from the complex 
$A^k_{\natural}:=\Hom(A^{\otimes k+1},\C)$ dual to the Hochschild complex, 
with the corresponding differential. Of course, in our case 
we use topological tensor products and consider only continuous functionals. 
Again using the $\hbar$-adic filtration, one observes that the $E^1$-term 
in the spectral sequence is given by the complex of de Rham currents 
$(\Omega_\bullet, b'_\Pi)$ on $NG$ with the dual Poisson differential. 
The isomorphism of \cite{byi:poisson} dualizes on each component $\orbit$ 
to give 
$(\Omega_\bullet,b'_\Pi)\cong (\Omega_{\dim\orbit-\bullet},
d'_\text{\tiny dRh})\cong (\Omega^\bullet,d_\text{\tiny dRh})$. 
As for homology, the spectral sequence collapses at the second stage and
one proves the first part of the following
\begin{theorem}
\label{cycliccohomology}
  The Hochschild and cyclic cohomology of $\A^\hbar\rtimes G$ are given by 
\begin{displaymath}
\begin{split}
  HH^\bullet (\A^{((\hbar))}\rtimes G)& \cong  
  H^\bullet_{\text{\tiny \rm orb}} (X,\C((\hbar))),\\ 
  HC^\bullet (\A^{((\hbar))}\rtimes G)& \cong\bigoplus_{k\geq 0}
  H^{\bullet-2k}_{\text{\tiny \rm orb}}(X,\C((\hbar))).
\end{split}
\end{displaymath}
Furthermore, the pairing between homology and cohomology is given by 
Poincar\'{e} duality for orbifolds.
\end{theorem}
\begin{proof}
 Observe that the isomorphism in Hochschild cohomology described above is 
 induced by the maps dual to the ``reduction to loops'' from 
 Eq.~(\ref{redloops}):
\begin{equation}
\label{restrloops}
  \Gamma_c(B^{(k)},\sigma_{k+1}^{-1}(\A^{((\hbar))})^{\boxtimes(k+1)})'
  \rightarrow(\A^{((\hbar))}\rtimes G)_\natural^k,
\end{equation}
where $(~)'$ denotes the topological dual. Notice that the 
sheaf of distributions $\A'$ has the property that 
$$
  \Gamma(X,\A')=\Gamma_c(X,\A)'.
$$ 
Therefore, we can realize the left hand side of Eq.~(\ref{redloops})
as the space of sections over $B^{(k)}$ of the pullback via 
$\sigma_{k+1}$ of the sheaf of distributional Laurent-series, dual to
$\mathcal A^{((\hbar))}$. On a groupoid, the sheaf of distributions has 
a natural $G^{\text{\tiny op}}$-structure. This property relates, as in 
\cite{crainic} for the sheaf of smooth functions, the left hand side of 
(\ref{redloops}) to the twisted cyclic cohomology of 
$(\mathcal A^{((\hbar))})'$ on $\Lambda G$. Dualizing 
Prop.~\ref{quantcycquasiiso}, one obtains a quasi-isomorphism from
$\big(\mathcal A^{((\hbar))}\big)'_{\natural,\text{\rm tw}}$
to $\big( \mathcal A^{((\hbar))}_{(0)} \big)'_\natural$ 
with the trivial $\theta$-cyclic structure. 
To compute the hypercohomology of the resulting sheaf complex 
$\operatorname{Tot}^\bullet \mathcal B^{\bullet,\bullet} 
\big( \A^{((\hbar))}_{(0)} \big)$ dual to the 
$\operatorname{Tot}_\bullet \mathcal B_{\bullet,\bullet} 
\big( \A^{((\hbar))}_{(0)} \big)$-complex,
we use the following quasi-isomorphism which is dual to the one of 
Prop.~\ref{quasi-iso} and which is obtained by computing the cyclic 
cohomology of the formal Weyl algebra of $\C^n$ with compact support:
$$
  \operatorname{Tot}^\bullet 
  \big( \mathcal B^{\bullet,\bullet} \big( \A^{((\hbar))}_{(0)}\big) \big)
  \cong \bigoplus_{k\in \mathbb{N}}\underline{\C}((\hbar))[2k].
$$ 
Taking the hypercohomology of this complex over $\Lambda G$ gives the result.
\end{proof}
\subsection{Examples}
\label{mainexamples}
In this section we will give some examples of the computations of the 
Hochschild and cyclic homology, as well as the noncommutative Poisson 
homology. All are related to so-called transformation groupoids, i.e., 
the groupoid $G= \Gamma \ltimes M \rightrightarrows M$ associated to a
proper action of a discrete group $\Gamma$ on a manifold $M$. In this case 
the underlying orbifold is simply the quotient $X=M/\Gamma$. An invariant 
Poisson structure on $M$ leads, by equation \eqref{def:pi}, to a 
noncommutative Poisson structure on $\mathcal C^\infty_c(G)=A_c\rtimes\Gamma$, 
where $A_c=\mathcal C_c^\infty(M)$, already stated in \cite{xu}. As a special 
case of the quantization procedure for Poisson groupoids \cite{ta:lemma}, a  
$\Gamma$-invariant deformation quantization $A^\hbar_c$ of $A_c$ defines a 
quantization of the noncommutative Poisson algebra $A_c\rtimes\Gamma$ by 
taking the crossed product $A^\hbar_c\rtimes\Gamma$. 
We discuss our computations in the following special cases:

\begin{example}
\label{ex:crossed} (Free action) Notice that in case of a free action 
$\Gamma$ must be a finite group, since the groupoid $G$ is assumed to be 
proper. Trivially, we have that $B^{(0)}=M$ and $\Lambda G= \Gamma \ltimes M$,
i.e., $\tilde{X}=X$ in this case. Using the computation of the Hochschild 
homology of $A_c\rtimes\Gamma$ in \cite{bn,crainic}, and the Leray spectral 
sequence associated to the groupoid morphism 
$\pi: \Gamma \ltimes M\to M/\Gamma$ sending a point of $M$ to its image in the 
quotient space $M/\Gamma$, cf. \eqref{d}, one obtains an isomorphism 
$\pi: HH_\bullet( A_\text{\rm c} \rtimes G)\to HH_\bullet(M/\Gamma)$. It is 
easy to check that under $\pi$, the Poisson differential on
$HH_\bullet(A_\text{\rm c})$ is mapped to the Poisson differential on
$HH_\bullet(M/\Gamma)$, since a section of a sheaf on $M/\Gamma$ can be
identified with an invariant section of the corresponding
$G$-sheaf on $\Gamma \ltimes M$. Therefore, the Poisson homology
of $A_\text{\rm c}\rtimes \Gamma$ is equal to the Poisson
homology of the quotient space $M/\Gamma$. When the Poisson
structure is in fact symplectic, the Poisson homology on $M/\Gamma$ is dual to
the de Rham cohomology of $M/\Gamma$ by Brylinski's result in
\cite{byi:poisson}. 
\end{example}
\begin{remark}
\label{morita}
We can also use Proposition \ref{prop:morita} to obtain this result. 
We know that 
$A_\text{\rm c}\rtimes \Gamma$ is Morita equivalent to 
$\mathcal C_\text{\rm c}^{\infty}(M/\Gamma)$ by the bimodule 
$\mathcal C^\infty_\text{\rm c} (M)$. As shown in \cite{xu}, under this
Morita equivalence the Poisson structure $\Pi$ on
$A_\text{\rm c}\rtimes \Gamma$ is mapped to the Poisson structure
$\tilde{\Pi}$ on $M/\Gamma$ coming from the projection 
$\pi: M\to M/\Gamma$. By Prop.~\ref{prop:morita}, the Poisson homology
of $(A_\text{\rm c}\rtimes \Gamma, \Pi)$ is isomorphic to the Poisson
homology of $(\mathcal C_\text{\rm c}^{\infty}(M/\Gamma), \tilde{\Pi})$.
\end{remark}
Next, we assume the Poisson structure to be symplectic, so that $M$ is a 
manifold of dimension $2n$. For a free action, the projection 
$\pi:M\rightarrow M/\Gamma=X$ gives rise to a quasi-isomorphism 
$\pi_!\Omega_M\cong \Omega_X\in\mathsf{Sh}(X)$, where the right hand side 
is the sheaf of de Rham forms. Then the computations of the Hochschild and 
cyclic homology of the algebra $A_c^{((\hbar))}\rtimes\Gamma$ give: 
\begin{eqnarray*}
 HH_\bullet(A_c^{((\hbar))}\rtimes\Gamma)&\cong& H^{2n-\bullet}_c(X,\C)
 \otimes\C((\hbar)),\\
 HC_\bullet(A_c^{((\hbar))}\rtimes\Gamma)&\cong& 
 \bigoplus_{k\geq 0}H^{2n+2k-\bullet}(X,\C)\otimes\C((\hbar)).
\end{eqnarray*}
As for the computation of the Poisson homology, this could have been 
deduced at once by the observation that the quotient $A_c^\hbar/\Gamma$ 
gives a quantization of $\mathcal C_c^\infty(X)$ which is Morita equivalent
to the crossed product $A^\hbar_c\rtimes\Gamma$ by the equivalence bimodule 
$A^\hbar_c$. Therefore, the Hochschild and cyclic homology may be computed 
from the deformation quantization of $\mathcal C_c^\infty(X)$, for which the 
computations in \cite{nt} give the results above. 

As can be seen from above, for a trivial group our computations reduce to 
the well-known statements in \cite{nt}. We leave the statements about 
Hochschild and cyclic cohomology to the reader. Finally, notice that the 
space of traces on $A^\hbar_c\rtimes\Gamma$, i.e., 
$HC^0(A^{((\hbar))}_c\rtimes\Gamma)$ or the dual of 
$HH_0(A^{((\hbar))}_c\rtimes\Gamma)$, is one dimensional: all traces are 
proportional to each other. (See section \ref{sectraces}.)
\begin{remark}
It is easy to generalize this discussion to groupoids with trivial isotropy
groups, for example those obtained from a covering of a manifold. 
Again by Morita invariance, the computations of Poisson, Hochschild and 
cyclic homology reduce to the results of \cite{nt} on the underlying 
smooth manifold.
\end{remark}

\begin{example}
(Proper action, cf.~\cite{bc} for the undeformed case) 
One easily checks that for a proper action of a discrete 
group, the transformation groupoid $\Gamma \ltimes M \rightrightarrows M$ is
proper and \'{e}tale. Since we no longer assume the action to be free, 
$B^{(0)}$ is usually bigger than $M$, more precisely one has 
$B^{(0)}= \{ (\gamma , x)\in \Gamma \ltimes M \mid \gamma x = x\}$. 
Now, $B^{(0)}$ has the following decomposition into sectors 
$$
  B^{(0)}=\coprod_{\left<\gamma\right>\in \operatorname{Conj} (\Gamma)}
  \mathcal O_{\left<\gamma\right>},
$$ 
where $\left<\gamma\right>$ denotes the conjugacy class of $\gamma$ in 
$\Gamma$, and
$$
  \mathcal O_{\left<\gamma\right>}:= 
  \coprod_{\gamma'\in\left<\gamma\right>}M^{\gamma'}.
$$
According to Theorem \ref{thm:locality}, we therefore get the decomposition:
\begin{displaymath}
  H^{\Pi}_\bullet(A_\text{\rm c}\rtimes \Gamma)=
  \bigoplus_{\left<\gamma\right>\in \operatorname{Conj} (\Gamma)}
  H^{\Pi}_\bullet (A_\text{\rm c}\rtimes 
  \Gamma)_{\mathcal O_{\left<\gamma\right>}}.
\end{displaymath}
For $\gamma \in \Gamma$, we define 
$Z_\gamma := \{ \gamma' \in \Gamma \mid \gamma' \gamma = \gamma \gamma' \}$
and $N_\gamma = Z_\gamma / \langle\gamma\rangle $. We have Morita equivalences
$\Lambda G_{\mathcal O_{\left<\gamma\right>}}\simeq Z_\gamma \ltimes M^\gamma$
and $NG_{\mathcal O_{<\gamma>}} \simeq N_\gamma \ltimes M^\gamma $, 
which induce isomorphisms between the corresponding Poisson homologies 
by a similar argument as above for a free action. Therefore, we get
\begin{displaymath}
H^{\Pi}_\bullet (A_\text{\rm c}\rtimes \Gamma)=
\bigoplus_{\left<\gamma\right>\in \operatorname{Conj} (\Gamma)}
H^{\Pi}_\bullet (N_\gamma \ltimes M^\gamma ),
\end{displaymath}
and the right hand can be computed from the Poisson differential on the 
invariant differential forms on $M^\gamma/N_\gamma$. The decomposition of 
$B^{0}$ above into sectors, together with the Morita equivalences, give a 
decomposition of the inertia orbifold $\tilde{X}$ of $X=M/\Gamma$ as 
$$
  \tilde{X}=\coprod_{\left<\gamma\right>\in\operatorname{Conj}(\Gamma)}
  M^\gamma/N_\gamma.
$$ 
Therefore, we have for the Hochschild homology 
$$
  HH_\bullet(A^{((\hbar))}_\text{c}\rtimes\Gamma)\cong 
  \bigoplus_{\left<\gamma\right>\in \operatorname{Conj} (\Gamma)}
  H^{\operatorname{dim}(M^\gamma)-\bullet}_\text{c}
  (M^\gamma/N_\gamma,\C)\otimes\C((\hbar)).
$$
\end{example}
\section{Traces on the deformed groupoid algebra}
\label{sectraces}
Traces on an algebra obtained by deformation quantization form an important
ingredient in index theory. Since such functionals are nothing but cyclic 
cocycles of degree zero, Theorem \ref{cycliccohomology} gives a complete 
classification of traces on the deformed groupoid algebra 
$\A^\hbar\rtimes G$, that means of maps 
$\tr:\A^\hbar\rtimes G\rightarrow \C((\hbar))$ such that 
\begin{equation}
\label{traceproperty}
  \tr(a\star_\text{c} b)=\tr(b\star_\text{c} a),\quad 
  \text{ for all $a,b\in\A^\hbar\rtimes G$}.
\end{equation}
In this section, we will be concerned with the actual construction of 
all traces. Our discussion 
somewhat parallels with the constructions in \cite{fst:conjecture}, and 
also uses in an essential way the paper \cite{fe:g-index}, however notice 
that \cite{fst:conjecture} is only concerned with the subalgebra of 
$\A^\hbar\rtimes G$ of invariant quantized functions on the underlying 
orbifold $X$. The full algebra $\A^\hbar\rtimes G$ contains more information
which we believe to be essential for index theory.
\subsection{Traces on finite transformation groupoids}
We work in the situation of Sec.~\ref{mainexamples} and use the notation from 
there. Additionally, we assume (for notational convenience only) that each 
fixed point manifold $M^\gamma \subset M$ with $\gamma \in \Gamma$  
has constant dimension. We then consider the 
crossed product algebra 
$(A^{\hbar}_\text{c} \rtimes \Gamma,\star_\text{c})$, 
where the star product $\star$ 
on $A_\text{c}$ has been obtained by a $\Gamma$-invariant Fedosov 
construction. As explained above, the cyclic cohomology group  
$HC^0 \big (A^{((\hbar))}_\text{c} \rtimes \Gamma \big)$    
determines the space of traces on $A^{\hbar}_\text{c} \rtimes \Gamma$ and     
is given as follows:
\begin{equation}   
\label{tracecyccohom}  
  HC^0 \big( A^{((\hbar))}_\text{c} \rtimes \Gamma \big)  
  \cong  H^0_{\text{\tiny \rm orb}} (M/\Gamma,\C)\, ((\hbar))
  = \bigoplus_{\left<\gamma\right>\in \operatorname{Conj} (\Gamma)}
  H^0 (M^\gamma /N_\gamma, \C)\,((\hbar)).  
\end{equation}
Hence, the space of traces has dimension
\begin{displaymath}
  \dim_{\C((\hbar))} H^0_\text{\tiny orb} (M/\Gamma, \C)\,((\hbar)) = 
  \# \operatorname{Comp}  (\widetilde{M/\Gamma} ), 
\end{displaymath}
the number of connected components of the inertia orbifold.

Let us now examine the space of traces on 
$A^{\hbar}_\text{c} \rtimes \Gamma$ in some more detail.
To this end we will use in the remainder of this section the 
following notation.
Like in Sec.~\ref{hceg}, Eq.~(\ref{repelconvalg}), we expand elements 
$a \in A_\text{c} \rtimes \Gamma$ as sums 
$\sum_{\gamma \in \Gamma} f_\gamma \delta_\gamma$ with 
$f_\gamma \in A_\text{c}$ 
and extend this decomposition to formal Laurent series
$a = \sum_k a_k \hbar^k \in A_\text{c}^{((\hbar))} \rtimes \Gamma$
with $a_k \in A_\text{c} \rtimes \Gamma$ as follows:
\begin{equation}
  a =  \sum_{\gamma \in \Gamma} f_\gamma \delta_\gamma , \quad 
 \text{ where } f_\gamma  = \sum_k f_{k,\gamma}  \hbar^k
 \text{ and } 
 a_k = \sum_{\gamma \in \Gamma} f_{k,\gamma} \delta_\gamma .
\end{equation}
Following \cite[Sec.~1]{fe:g-index} and \cite{fst:conjecture} we
now consider a family $( \tau_\gamma )_{\gamma\in \Gamma}$ of
linear forms on $A_\text{c}^{((\hbar))}$ with the following properties:
\begin{eqnarray}
 \tau_\gamma (f  \star f') &=&
 \tau_{\gamma} ( f' \star  \gamma f ) \quad 
 \textrm{for all  $\gamma\in \Gamma$ and $f, f' \in A_\text{c}$},
 \label{TrAxAverInv}\\
 \tau_{\gamma'}(f ) &=& \tau_{\gamma\gamma'\gamma^{-1}}(\gamma f)
 \quad \textrm{for all $f\in A_\text{c}$ and 
 $\gamma,\gamma' \in \Gamma$}.
 \label{TrAxConjInv}
\end{eqnarray}
The following result can now be verified by a straightforward 
computation.
\begin{proposition}\label{CrossedProp}
Under the assumptions stated above let 
$(\tau_\gamma)_{\gamma \in \Gamma}$ be a family of linear
forms on $A_\text{\rm c}^{((\hbar))}$ which satisfies the assumptions
(\ref{TrAxAverInv}) and (\ref{TrAxConjInv}). Then
the functional 
\begin{equation}
\label{deftrace}
 \tr : A_\text{\rm c}^{((\hbar))} \rtimes \Gamma \to \C((\hbar)) , \quad
 a \mapsto \sum_{\gamma\in \Gamma} \tau_\gamma( f_\gamma)
\end{equation}
is a trace on $A_\text{\rm c}^{((\hbar))} \rtimes \Gamma$. 
Vice versa, given a trace 
$\tr :\: A_\text{\rm c}^{((\hbar))} \rtimes \Gamma \rightarrow \C((\hbar))$,
one obtains a family
$(\tau_\gamma)_{\gamma\in \Gamma}$ satisfying the above 
conditions by defining
\begin{equation}
 \tau_\gamma (f) := \tr ( f \delta_\gamma )  \quad 
  \text{ for all $f \in A_\text{\rm c}^{((\hbar))}$}.
\end{equation}
Finally, the trace corresponding to the so
defined family of linear forms coincides with the originally given
one.
\end{proposition}

In order to construct all traces on the quantized convolution algebra
we thus have to find functionals $(\tau_\gamma)_{\gamma\in\Gamma}$ 
satisfying the assumptions made above. In \cite{fe:g-index}
Fedosov has explicitly constructed such functionals. Let us recall
Fedosov's construction. To this end we restrict our assumptions further and
assume that $M$ is an open $\Gamma$-invariant convex neighborhood
of the origin of some symplectic vector space $V$. Then it is
well-known that over $M$, the star product $\star$ on $A^\hbar$ 
is equivalent to the Weyl star product $\star_\text{W}$ coming from $V$.
Let 
$S = 1 +\sum_{k=1}^\infty S_k \hbar^k : \mathcal C^\infty (M)[[\hbar]]
\rightarrow C^\infty (M)[[\hbar]]$ be an equivalence from 
$\star$ to the Weyl star product $\star_\text{W}$. Next choose a
$\Gamma$-invariant complex structure $J$ on $V$, and consider the 
Hermitian product induced by $J$ on the symplectic vector space $V$.
Since then, $\Gamma$ acts unitarily on $V$, one has for every 
$\gamma \in \Gamma$ a decomposition 
$V = V^\gamma \oplus V^\perp$, where $V^\gamma$ is the fixed point subspace
of $\gamma$, and $V^\perp$ its orthogonal complement. 
Now, $\gamma$ leaves $V^\perp$ invariant and acts on 
$V^\perp$ via a matrix $\gamma_\perp$.
Finally, choose complex unitary coordinates $z_\text{inv}$ of $V^\gamma$, 
$z_\perp$ of $V^\perp$, and put $z=(z_\text{inv},z_\perp)$.
With these notations, one can define functionals 
$\tau^\text{W}_\gamma : \mathcal C^\infty (M) ((\hbar)) 
 \rightarrow \C{((\hbar))}$
as follows 
(with integration induced by the real part of the Hermitian product):
\begin{equation}
\label{FedosovWgammatr}
  \tau^\text{W}_\gamma ( f) = \frac{1}{(2\pi \hbar)^k}
  \int_{V^\gamma} \frac{1}{\det (1 - \gamma_\perp^{-1})} \left.\exp
  \left(\hbar \frac{\partial}{\partial z_\perp}
  \frac{1 + \gamma_\perp^{-1}}{1 -\gamma_\perp^{-1}}
  \frac{\partial}{\partial z_\perp^*}\right) f
  \right|_{z_\perp = \overline{z}_\perp =0} 
  dz_\text{inv} d\overline{z}_\text{inv} .
\end{equation}
According to Eq.~(2.17) and Prop.~2.5 of \cite{fe:g-index}, the thus 
obtained family 
$(\tau_\gamma^\text{W})_{\gamma \in \Gamma}$ satisfies properties
(\ref{TrAxAverInv}) and (\ref{TrAxConjInv}) above with respect to the Weyl 
star product. The following result is essentially a reformulation 
of \cite[Prop.~2.5]{fe:g-index} and \cite[Cor.~7.5]{fst:conjecture}.
\begin{proposition}
  Let $\Gamma$, $V$, $M\subset V$ and $(A^\hbar,\star)$ as before, and
  choose $S$, $J$ and $z=(z_\text{\rm inv},z_\perp)$ like above. 
  Then the family $(\tau_\gamma^\text{F})_{\gamma \in \Gamma}$ 
  of functionals on $A_\text{\rm c}^{((\hbar))}$ defined  by
  \begin{equation}
  \label{FedosovFgammatr}
     \tau^\text{F}_\gamma (f) := \tau^\text{W}_\gamma (S f) \quad
     \text{ for all $f \in  A_\text{\rm c}$}
  \end{equation}
  satisfies the conditions (\ref{TrAxAverInv}) and (\ref{TrAxConjInv}) 
  above. Moreover, for every $\gamma$, the restriction of the functional
  $\tau^\text{F}_\gamma$ to the $\Gamma$-invariant elements 
  of $A_\text{\rm c}^{((\hbar))}$ does not depend on the choice 
  of $J$, $z$ and $S$. 
\end{proposition}
Now consider a family 
$\kappa = 
 (\kappa_{\left< \gamma\right>})_{\left<\gamma\right> \in \operatorname{Conj} 
 (\Gamma)}$  
of complex coefficients. Then the functional
\begin{equation}
\label{deftraces}
  \tr_\kappa :\:   A_\text{\rm c}^{((\hbar))} \rtimes \Gamma \rightarrow
  \C((\hbar)), \quad
  a =\sum_{\gamma \in \Gamma} f_\gamma \delta_\gamma \mapsto
  \sum_{\gamma \in \Gamma} \, 
  \kappa_{\left<\gamma\right>} \,  \tau^\text{F}_{\gamma} (f_{\gamma})
\end{equation}
has to be a trace on the deformed crossed product algebra by the preceding
propositions, since the family 
$(\kappa_{\left< \gamma \right>} \tau^\text{F}_\gamma )_{\gamma \in \Gamma}$
satisfies conditions (\ref{TrAxAverInv}) and (\ref{TrAxConjInv}).
One  even has more.
\begin{corollary}
\label{trfintransgrp}
With notations from above, the functionals $\tr_\kappa$ 
have the following properties.
\begin{enumerate}
\item[(1)]
Every trace on $A_\text{\rm c}^{((\hbar))}  \rtimes \Gamma$ is of the form
$\tr_\kappa$ with a uniquely determined family 
$\kappa \in \C^{\operatorname{Conj} (\Gamma)}$.
\item[(2)]
The traces $\tr_\kappa$ are invariant in the following sense.
Let $V'$ be another symplectic vector space, 
$\Gamma'$ a finite group acting by linear symplectomorphisms on $V'$, and 
let $\star'$ be a Fedosov star product on $V'$. 
Assume further that $F : M \hookrightarrow V' $
is an open embedding with the following properties:
\begin{enumerate}
\item[(a)] 
  $F$ is equivariant with respect to an injective 
  homomorphism $\iota : \Gamma \rightarrow \Gamma'$,
\item[(b)]
  $F$ is symplectic,
\item[(c)]
  the pull-back via $F$ induces a homomorphism of star product algebras
  $F^* : \big( \mathcal C^\infty (V')[[\hbar]] , \star' \big)
  \rightarrow A^\hbar$,
\item[(d)]
  the induced quotient map $\overline{F} : M/\Gamma \rightarrow V'/\Gamma'$
  is an open embedding. 
\end{enumerate}
Then  for every 
$a = \sum_{\gamma \in \Gamma} f_\gamma \delta_\gamma \in 
A_\text{\rm c}^{((\hbar))}  \rtimes \Gamma$ the equality
\begin{equation}
\label{scdclaim}
   \tr_{\iota_* \kappa} ( F_* a)  = \tr_\kappa (a)
\end{equation}
holds true, where $F_* (a) = \sum_{\gamma \in \Gamma} 
 (f_\gamma \circ F^{-1}) \, \delta_{\iota (\gamma)}$ and
 $\iota_* \kappa $ is the family 
 $\big( \kappa'_{\left< \gamma' \right>} 
  \big)_{\left< \gamma' \right> \in \operatorname{Conj} (\Gamma')}$ defined by 
 \begin{equation}
 \label{coeff}
   \kappa'_{\left< \gamma' \right>} =
   \begin{cases}
     \kappa_{\left< \gamma \right>}, &  
     \text{if $\gamma' \sim \iota (\gamma) $ for some $\gamma \in \Gamma$},\\ 
     0 , & \text{else}. 
   \end{cases} 
 \end{equation}
\end{enumerate}
\end{corollary}
\begin{proof}
 Observe that by Eq.~(\ref{deftraces}) the map 
 $\C^{\operatorname{Conj}(\Gamma)} \rightarrow 
 (A_\text{\rm c}^{((\hbar))}  \rtimes \Gamma)'$,
 $ \kappa \mapsto \tr_\kappa $ is
 injective, since each $\tau^\text{F}_\gamma$ is non-zero. 
 Since the dimension of the space of traces is 
 $\# \operatorname{Conj} (\Gamma)$, the first claim follows.
 
 For the second claim note first that $\iota_* \kappa $ is well-defined 
 indeed, since by the assumption (d)
 on the equivariant embedding $(F,\iota)$, the induced map 
 $\overline{\iota} :\operatorname{Conj} (\Gamma) \rightarrow 
 \operatorname{Conj} (\Gamma')$ has to be injective.
 Next we conclude from \cite[Cor.~7.5]{fst:conjecture} that  
 $$
     \tau^\text{F}_\gamma (f) =  
     \tau^\text{F}_{\iota(\gamma)} (f\circ F^{-1} ) 
 $$
 holds for all $\Gamma$-invariant $f \in A_\text{\rm c}^{((\hbar))} $
 and all $\gamma \in \Gamma$. 
 By the first claim and the definition of $\iota_* \kappa$, 
 this entails Eq.~(\ref{scdclaim}).
\end{proof}

\subsection{Traces in the general case}
Let now $G$ be an arbitrary proper \'etale Lie groupoid with a 
symplectic structure and let $\star$ be an invariant star product
on $\mathcal A$. We now want to construct all traces
on the crossed product $\mathcal A^{((\hbar))} \rtimes G$.
To this end first fix a dense countable family $(x_i)_{i\in I}$ of points
of $G_0$ and an open covering 
$\mathcal U = ( U_i )_{i\in I}$ of $G_0$
such that $x_i \in U_i$ for all $i\in I$ and such that one has isomorphisms 
$G_{|U_i}\cong \Gamma_i \ltimes U_i$, where $\Gamma_i$ is the isotropy group 
$G_{x_i}$. 
By appropriate choices we can even achieve, that each $U_i$ 
is symplectomorphic to an open ball around the origin of some symplectic 
$\Gamma_i$-representation space $V_i$, that $x_i$ corresponds to the origin 
under this symplectomorphism, and that for every pair $U_i,U_j$ 
with $\pi( U_i )  \cap \pi (U_j )\neq \emptyset$ there exists an open 
connected subset $W_{ij} \subset G_0$ and a finite isotropy 
group $\Gamma_{ij} := G_{x_{ij}}$ (with $x_{ij}\in W_{ij}$)
acting symplectically
on $W_{ij}$ such that $G_{|W_{ij}}\cong \Gamma_{ij} \ltimes W_{ij}$
and such that $\pi (U_i ), \pi (U_j ) \subset \pi (W_{ij})$. 
Moreover, we can assume that $W_{ij} = W_{ji}$ and that $W_{ij}$
is an open invariant set of some symplectic $\Gamma_{ij}$-representation 
space  $V_{ij}$. Finally, we can choose the sets $U_i$ so small such that 
for all $i,j\in I$ with $\pi( U_i )  \cap \pi (U_j )\neq \emptyset$
there exist bisections $s_{ij}: U_i \rightarrow W_{ij}$ with
$t\circ s_{ij} (x) \in W_{ij}$ for all $x\in U_i$. Now put 
$F_{ij} := t\circ s_{ij}$. Using that the underlying groupoid is proper 
\'etale and the assumptions on the covering $( U_i )_{i\in I}$ one then
immediately checks the following properties of the maps $F_{ij}$:
\begin{enumerate}
\item[(a)] 
  each $F_{ij}$ is an embedding and equivariant with respect to the  
  monomorphism $\Gamma_i \rightarrow \Gamma_{ij}$ induced by the
  composition on $G_1$,
\item[(b)]
  each $F_{ij}$ is symplectic,
\item[(c)]
  since the star product on $G_0$ is $G$-invariant, 
  the pull-back via $F_{ij}$ induces a homomorphism 
  $F_{ij}^* : \big( \mathcal C^\infty (W_{ij})[[\hbar]] , \star \big)
  \rightarrow \big( \mathcal C^\infty (U_i)[[\hbar]] , \star \big) $,
\item[(d)]
  the induced quotient map 
  $\overline{F}_{ij} : U_i/\Gamma_i \rightarrow W_{ij}/\Gamma_{ij}$
  is the natural inclusion of open subsets of $X$. 
\end{enumerate}
Using Cor.~\ref{trfintransgrp} (2), these properties will later guarantee that 
one can glue together local traces on $ \mathcal A^{((\hbar))} \rtimes G$.

Associated to the covering $\mathcal{U}$ is the groupoid 
$G_\mathcal{U}$ with objects and morphisms given by
\begin{equation}
\label{covgpd}
  (G_\mathcal{U})_0=\coprod_{i\in I}U_i,
  \quad (G_\mathcal{U})_1=\coprod_{i,j\in I} s^{-1}(U_i)\cap t^{-1}(U_j).
\end{equation}
The obvious morphism $G_\mathcal{U}\rightarrow G$ then is a weak
equivalence. The sheaf of quantum algebras $\A^{((\hbar))}$ restricts 
over every open subset $U_i$ to define a sheaf, also
denoted by $\A^{((\hbar))}$ and therefore defines a crossed product
$\A^{((\hbar))}\rtimes G_\mathcal{U}$.  Let us write elements in
$\A^{((\hbar))}\rtimes G_\mathcal{U}$ as $a=(a_{ij})_{i,j\in I}$, etc.
Denote by $\star_\mathcal{U}$ the multiplication on 
$\A^{((\hbar))}\rtimes G_\mathcal{U}$ 
obtained by combining the $\star$-product with the convolution
product. This product then reads as follows
(using germs in the notation):
\begin{equation}
\label{multiplication}
 [(a\star_\mathcal{U} b)_{ij}]_g=
 \sum_{k}\sum_{\substack{ {g_1g_2=g} \\ {s(g_1)=t(g_2)\in U_k}}}
 [a_{ik}]_{g_1} g_2 \star [b_{kj}]_{g_2},\quad s(g)\in U_i,~t(g)\in U_j.
\end{equation}
We will now construct an injective homomorphism
$\A^{((\hbar))}\rtimes G \rightarrow \A^{((\hbar))}\rtimes G_\mathcal{U}$.
Consider a partition of unity $(\varphi_i)_{i\in I}$ subordinate to
$\mathcal{U}$, satisfying $\sum_i\varphi_i^2=1$. Define the following
formal power series on $G_0$:
$$
  \Phi_i=\left(\sum_k\varphi_k\star\varphi_k\right)^{-1/2}\star
  \varphi_i,~\hspace{1cm}\Psi_i=\varphi_i\star
  \left(\sum_k\varphi_k\star\varphi_k\right)^{-1/2}.
$$ 
Notice that the inverse of the square root exists, since 
$\sum_k\varphi_k\star\varphi_k -  1 \in \hbar \mathcal A^\hbar (G_0)$. 
By construction, we have $\operatorname{supp}(\Phi_i)\subset U_i$, 
$\operatorname{supp} (\Psi_i)\subset U_i$, and
$\sum_i\Psi_i\star\Phi_i=1$. 
From these properties it is easy to deduce that for the 
``convolution $\star$-product'' in $\A^{((\hbar))}\rtimes G$ we have
$$
  \sum_i \left( \Psi_i \star_c \Phi_i\right)= \delta_\text{u},
$$ 
where the $\Psi_i$ and $\Phi_i$ have been extended by $0$ outside 
$G_0$, and where $\delta_\text{u}$ is the ``unit in the convolution 
algebra'' from Step 2., Sec.~\ref{hceg}. 
By inspection of the multiplication \eqref{multiplication} 
it follows that the map 
$$
 a\mapsto 
 \left( \Phi_i \star_{c} a \star_c \Psi_j \right)_{ij}
$$ 
defines a homomorphism $\Phi_\mathcal{U} : \: \A^{((\hbar))}\rtimes
G\rightarrow\A^{((\hbar))}\rtimes G_\mathcal{U}$.
With this notation, we have the following final result. 
\begin{theorem}
  Let $G, \mathcal U, \mathcal A^\hbar$ be like above and 
  $\kappa :B^{(0)} \rightarrow \C$ a locally constant $G$-invariant 
  function.
  Then the restriction of $\kappa$ to $B^{(0)} \cap s^{-1} (U_i)$ induces 
  for every $i\in I$ a family 
  $\kappa_i \in \C^{\operatorname{Conj} (\Gamma_i)}$. Moreover,
  the formula
  \begin{equation}
    \tr_\kappa (a) = \sum_i \tr_{\kappa_i} 
    \big( \Phi_i \star_\text{c} a \star_\text{c} \Psi_i \big)
  \end{equation}
  defines a trace on $\mathcal A^{((\hbar))} \rtimes G$, and
  every trace on $\mathcal A^{((\hbar))} \rtimes G$ is equal 
  to such a $\tr_\kappa$ with unique $\kappa$.
\end{theorem}
\begin{proof}
  The $\kappa_i$ induce traces $\tr_{\kappa_i}$ on 
  $A_{\text{c},i}^{((\hbar))}\rtimes \Gamma_i $,
  $A_{\text{c},i} := \mathcal C^\infty_\text{c} (U_i)$ by 
  Cor.~\ref{trfintransgrp}.
  Since $\Phi_\mathcal{U}$ is a homomorphism of algebras, 
  $\tr_\kappa$ is proved to be a trace, if the functional
  \begin{equation}
  \label{trfct}
    \tr_\mathcal{U} :\: \A^{((\hbar))}\rtimes G_\mathcal{U}
    \rightarrow \C((\hbar)),\quad
    (a_{ij}) \mapsto \sum_i \, \tr_{\kappa_i} (a_{ii}) 
  \end{equation}
  is a trace. To this end it suffices to show that for all 
  $(a_{ij}), (b_{ij}) \in  \A^{((\hbar))} \rtimes G_\mathcal{U}$ one has
  \begin{equation}
  \label{comrel}
    \tr_{\kappa_i} (a_{ij} \star_\text{c} b_{ji}) 
    = \tr_{\kappa_j} (b_{ji} \star_\text{c} a_{ij}) ,
  \end{equation}
  if $\pi (U_i) \cap \pi (U_j) \neq \emptyset$.  
  For the proof of this equality we use the equivariant embeddings
  $F_{ij}: U_i \rightarrow W_{ij}$ constructed above and apply 
  Cor.~\ref{trfintransgrp}. 
  More precisely, let $\tr_{ij}$ be the trace on 
  $A_{\text{c},ij}^{((\hbar))} \rtimes \Gamma_{ij}$,
  $A_{\text{c},ij} := \mathcal C^\infty_\text{c} (W_{ij})$ 
  induced by the restriction of $\kappa$ to 
  $B^{(0)} \cap s^{-1} (W_{ij})$. Cor.~\ref{trfintransgrp} (2) 
  entails that the left hand side of Eq.~(\ref{comrel}) coincides
  with $\tr_{ij} ({F_{ij}}_*a_{ij} \star_\text{c} {F_{ij}}_*b_{ji} )$, 
  and the right hand side with 
  $\tr_{ij} ( {F_{ij}}_* b_{ji} \star_\text{c} {F_{ij}}_* a_{ij} )$.
  By the trace property of $\tr_{ij}$,  Eq.~(\ref{comrel}) follows,
  and $\tr_\kappa$ is a trace indeed. 
    
  Since the map $\kappa \mapsto \tr_\kappa$ clearly is injective,
  the second claim now follows easily from the fact that
  $HC^0 (\A^{((\hbar))}\rtimes G ,\C )$ has dimension (over $\C((\hbar))$)
  equal to the number of components of the inertia orbifold,
  and the fact that the latter number gives also the complex dimension 
  of the space of locally constant invariant functions from $B^{(0)}$ to $\C$.
\end{proof}
\subsection{On a conjecture of Fedosov, Schulze and Tarkhanov}
Unlike in the case of a (connected) symplectic manifold, where the space of
traces on a deformed algebra of compactly supported smooth functions is
one dimensional, the space of traces on $\mathcal A^{((\hbar))}\rtimes G$
has dimension $>1$, since by Theorem \ref{cycliccohomology} this dimension 
is given by the number of connected components of the inertia orbifold 
$\widetilde X$. In \cite{fst:conjecture}, Fedosov, Schulze and Tarkhanov 
show that a certain abelian group of isomorphism classes of line bundles 
on a symplectic orbifold acts nontrivially on the space of traces of the 
deformed convolution algebra and conjecture that 
``this ambiguity in traces is the only possible one''. 
In our framework such type of questions can be answered naturally. 

We start with a different view on orbifold cohomology. Consider the 
representation ring sheaf $\underline{R}_\C$ on $X$ whose stalk at $x\in X$ 
is given by the complexified representation ring 
$R_\C(G_x) = R (G_x) \otimes_\Z \C$ of the isotropy group
$G_x$, a finite group. As explained in \cite[Sec.~6.4]{moerdijk},
the Leray spectral sequence associated to the morphism of groupoids 
$\beta:\Lambda G \rightarrow G$ yields an isomorphism 
\begin{equation}
\label{repring}
  H^k_{\rm\tiny orb}(X,\C)\cong H^k(X,\underline{R}_\C),
\end{equation}
where the right hand side is simply sheaf cohomology on the space $X$. 
Consider now the abelian group $P_G := H^1(G,S^1) / H^1(X,S^1)$. 
By Sec.~\ref{sheaves}, $H^1(G,S^1)$ classifies the $G$-line bundles on 
$G$, whereas $H^1(X,S^1)$ gives the set of isomorphism classes of line bundles 
on $X$, which, by pull-back along the projection $\pi:G\rightarrow X$, 
identifies with the set of isomorphism classes of $G$-line 
bundles with trivial action of the isotropy groups. Thus, $P_G$
is the Picard group as defined in \cite{fst:conjecture},
however we do not use this terminology here in view of the Picard group in 
Poisson geometry (cf.~\cite{bw:picard}) which is a completely different group.
Now, there is a natural homomorphism from $P_G$ into
the group of units of the ring $H^0(X,\underline{R}_\C)$.
For its construction observe that every $G$-line bundle gives rise to a 
representation of the isotropy $G_x$ for every $x\in X$,
hence, by taking the character at every point, there is a canonical map 
$H^1(G,S^1) \rightarrow H^0(X,\underline{R}_\C)$.
As its kernel is given by $H^1(X,S^1)$, the existence
of the injection $P_G \rightarrow H^0(X,\underline{R}_\C)$ follows.

Since traces on an algebra are nothing but cyclic $0$-cocycles, 
Theorem \ref{cycliccohomology} entails that the space of traces on
$\mathcal A^{((\hbar))}\rtimes G$ is isomorphic 
to $H^0_{\rm\tiny orb}(X,\C)\otimes\C((\hbar))$. By Eq.~(\ref{repring}),
the conjecture of \cite{fst:conjecture} can now be reformulated as 
the statement that the image of $P_G$ in $H^0(X,\underline{R}_\C)$ forms a 
basis. But since irreducible representations of a finite group are necessarily 
one dimensional only if the underlying group is abelian,
the claim holds true in general, if and only if every isotropy
group $G_x$ is abelian.
Therefore, the conjecture in \cite{fst:conjecture} is true for $G$ with
abelian isotropy groups, but not otherwise.

Strictly speaking, the paper \cite{fst:conjecture} is only concerned with 
the algebra $\Gamma_{\text{\rm\tiny inv,c}}(\A^\hbar)$ of invariant sections 
of $\A^\hbar$, or in other words with the deformation quantization of 
$\mathcal{C}^\infty_c(X)$ constructed in \cite{P2}, 
cf.~Sec.~\ref{sec:intr-ncq}. However, the conclusion 
above remains true also in this case in view of the following:
\begin{proposition}
 In case the proper \'etale Lie groupoid $G$ is reduced, i.e.~if each
 isotropy group $G_x$ acts faithfully on a neighborhood of $x\in G_0$
 {\rm (}cf.~\cite[Sec.~1.5]{moerdijk}{\rm )}, then the 
 algebras $\A^{((\hbar))}\rtimes G$ and 
 $\Gamma_{\text{\rm\tiny inv,c}}(\A^{((\hbar))})$ are Morita equivalent.
\end{proposition}
\begin{sketchproof}
The equivalence bimodule is given by $\A^{((\hbar))}(G_0)$, the quantization 
of the symplectic manifold $G_0$. To prove that this really defines a Morita
equivalence one first observes that it suffices to prove the claim locally.
One can check this for example by using the groupoid
$G_\mathcal{U}$ and the partition of unity associated to a covering of
$G_0$. Using a covering by open
subsets over which the restricted groupoid is isomorphic to a translation
groupoid by a faithful action of a finite group, 
the claim is proved locally as in \cite{dolget} 
by the fact that the deformed algebra $\A^{((\hbar))}(M)$ of a symplectic 
manifold is simple. The latter holds true since there are no nontrivial 
Poisson ideals in the ring of smooth functions on a  symplectic manifold.
\end{sketchproof}
\newpage 
\bibliographystyle{alpha}

\end{document}